\documentclass[10pt]{article}
\newcommand{\version}{\today}
\newcommand{\lanbox}{\hfill \hbox{$\, \vrule height 0.25cm width 0.25cm depth 0.01cm
\,$}}
\usepackage{amsthm,amsfonts,amsmath, amscd}

\swapnumbers
\theoremstyle{plain}
\newtheorem{thm}{THEOREM}[section]
\newtheorem{lm}[thm]{LEMMA}
\newtheorem{cl}[thm]{COROLLARY}
\newtheorem{prop}[thm]{PROPOSITION}
\newtheorem{remark}[thm]{Remark}

\theoremstyle{definition}
\newtheorem{defi}[thm]{DEFINITION}
\theoremstyle{remark}
\newcommand{\upchi}{\raise1pt\hbox{$\chi$}}
\newcommand{\R}{{\mathord{\mathbb R}}}

\newcommand{\Z}{{\mathord{\mathbb Z}}}
\newcommand{\N}{{\mathord{\mathbb N}}}

\newcommand{\hn}{{\mathord{\widehat{n}}}}

\renewcommand{\AA}{{A}}

\renewcommand{\|}{{\Vert}}
\numberwithin{equation}{section}
\newcommand{\Fhls}{{\mathcal F}_{\rm PKS}}

\newcommand{\Fcfd}{{\mathcal H}_{\lambda}}
\newcommand{\Fcfdr}{{\mathcal H}_{\lambda,\delta}}

\newcommand{\RR}{{\mathbb R}}
\newcommand{\NN}{{\mathbb N}}

\newcommand{\PP}{{\cal P}}

\newcommand{\W}{{\rm W}}
\newcommand{\dd}{\:{\rm d}}
\newcommand{\dx}{\;{\rm d}x}

\newcommand{\ds}{\;{\rm d}s}

\renewcommand{\ln}{\log}
\newcommand{\eps}{\varepsilon}
\newcommand{\un}{{\rm 1\kern -2.5pt l}}
\newcommand{\id}{{\rm Id}}
\newcommand{\rhohls}{{\varrho_{\lambda}}}

\newcommand{\crho}{C_{\rho_0}}

\begin{document}
\markboth{\scriptsize{BCC \version}}{\scriptsize{BCC \version}}
\def\mn{{\bf M}_n}
\def\hn{{\bf H}_n}
\def\hnp{{\bf H}_n^+}
\def\hmnp{{\bf H}_{mn}^+}
\def\h{{\cal H}}
\title{{\sc Functional inequalities, thick tails and asymptotics \\
for the critical mass Patlak-Keller-Segel model}}
\author{
\vspace{5pt} Adrien Blanchet$^{1}$, Eric A. Carlen$^{2}$  and Jos\'e A. Carrillo$^{3}$ \\
\vspace{5pt}\small{$^{1}$ TSE (GREMAQ),}\\[-6pt]
\vspace{5pt}\small{Universit\'e de Toulouse, 21 All\'ee de Brienne, F-31000 Toulouse, France}\\[-6pt]
\vspace{5pt}\small{$^{2}$ Department of Mathematics, Hill Center,}\\[-6pt]
\small{Rutgers University,
110 Frelinghuysen Road
Piscataway NJ 08854-8019 USA}\\
\vspace{5pt}\small{$^{3}$ Instituci\'o Catalana de Recerca i Estudis Avan\c cats and Departament de Matem\`atiques}\\[-6pt]
\small{Universitat Aut\`onoma de Barcelona, E-08193 Bellaterra, Spain}\\[-6pt]
 }
\date{\today}
\maketitle
\footnotetext [1]{Work partially supported by the project EVaMEF ANR-09-JCJC-0096-01.}

\footnotetext [2]{Work partially supported by U.S.
National Science Foundation grant  DMS 0901632.}

\footnotetext [3]{Work partially supported by the projects
MTM2011-27739-C04-02. DGI-MCI (Spain) and
2009-SGR-345 from AGAUR-Generalitat de Catalunya. \\
\copyright\, 2010 by the authors. This paper may be reproduced, in
its entirety, for non-commercial purposes.}

\begin{abstract}
We investigate the long time behavior of the critical mass
Patlak-Keller-Segel equation. This equation has a one parameter
family of steady-state solutions $\rhohls$,
$\lambda>0$, with thick tails whose second moment is not bounded.
We show that these steady state solutions are stable, and find
basins of attraction for them using an entropy functional $\Fcfd$
coming from the critical fast diffusion equation in $\R^2$. We
construct solutions of Patlak-Keller-Segel equation satisfying an
entropy-entropy dissipation inequality for $\Fcfd$. While the
entropy dissipation for $\Fcfd$ is strictly positive, it turns out
to be a difference of two terms, neither of which need to be small
when the dissipation is small. We introduce a strategy of {\em
controlled concentration} to deal with this issue, and then use
the regularity obtained from the entropy-entropy dissipation
inequality to prove the existence of basins of attraction for each
stationary state composed by certain initial data converging
towards $\rhohls$. In the present paper, we do not
provide any estimate of the rate of convergence, but we discuss
how this would result from a stability result for a certain sharp
Gagliardo-Nirenberg-Sobolev inequality.
\end{abstract}

\medskip
\centerline{Mathematics subject classification numbers:  15A45, 49M20}

\section{Introduction}

\subsection{The PKS  system and its critical mass}

The Patlak-Keller-Segel system~\cite{patlak, Keller-Segel-70} is
one of the simplest models of {\it chemotaxis}, describing the
evolution of the population density of a cell colony which is
diffusing across a two dimensional surface. In addition to the
diffusion, as the cells move across the surface, they continually
emit a chemical attractant, which itself diffuses across the
surface. The cells tend to move towards higher concentrations of
the attractant, and this induces a drift term tending to
concentrate the population, and countering the spreading effects
of the diffusion. A model organism for this type of behavior is
the {\it Dictyostelium Discoideum} which segregates {\it cyclic
adenosine monophosphate}, another important example of chemotactic
movement are endothelial cells who react to VEGF to form blood
vessels. See  \cite{HP,Perthame} for recent reviews on chemotaxis
models.

The model is mathematically interesting on account of this
competition between the concentrating effects of the drift induced
by the chemical attractant and the spreading effects of the
diffusion, and there is a {\em critical value} to the total mass
of the initial data, so that for masses above this value, the
concentration wins, and the density collapses in a finite time.
However, for masses below this critical mass, diffusion dominates,
and the colony smoothly diffuses off to infinity. At the critical
mass, there is a continuous family of stationary solutions, and
this paper is concerned with determining their stability
properties, and since they all turn out to be stable, basins of
attraction for each of them. We begin  by introducing the model
and the critical mass associated with it. If $\rho$ denotes the
population density, and $c$ the concentration of the chemical
attractant,  the system of equations is
\begin{equation}\label{eq:sp}
\begin{cases}
{\displaystyle \frac{\partial \rho}{\partial t}(t,x) ={\rm div } \left[ \nabla \rho(t,x)-\,\rho(t,x)\nabla c(t,x)\right]}\qquad & t>0\,,\;x\in\R^2\;,\vspace{.2cm}\\
{\displaystyle c(t,x)=-\frac{1}{2\,\pi}\int_{\R^2}\log|x-y| \rho(t,y) \dd y}\;,\qquad &t>0\,,\;x\in\R^2\;,\vspace{.2cm}\\
\rho(0,x)=\rho_0(x)\geq 0\qquad &x\in\R^2\, ,
 \end{cases}
\end{equation}
with an appropriate choices of units, so that all dimensional
constants are unity.

In most of this paper, we consider initial data $\rho_0$ that
belongs to $L^1(\R^2,\log(e+|x|^2)\dd x)$, and such that
$\rho_0\log \rho_0$ is integrable. The relevance of these conditions
shall be explained shortly, but at the very least, they insure
that $c(0,x)$ is well defined. It will sometimes be convenient to write the second
equation in~\eqref{eq:sp} in the compact form $c(t,x) =
G*\rho(t,x)$ where $G(x) = -{1}/{(2\,\pi)}\log|x|$ is the Green's
function for $-\Delta$ in $\R^2$. That is, $-\Delta c=\rho$.

Also throughout the paper, the term {\em density} shall always
refer to a non-negative integrable function on $\R^2$, and we
shall use the term  {\em  mass} to refer to the total integral of
a density $\rho$. Because of the divergence form structure of the
system, solutions formally satisfy the  conservation of  mass
\begin{equation*}
\int_{\R^2}\rho(t,x)\dd x = \int_{\R^2}\rho_0(x)\dd x:=M
\end{equation*}
for all $t\ge 0$; {\it i.e.}, the mass $M$ is conserved in time.

The PKS system can be rewritten advantageously as follows:
Introduce the {\it free energy functional} $\Fhls$
\begin{equation*}
\Fhls[\rho]  = \int_{\R^2} \rho(x) \log \rho(x)\dd x  +
\frac{1}{4\pi} \iint_{\R^2 \times \R^2} \rho(x) \log |x-y| \rho(y)\dd x \dd y \ .
\end{equation*}
The first integral  is  well defined if $\rho\log \rho$
is integrable, and the positive part of $\rho(x)\log|x-y|\rho(y)$ is
integrable when  $\rho$ belongs to  $L^1(\R^2,\log(e+|x|^2)\dd
x)$, so that the second integral is at least well-defined under
this condition.

Now suppose that the density $\rho$ belongs to
$L^1(\R^2,\log(e+|x|^2)\dd x)$, and moreover, $\rho \log \rho$ is
integrable. Then a simple formal calculation shows that for  all
$u \in C_c^\infty(\R^2)$ with zero mean,
$$
\lim_{\epsilon\to 0}\frac{1}{\epsilon}\left(\Fhls[\rho+\epsilon u]
-\Fhls[\rho]\right) = \int_{\R^2}\frac{\delta
{\Fhls}[\rho]}{\delta  \rho}(x) \, u(x)\dd x
$$
where
\begin{equation*}
\frac{\delta {\Fhls}[\rho]}{\delta  \rho}(x) := \log \rho(x) +
\frac{1}{2\pi} \int_{\R^2}\log |x-y|\rho(y)\dd y = \log \rho(x)
-G*\rho(x) \ .
\end{equation*}
It is then easy to see that the evolution equation in~\eqref{eq:sp} can be rewritten as
\begin{equation}\label{wasgrad}
\frac{\partial \rho}{\partial t}(t,x)  = {\rm
div}\left(\rho(t,x)\nabla \left[ \frac{\delta
{\Fhls}[\rho(t)]}{\delta  \rho}(x) \right]\right)\ .
\end{equation}
It follows that at least along well-behaved classical solutions (for which we may integrate by parts),
\begin{equation}\label{enerdis}
\frac{{\rm d}}{{\rm d}t}{\Fhls}[\rho(t)] = - \int_{\R^2}
\rho(t,x)\left|\nabla  \frac{\delta {\Fhls}[\rho(t)]}{\delta
\rho}(x)\right|^2\dd x \ .
\end{equation}
In particular, along such solutions, $t \mapsto {\Fhls}[\rho(t)]$
is monotone non--increasing. The key to exploiting this monotonicity,
as discovered in~\cite{DP}, is the sharp logarithmic Hardy--Littlewood--Sobolev  (Log HLS) inequality
\cite{Be,CarLos}:
\begin{lm}[Logarithmic Hardy-Littlewood-Sobolev inequality]
\label{lm:loghls} Let $f$ be a non-negative measurable function in
$\R^2$ such that $f\ln f$ and $f\ln (e+|x|^2)$ belong  to
$L^1(\R^2)$. Let  $M:=\int_{\R^2}f\dd x$. Then
\begin{equation}\label{loghls}
\int_{\R^2}f(x)\log f(x)\dd x+\frac{2}{M}
\iint_{\R^2\times\R^2}f(x) \log|x-y| f(y)\dd x
\dd y\ge -\;C(M)\,,
\end{equation}
with $C(M):=M\,(1+\log\pi-\log(M))$. There is equality if and only if
$f(x) = \rhohls(x-x_0)$ for  some $\lambda > 0$ and
some $x_0\in \R^2$, where
\begin{equation*}
 \rhohls(x):=\frac{M}{\pi}\frac{\lambda}{{\left(\lambda +|x|^2
 \right)^2}}\, .
\end{equation*}
\end{lm}
\

\noindent Following \cite{DP}, one may apply sharp log HLS
inequality~\eqref{loghls}  to deduce that
\begin{align}\label{bdp}
\Fhls[\rho] &= \frac{M}{8\pi}\left( \int_{\R^2}\rho(x)\,\log
\rho(x)\dd x+\frac{2}{M} \iint_{\R^2\times\R^2}\rho(x) \log|x-y| \rho(y)\dd x \dd y\right)\nonumber\\
 &\quad+ \left(1 -   \frac{M}{8\pi}\right)\int_{\R^2}\rho(x)\,\log \rho(x)\dd x\nonumber\\
 &\ge - \frac{M}{8\pi}\,C(M) + \left(1 -   \frac{M}{8\pi}\right)\int_{\R^2}\rho(x)\,\log \rho(x)\dd x\ .
\end{align}
It follows from this and the monotonicity of $\Fhls[\rho(t)]$ that
for solutions $\rho$ of the PKS system for which $M< 8\pi$,
$$
{\cal E}[\rho(t)]:=\int_{\R^2}\rho(t,x)\log \rho(t,x)\dd x \le
\frac{8\pi {\cal F}[\rho_0]  - M\,C(M)}{8\pi - M}\ .
$$
Therefore, for $M< 8\pi$, the {\rm entropy} ${\cal E}[\rho(t)]$
stays bounded from above, uniformly in time. This precludes the
collapse of mass into a point mass for such initial data.
In~\cite{DP,BDP}, this formal analysis is made rigorous, and the
global existence of solutions below the critical value $8\pi$ and
a number of their properties as well are established.

Previous work in this direction, by J\"ager and
Luckhaus~\cite{JL}, had shown that for initial data of
sufficiently small mass, the entropy ${\cal E}[\rho(t)]$ stayed
bounded uniformly in $t$. Their analysis used the
Gagliardo--Nirenberg-Sobolev inequality for functions $f$ in
$\R^2$ that bounds $\|f\|_4$ in terms of $\|\nabla f\|_2$ and
$\|f\|_2$, and not the Log HLS inequality, but their global
existence result requires the mass to lie below a threshold that
is strictly less than $8\pi$.

That $8\pi$ is the actual critical value at which diffusive and
concentrating effects are balanced, and not only a better lower
bound, can be seen by computing moments: When the initial data has
a finite second moment, and $M> 8\pi$ such collapse, or
``blow-up'' does indeed occur in a finite time. To see this, we
first note a weak formulation of our the PKS evolution equation
that will useful to us later on. Let $\psi$ be any test function.
Then
\begin{align}
 \frac{{\rm d}}{{\rm d}t}
 \int_{\R^2}\psi(x)\rho(t,x)\dd x =\,&  \int_{\R^2}\Delta \psi(x)\rho(t,x)\dd x \nonumber\\
 &- \frac{1}{4\pi}\int_{\R^2 \times \R^2}\rho(t,x)\frac{\left(\nabla \psi(x) - \nabla \psi(y)\right)\cdot (x-y)}{|x-y|^2}\,\rho(t,y)\dd x \dd y\ .
 \label{weakform}
\end{align}
In addition to the usual integration by parts, we have {\em
symmetrized} the second term on the right in $x$ and $y$. Fixing
any $a\in \R^2$ and taking $\psi(x) = a\cdot x$, we see
from~\eqref{weakform} that ${\displaystyle  \frac{{\rm d}}{{\rm
d}t}
 \int_{\R^2} x \rho(t,x)\dd x = 0}$; i.e.,
the center of mass is
conserved.
Due to the translational invariance,
we henceforth assume zero center of mass.
More interestingly, taking $\psi(x) = |x|^2$, so that $
\left(\nabla \psi(x) - \nabla \psi(y)\right)\cdot (x-y) = 2$, we find
\begin{equation*}
 \frac{{\rm d}}{{\rm d}t}
 \int_{\R^2} |x|^2\rho(t,x)\dd x =  4M - \frac{1}{2\pi}M^2 = 4M\left(1 - \frac{M}{8\pi}\right)\ .
\end{equation*}
Thus, if $M> 8\pi$, the right hand side is strictly negative, and
this shows that the second moment of $\rho(t)$ reaches zero in a
finite time if {\it initially bounded}, or else some sort of
singularity develops that would invalidate the formal calculation
we have just made.

Thus, the mass value $M = 8\pi$ is the {\em critical mass} for the
PKS system: For $M < 8\pi$, one has global solutions for which
diffusion dominates so that all of the mass tends to infinity as
the time tends to infinity, see~\cite{BDP}, while for $M>8\pi$, solutions develop
singularities, see~\cite{JL}.

Our focus in this paper is on the case $M = 8\pi$.
Notice that for $M = 8\pi$, $\Fhls$ is exactly the functional that is on
the left hand side in~\eqref{loghls}. Since the densities
$\rhohls$ are minimizers of $\Fhls$ for $M = 8\pi$,
it follows that
$$
\frac{\delta {\Fhls}[\rhohls]}{\delta  \rho}(x) = 0\ ,
$$
and then from~\eqref{wasgrad} that each $\rhohls$ -- and each of
their translates --  is a stationary solution of~\eqref{wasgrad};
{\it i.e.}, of~\eqref{eq:sp}. Of course, this can also be checked
directly. {\em Our main goal in this paper is to determine the
stability of these solutions, and to determine basins of
attraction for them. In achieving this goal, we develop several
novel functional inequalities, and a strategy of concentration
control that may be useful elsewhere, and may be the main
contribution of the paper.}

Each of the $\rhohls$ has an infinite second moment, and so shall
all of the functions in the basins of attraction that we find for
them. This must be the case according to previous work \cite{BCM}
on the case $M=8\,\pi$ for initial data with a finite second
moment. The paper \cite{BCM}  proves the global existence  of weak solutions
with finite second moment that satisfy the   free energy
dissipation  inequality
\begin{equation*}
{\Fhls}[\rho(T)]  + \int_0^T\left[ \int_{\R^2}
\rho(t,x)\left|\nabla  \frac{\delta {\Fhls}[\rho(t)]}{\delta
\rho}(x)\right|^2\dd x\right]\dd t \leq  {\Fhls}[\rho(0)] \ ,
\end{equation*}
which is what one would guess should hold from (\ref{enerdis}).
Moreover, \cite{BCM} proves that every such solution blows up at
{\em infinite time}. That is, the $t\to\infty$ limit of any such
solution is a Dirac  mass $8\,\pi$ at the center of mass of the
initial data. Furthermore, a  point mass of mass $M$  is a
stationary measure-valued solution in the sense introduced in
\cite{DS} if and only if $M\geq 8 \pi$. Let us finally mention
that an analysis of basins of attraction of $\rhohls$ in the
radial case was done in \cite{BKLN} in terms of certain relative
moment conditions.

From this point of view, the solutions in the critical
mass case with  finite  initial second moment are choosing in their
large time asymptotics the only possible stationary state  with a finite second
moment.

Let us finally comment that the family of stationary solutions
$\rhohls$ play a role too in the conjectured profile
of blow-up for any point singularity of the solutions for masses
$M>8\pi$. Velazquez has proved  \cite{Velaz,Velaz1}
that the inner part of the matched-asymptotics expansion for the
blow-up profile is given by these stationary solutions for the
critical mass value.


\subsection{The second Lyapunov functional}

The essential tool in our construction and analysis of solutions
of the critical mass PKS system is an interesting and somewhat
surprising interplay between the PKS system and another evolution
equation which also has the $\rhohls$ as stationary solutions --
the Fokker-Planck version of the fast diffusion equation in $\R^2$
with exponent $1/2$:
\begin{equation}\label{pmfp1}
\begin{cases}
 {\displaystyle \frac{\partial u}{\partial t}(t,x) =
 \Delta \sqrt{u(t,x)}  +
 2\sqrt{\frac{\pi}{\lambda M}} \,{\rm div}(x\,u(t,x))}\qquad & t>0\,,\;x\in\R^2\;,\vspace{.2cm}\\
 u(0,x)=u_0(x)\ge 0\qquad &x\in\R^2\, ,
 \end{cases}
\end{equation}
corresponding to the {\em fast diffusion equation} ${\displaystyle \frac{\partial
u}{\partial t} = \Delta \sqrt{u}}$ by a self-similar
change of variable, see \cite{VaPME}. In the interest of brevity we refer
to~\eqref{pmfp1} as the fast--diffusion equation.

This equation can also be written in a form analogous
to~\eqref{wasgrad}: for $\lambda > 0$,  define the  functional
$\Fcfd$ on the non-negative functions in $L^1(\R^2)$ by
\begin{equation*}
 \Fcfd[u]:= \int_{\R^2} \left(\sqrt{u(x)}-\sqrt{\rhohls(x)}\right)^2 \rhohls^{-1/2}(x)\,\dd
 x\, .
\end{equation*}
This functional is the relative entropy of the fast diffusion
equation with respect to the stationary solution
$\rhohls$. The unique minimizer of $\Fcfd$ is
$\rhohls$, and a simple formal computation yields
\begin{equation}\label{cfdvar}
\frac{\delta  \Fcfd[u]}{\delta u}  = \frac{1}{ \sqrt{\rhohls}} - \frac{1}{\sqrt{u}}\ ,
\end{equation}
from which one sees that~\eqref{pmfp1} can be rewritten as
\begin{equation}\label{half}
\frac{\partial u}{\partial t}(t,x) = {\rm
div}\left(u(t,x)\nabla \frac{\delta  \Fcfd[u(t)]}{\delta
u}(x)\right)\ ,
\end{equation}
It follows that for classical solutions $u$ of~\eqref{pmfp1} for
which one can integrate by parts,
\begin{equation*}
\frac{\dd}{\dd t} \Fcfd[u(t)] = -\int_{\R^2}u(t,x)\left|\nabla
\frac{\delta \Fcfd[u(t)]}{\delta u}(x)\right|^2\dd x= -
\int_{\R^2} \left| \frac{1}{2}\nabla \log u(t,x)    +
2\,\sqrt{\frac{\pi}{\lambda\,M}}\,x\sqrt{u(t,x)}\right|^2\dd x\
.\nonumber
\end{equation*}
As one sees from~\eqref{cfdvar} and~\eqref{half}, the densities $\rhohls$
are stationary solutions of the  fast diffusion
equation~\eqref{pmfp1}, as well as the PKS system~\eqref{eq:sp}.
This is much more than a coincidence, and there are very close
connections between the two evolution equations.

Unlike the functional $\Fhls$, the functional $\Fcfd$ is convex.
In fact, it is not only convex in the usual sense, but is also
{\em displacement convex}.  This second type of convexity, arising
in the theory of optimal mass transportation, will play a basic
role in our analysis. We shall explain the relevant aspects of
displacement convexity in Section 2 of this paper.  For now, we
return to the convexity of $\Fcfd$ in the usual sense. Making
simple computations, one finds that
\begin{equation}\label{Hlam3}
 \Fcfd[u]:=  \sqrt{\pi M\lambda} + \int_{\R^2} \left[  \sqrt{\frac{\pi}{M\lambda}} |x|^2u(x) - 2\sqrt{u(x)} + \sqrt{\rhohls(x)}\right]\,\dd x\ .
\end{equation}
Let us define the functionals
$$
{\displaystyle {\cal G}_1[u]:=\int_{\R^2}|x|^2u(x)\dd x} \qquad
\mbox{and} \qquad {\displaystyle {\cal G}_2[u]:=-\int_{\R^2}
\sqrt{u(x)}\dd x}\, .
$$
Since ${\cal G}_1[u]$ is affine on its domain of definition, and
since ${\cal G}_2[u]$ is convex on its domain of definition, one
might formally conclude the convexity of $\Fcfd$ on its domain of
definition. In fact, those who are familiar with displacement
convexity will recognize that functionals ${\cal G}_1[u]$ and
${\cal G}_2[u]$ are displacement convex on their domains of
definition.

Unfortunately, separate consideration of ${\cal G}_1[u]$ and
${\cal G}_2[u]$ is not relevant in our context: Note that
$\sqrt{\rhohls(x)}$ is not integrable, and thus if   $\Fcfd[u]$ is
to be well defined, $\sqrt{u(x)}$ cannot be integrable either.
Furthermore, since $|x|^2\rhohls(x)$ is not integrable, it is
clear that $|x|^2u(x)$ also will not be integrable on the whole
domain of definition of $\Fcfd$.  Thus, {\em cancelations are
crucial} to the definition of $\Fcfd$, and the integral in
(\ref{Hlam3}) cannot be split into a sum of three integrals to be
analyzed separately.

As far as the convexity (in the usual sense) of $\Fcfd$ is
concerned, it is easy to give a rigorous proof: Indeed, $\Fcfd[u]$
can be written as
\begin{equation}\label{Hlam2}
 \Fcfd[u]:= \int_{\R^2} \left[\Phi(u(x))-\Phi(\rhohls(x))-
 \Phi'(\rhohls)(u(x)-\rhohls(x))\right]\,\dd x
\end{equation}
with $\Phi(s)=-2\sqrt{s}$, which is a convex function. However,
displacement convexity  is essential to our strategy, and
even after we have properly introduced the notion of displacement
convexity, we shall have to work much harder to prove that $\Fcfd$
is in fact strictly displacement convex on its domain of
definition.

The convexity properties of $\Fcfd$ are relevant to the analysis
of the PKS system due to the perhaps surprising  fact that $\Fcfd$
is monotone decreasing also along solutions of the critical mass
PKS system~\eqref{eq:sp}, and not only along solutions of  the
fast diffusion equation~\eqref{pmfp1}. {\em This gives us a second
Lyapunov function for the critical mass PKS system}. To see why
this should be so, we make a formal calculation that we shall
revisit in full rigor later on: Let $\rho$ be a sufficiently nice
solution of the PKS system. Then
\begin{align}
\frac{\dd}{\dd t} \Fcfd[\rho(t)] =\,& \int_{\R^2} \frac{\delta
\Fcfd[\rho]}{\delta \rho} {\rm div}\left( \rho(t,x) \nabla \left[
\frac{\delta  \Fhls[\rho]}{\delta \rho}\right] \right)\dd x =
-\int_{\R^2} \rho  \nabla \left[ \frac{\delta \Fcfd[\rho]}{\delta
\rho}\right]
\cdot \nabla \left[ \frac{\delta  \Fhls[\rho]}{\delta \rho}\right] \dd x\nonumber\\
=\,& -\int_{\R^2} \rho  \nabla \left[ \frac{1}{
\sqrt{\rhohls}} - \frac{1}{\sqrt{\rho}} \right]
\cdot \nabla \left[ \log \rho  -G*\rho \right] \dd x\nonumber\\
=\,& -\int_{\R^2}    \left[ 2\,\sqrt{\frac{\pi}{
\lambda\,M}}\,x\,\rho +\nabla \sqrt{\rho} \right] \cdot \nabla
\left[ \log \rho  -G*\rho \right] \dd x \label{keycalc}
\end{align}

Integrating by parts once more on the term involving the
Green's function,
$$
\int_{\R^2}\nabla \sqrt{\rho} \cdot \nabla \left[ \log \rho
-G*\rho \right] \dd x = \frac{1}{2}\int_{\R^2}\frac{|\nabla
\rho|^2}{\rho^{3/2}} + \int_{\R^2}\sqrt{\rho}\,\Delta G*\rho
=\frac{1}{2}\int_{\R^2}\frac{|\nabla \rho|^2}{\rho^{3/2}} -
\int_{\R^2}\rho^{3/2}\dd x\ .
$$
Also, $\int_{\R^2} x \rho\cdot \nabla \log\rho \dd x = -2M$ and, making
the same symmetrization that led to~\eqref{weakform},
\begin{equation}\label{key2}
\int_{\R^2} \rho(x) \, x \cdot \nabla G*\rho(x) \dd x = \frac{1}{4\pi}\int_{\R^2 \times \R^2}\rho(t,x)\,(x-y)\cdot\frac{x-y}{|x-y|^2} \,\rho(t,y)\dd x \dd y = \frac{M^2}{4\pi}\ .
\end{equation}

Using the last three calculations in~\eqref{keycalc}, we find
\begin{equation*}
\frac{\dd}{\dd t} \Fcfd[\rho(t)]  =
-\frac{1}{2}\int_{\R^2}\frac{|\nabla \rho|^2}{\rho^{3/2}}\dd x +
\int_{\R^2}\rho^{3/2}\dd x  +
4\sqrt{\frac{M\,\pi}{\lambda}}\left(1- \frac{M}{8\pi}\right)\ .
\end{equation*}

Notice that the constant term vanishes in critical mass case $M =
8\pi$. Thus, in the critical mass case, formal calculation yields
that for all $T>0$,
\begin{equation}\label{key3}
 \Fcfd[\rho(T)]  + \int_0^T\left[ \frac{1}{2}\int_{\R^2}\frac{|\nabla \rho|^2}{\rho^{3/2}}(t,x)\dd x
- \int_{\R^2}\rho^{3/2}(t,x)\dd x \right] \dd t \leq \Fcfd[\rho_0]
\ .
\end{equation}
In fact, the formal computation yields equality instead of merely
inequality in (\ref{key3}), but it is this inequality that is
useful to us, and this is what we shall actually  prove
for the solutions that we construct here.

The key to exploiting~\eqref{key3} is a particular case of the
Gagliardo-Nirenberg-Sobolev  (GNS) inequalities for which the sharp form
was found by Del Pino and Dolbeault \cite{DD}.

\begin{lm}[Gagliardo-Nirenberg-Sobolev inequality]\label{spgn}
For all functions $f$ in $\R^2$ with a square integrable
distributional gradient $\nabla f$,
\begin{equation*}\label{gns}
\pi \int_{\R^2} |f|^6 \dd x \le \int_{\R^2} |\nabla f|^2 \dd x \int_{\R^2}
|f|^4 \dd x \ ,
\end{equation*}
and there is equality if and only if $f$ is a multiple of a translate of
$\rhohls^{1/4}$ for some $\lambda>0$.
\end{lm}

To apply this, note that at least for strictly positive densities $\rho$,
$$\int_{\R^2}|\nabla \rho^{1/4}(x)|^2\dd x = \frac{1}{16}\int_{\R^2}\frac{|\nabla \rho|^2}{\rho^{3/2}}(x)\dd x\ .$$
Therefore, we define:

\begin{defi}[Entropy dissipation functional]
For any density $\rho$ of total mass $8\pi$ such that $\rho^{1/4}$
has a square-integrable distributional gradient, we define the
{\em entropy dissipation functional} ${\mathcal D}[\rho]$ by
$$
{\mathcal D}[\rho] = 8\int_{\R^2}|\nabla \rho^{1/4}(x)|^2\dd x - \int_{\R^2} \rho^{3/2}(x)\dd x\ ,
$$
and we define ${\mathcal D}[\rho]= \infty$ in all other cases.
\end{defi}

\begin{lm}[Dissipation of $\Fcfd$]
\label{lm:diss}
For all densities $\rho$ of mass $M = 8\pi$,
$$ {\mathcal D}[\rho] \geq 0\ ,$$
and moreover, there is equality if and only $\rho$ is a translate
of $\rhohls$ for some $\lambda>0$.
\end{lm}

\noindent{\bf Proof:} Let $f = \rho^{1/4}$ and note that
$\int_{\R^2} f^4(x)\dd x = 8\pi$. Multiplying $ {\mathcal D}[\rho]
$ through by $\pi$, the claim follows directly from
Lemma~\ref{spgn}. \lanbox

\subsection{The main results on the PKS equation}

The formal result (\ref{key3}) may now be written as
$$
{\displaystyle \Fcfd[\rho(T)]  + \int_0^T{\mathcal D}[\rho(t)] \dd
t \leq \Fcfd[\rho_0] }\,.
$$
Since $\Fcfd[\rho(T)]\ge 0$, this suggests at the very least that
$$
\lim_{T\to\infty}\frac{1}{T}\int_0^T{\mathcal D}[\rho(t)] \dd
t\leq \lim_{T\to\infty}\frac{1}{T} \Fcfd[\rho_0]= 0\ ,
$$
and then Lemma~\ref{lm:diss} suggests  that for all large $t$,
$\rho(t)$ must be close to $\varrho_{\mu}$ for some $\mu>0$.
However, an easy calculation, see Remark~\ref{basinattraction},
using the fact that $\int_{\R^2}|x|^2 \rhohls(x)\dd x = \infty$,
shows that $\Fcfd(\varrho_{\mu}) = \infty$ for $\mu\ne \lambda$.
Therefore, since  $\Fcfd[\rho(t)]$ is non-increasing, one expects
that $\mu = \lambda$. In short, the formal calculations made so
far suggest that for solutions $\rho$ of the PKS system with
initial data $\rho_0$ satisfying $\Fcfd[\rho_0] < \infty$,
$\lim_{t\to\infty}\rho(x,t) = \rhohls(x)$. We now make one more
definition, and then state our main results on the PKS equation:

\begin{defi}[Properly dissipative weak solutions of the PKS equation]\label{prdisso} Let
Let $\rho_0$ be any density on $\R^2$ with mass $8\pi$, such that
for some $\lambda>0$, $\Fcfd[\rho_0]  < \infty$. Let $\rho:
[0,\infty)\to L^1(\R^2)$ satisfy:

\begin{itemize}
\item[(1.5.1)] For each $t\geq 0$, $\rho(t)$ is a continuous curve
of densities of mass $8\pi$ in the sense that for each bounded and
globally Lipschitz function $\psi$ on $\R^2$,
${\displaystyle t \mapsto \int_{\R^2}\psi(x)\,\rho(t,x)\dd x}$
is continuous with $\rho(0) = \rho_0$.

\item[(1.5.2)] For each $T> S\geq 0$, and each smooth and
compactly supported function $\psi$ on $\R^2$,
\begin{align*}
 \int_{\R^2}\psi(x)\rho(T,x)\dd x = &\,\int_{\R^2}\psi(x)\rho(S,x)\dd x + \int_S^T \int_{\R^2}\Delta \psi(x)\rho(t,x)\dd x\dd t \\
 &\,- \frac{1}{4\pi}\int_S^T\int_{\R^2 \times \R^2}\rho(t,x)\frac{\left(\nabla \psi(x) - \nabla \psi(y)\right)\cdot (x-y)}{|x-y|^2}\,\rho(t,y)\dd x \dd y \dd t\ .
\end{align*}

\item[(1.5.3)] For each $T> 0$,
\begin{equation}\label{enendi}
\Fcfd[\rho(T)] + \int_0^T {\mathcal D}[\rho(t)] \dd t \leq
\Fcfd[\rho_0] \ ,
\end{equation}
so that $\rho$ satisfies the  {\em entropy--entropy dissipation}
inequality expected for solutions of the PKS equation.
\end{itemize}
Then $\rho$ is a {\em properly dissipative weak solution} of the
PKS equation \eqref{eq:sp} with initial data $\rho_0$.
\end{defi}

\begin{thm}[Existence and regularity of properly dissipative weak solutions]\label{main}
Let $\rho_0$ be any density on $\R^2$ with mass $8\pi$, such that
$\Fhls[\rho_0] < \infty$, and for some $\lambda>0$, $\Fcfd[\rho_0]
< \infty$.Then there exists a properly dissipative solution of the
PKS equation \eqref{eq:sp} with initial data $\rho_0$. Moreover,
the solutions we construct have additional regularity properties,
including:
\begin{itemize}
\item[{\rm (1.6.1)}] For any $S>0$ and any $p$ with $1 < p <
\infty$, there is a constant $C$ depending only on $S$, $p$,
$\lambda$ and $\Fcfd[\rho_0]$ such that for all $t\geq S$,
$\|\rho(t)\|_p \leq C$.

\item[{\rm (1.6.2)}] The distributional gradient of $\rho^{1/4}$
is square integrable over $[0,\infty)\times\R^2$, and in fact,
$$
\int_0^\infty\!\!\int_{\R^2}|\nabla \rho^{1/4} (t,x)|^2\dd x \dd t
\leq \Fcfd[\rho_0]\ .
$$
\item[{\rm (1.6.3)}]  $\Fhls[\rho(t)]$ is well defined for each
$t$, and is monotone decreasing:
 $\Fhls[\rho(t)]\leq
\Fhls[\rho(s)]$ for all $0\leq s < t$.
\end{itemize}
\end{thm}

\begin{thm}[Basins of attraction]\label{basinat}
Let $\rho_0$ be any density on $\R^2$ with mass $8\pi$, such that
$\Fhls[\rho_0] < \infty$, and for some $\lambda>0$, $\Fcfd[\rho_0]
< \infty$.  Let $\rho$ be any properly dissipative weak solution
of the PKS equation \eqref{eq:sp} with initial data $\rho_0$
satisfying the additional regularity properties {\rm (1.6.1)},
{\rm (1.6.1)} and {\rm (1.6.3)} of Theorem~{\rm\ref{main}}. Then
$$\lim_{t\to\infty}\Fhls[\rho(t)] = \Fhls[\rhohls]\qquad{\rm and}\qquad   \lim_{t\to\infty}\|\rho(t) - \rhohls\|_1 = 0\ .
$$
\end{thm}

Let ${\mathcal B}_\lambda$ denote the set of densities $\rho_0$
of mass $8\pi$ for which $\Fhls[\rho] < \infty$, and
$\Fcfd[\rho] < \infty$.  According to Theorem~\ref{basinat},
${\mathcal B}_\lambda$ is a basin of attraction for
$\varrho_\lambda$ under the PKS evolution in the sense that any
properly dissipative weak solution with initial data in ${\mathcal
B}_\lambda$, and the regularity produced here,  converges strongly to $\varrho_\lambda$ in
$L^1(\R^2)$.

\subsection{Controlled concentration  inequalities}

The proof of the additional regularity in Theorem~\ref{main}, and then Theorem~\ref{basinat}, might at first appear to be
possible by a standard application of entropy-entropy dissipation
methods, given the entropy-entropy dissipation inequality
(\ref{enendi}).  However, this is not the case. The essential
point is that ${\mathcal D}[\rho]$ is {\em not} a convex function
of $\rho$, and even worse, it is a {\em difference} of two
functionals of $\rho$ that can each be arbitrarily large even when
${\mathcal D}[\rho]$ is very close to zero. Indeed, for $M = 8\pi$
and each $\lambda>0$, ${\mathcal D}[\rhohls] = 0$ while
$$
\lim_{\lambda\to 0}\|\rhohls\|_{3/2} = \infty \, , \quad
\lim_{\lambda\to 0}\|\nabla\rhohls^{1/4}\|_{2} = \infty \, , \quad
\mbox{and} \quad \lim_{\lambda\to 0}\rhohls = 8\pi\delta_0\, .
$$
the point mass of $8\pi$ at 0. It follows that the level sets of
${\mathcal D}$ cannot be weakly compact in $L^1(\R^2)$. Likewise,
$\Fhls[\varrho_\lambda]= 8\pi(\log 8 -1)$ for all $\lambda>0$
while $\lim_{\lambda\to 0}\rhohls = 8\pi\delta_0$. Thus, a family
of densities of mass $8\pi$ on which $\Fhls$ is uniformly bounded
need not be weakly compact in $L^1(\R^2)$

In these examples of non-compactness for level sets of ${\mathcal
D}$ and $\Fhls$, we have a family of densities, which, in the
limit,  concentrate {\em all of their mass} at  single point.  {\em
We shall show here that this is essentially the only way
compactness can fail for a family of densities of mass $8\pi$ on
which ${\mathcal D}$ or $\Fhls$  is uniformly bounded.}

Compactness of level sets of $\Fcfd$ fails for more mundane
reasons: A glance at (\ref{Hlam2}) is enough to see that for all
$h>0$, one can  construct a sequence of $\{\rho_k\}$ small
perturbations of $\varrho_\lambda$ lying in $\{\rho\ :\
\Fcfd[\rho] \leq h\}$ such that for some $\epsilon>0$ depending on
$h$, $\liminf_{k\to\infty}\int_{\{|x| < 1/k\}}\rho_k(x) \dd x \geq
\epsilon$. Thus, level sets of $\Fcfd$ are not uniformly
integrable, and thus not even weakly compact in $L^1(\R^2)$.

However, the densities  in level sets of $\Fcfd$  do have a
crucial property: They must have ``thick tails'' and for this
reason, they cannot concentrate more than a limited fraction of
their mass on any given small set. Thus, the examples of
non-compact subsequences in level sets of ${\mathcal D}$ and
$\Fhls$ that we exhibited above do not lie in any level set of
$\Fcfd$, and as we shall show, neither do any other non-compact
subsequences. In Section 3 we prove:

\begin{thm}[Thick Tails]\label{thicka}
Let $\rho$ be any  density of mass $M$ such that $\Fcfd[\rho] <
\infty$. Then  for $\eta_* := \frac15e^{-1/5}$
and any $s>1$
\begin{equation*}
\int_{|x|^2\geq \lambda s^2}  \rho(x)\dd x  \geq \eta_*\,
e^{-\frac{4}{\sqrt{ \pi M \lambda}}\Fcfd[\rho]} \int_{|x|^2\geq
\lambda s^2} \varrho_\lambda(x)\dd x\ = \frac{M\eta_*}{1+s^2}
e^{-\frac{4}{\sqrt{ \pi M \lambda}}\Fcfd[\rho]} .
\end{equation*}
\end{thm}

Though the statement of this theorem makes no reference to optimal
mass transportation, the proof we give in Section 3 relies heavily
on the optimal mass transportation results we present in Section
2, including a new Talagrand type inequality, Theorem~\ref{tal},
involving $\Fcfd$.

Using the {\em very mild} control on concentration  provided by
the Thick Tails Theorem, which says that densities $\rho$ with
$\Fcfd$ finite for some $\lambda$ cannot possibly concentrate {\em
most} of their mass near any one point, we prove two ``compactness
via controlled concentration'' theorems.  Of course, some sort of
concentration control is inherent in {\em any} compactness theorem
for densities, but the point here is that given only the very mild
limit on concentration provided by the Thick Tails Theorem, we are
able, in Section 3,  to prove compactness for the level sets of
$\Fhls$ and ${\mathcal D}$. The first of these theorems  concerns
$\Fhls$.

\begin{thm}[Concentration control for $\Fhls$]\label{both}
Let $\rho$ be any density with mass $M = 8\pi$, with $\Fcfd[\rho]
< \infty$ for some $\lambda > 0$.   Then there exist positive
computable constants $\gamma_1$ and $C_{{\rm CCF}}$ depending only
on $\lambda$ and $\Fcfd[\rho]$ such that
\begin{equation*}
\gamma_1\int_{\R^2}\rho \log_+ \rho \dd x  \le \Fhls[\rho]  +
C_{{\rm CCF}}\,.
\end{equation*}
\end{thm}

Our second `compactness via controlled
concentration''  theorem concerns ${\cal D}$:

\begin{thm}[Concentration control for ${\cal D}$]\label{both2}
Let $\rho$ be any density  with mass $8\pi$,
$\Fhls[\rho]$ finite, and $\Fcfd[\rho]$ finite for some $\lambda>
0$.  Then there exist positive computable  constants $\gamma_2$
and $C_{{\rm CCD}}$ depending only on $\lambda$,  $\Fcfd[\rho]$ and
$\Fhls[\rho]$ such that
\begin{equation*}
\gamma_2\,\int_{\R^2}|\nabla \rho^{1/4}|^2 \dd x \le \pi{\cal
D}[\rho] +C_{{\rm CCD}}\ \, .
\end{equation*}
Moreover, the same conclusion holds if we replace the assumption
that  $\Fhls[\rho]$ finite with the assumption that the entropy
${\mathcal E}[\rho]$ is finite, except that now the constant
depends on the bound on ${\mathcal E}[\rho]$ instead of the bound
on  $\Fhls[\rho]$.
\end{thm}

Theorems~\ref{both} and ~\ref{both2} give us the ``vertical
control'' needed for a compactness result. The horizontal control
is provided by $\Fcfd$ alone. Not only does a bound on
$\Fcfd[\rho]$ ensure that $\rho$ has thick tails, it also ensures
that the tails are not {\em too thick}: A bound on $\Fcfd[\rho]$
provides a bound on all moments of $\rho$ up to but not including
order $2$. Unlike the Thick Tails Theorem, this result is
elementary:

\begin{lm}[Localization]\label{controlcenter}
For all densities $\rho$ with mass $M$ and all $\lambda > 0$, and all $0< q < 2$, there is an explicitly computable
constant $C$ depending only on $q$, $\lambda$ and $M$ so that
\begin{equation*}
\int_{\R^2}|x|^q \rho(x)\dd x \leq C\left(1+
\Fcfd[\rho]\right)^{q/2}\ .
\end{equation*}
\end{lm}

\noindent{\bf Proof:} By the Cauchy-Schwarz inequality
\begin{align}
\int_{\R^2}|\rho - \varrho_\lambda|\varrho_\lambda^{-1/4}\dd x
&\leq  \int_{\R^2}|\sqrt{\rho} - \sqrt{\varrho_\lambda}|\varrho_\lambda^{-1/4}|\sqrt{\rho} + \sqrt{\varrho_\lambda}| \dd x\nonumber\\
&\leq \sqrt{\Fcfd[\rho]} \, \|\sqrt{\rho} +
\sqrt{\varrho_\lambda}\|_2 \leq \sqrt{2M\Fcfd[\rho]} \ .\nonumber
\end{align}
Since $\varrho_\lambda^{3/4}$ is integrable, there is a constant
$C$, depending only on $\lambda$ and $M$, whose explicit form is
easily worked out, for which
$$
\int_{\R^2}\varrho_\lambda^{-1/4} \rho\dd x \leq
\int_{\R^2}\varrho_\lambda^{-1/4} \varrho_\lambda\dd x +
\sqrt{2M\Fcfd[\rho]} \leq C\left(1+ \Fcfd[\rho]\right)^{1/2}\,.
$$

Now repeat the Cauchy-Schwarz inequality, but this time with
$\varrho_\lambda^{-3/8}$ in place of $\varrho_\lambda^{-1/4}$,
then there is a (different) constant $C$, depending only on
$\lambda$ and $M$ such that:
\begin{align}
\int_{\R^2}|\rho - \varrho_\lambda|\varrho_\lambda^{-3/8}\dd x
&\leq  \int_{\R^2}|\sqrt{\rho} - \sqrt{\varrho_\lambda}|\varrho_\lambda^{-1/4}|\sqrt{\rho} + \sqrt{\varrho_\lambda}|\varrho_\lambda^{-1/8} \dd x\nonumber\\
&\leq \sqrt{\Fcfd[\rho]}\, \|\sqrt{\rho}\varrho_\lambda^{-1/8} +
\sqrt{\varrho_\lambda}\varrho_\lambda^{-1/8}\|_2 \leq
\sqrt{\Fcfd[\rho]} \sqrt{C\left(1+ \Fcfd[\rho]\right)^{1/2}}\
.\nonumber
\end{align}
Since $\varrho_\lambda^{5/8}$ is integrable, then by changing the
constant $C$ accordingly, whose explicit form is easily worked
out, we deduce
$$
\int_{\R^2}\varrho_\lambda^{-3/8} \rho\dd x \leq
\int_{\R^2}\varrho_\lambda^{-3/8} \varrho_\lambda\dd x + C\left(1+
\Fcfd[\rho]\right)^{3/4}\,.
$$
The obvious iteration of this argument leads to
$$
\int_{\R^2}\varrho_\lambda^{-r} \rho\dd x \leq
\int_{\R^2}\varrho_\lambda^{-r} \varrho_\lambda \dd x + C\left(1+
\Fcfd[\rho]\right)^{2r}
$$
for each $r$ of the form $1/2 - (1/2)^k$ for $k\in \NN$, and then
by interpolation, for all $0\leq  r < 1/2$, where of course $C$
depends on $r$ as well as $\lambda$ and $M$. To conclude, note
that $\varrho_\lambda^{-r}(x) \sim |x|^{4r}$ for large $|x|$.
\lanbox

\

Lemma~\ref{controlcenter} shows in particular that when
$\Fcfd[\rho] < \infty$, then $\rho\log(e+|x|^2)\in L^1(\RR^2)$, so
that the Newtonian potential of $\rho$ is well defined.  Also,
Lemma~\ref{controlcenter} together with Theorem~\ref{both} shows,
via the Dunford-Pettis Theorem that the intersections of level
sets of $\Fcfd$ and $\Fhls$ are at least weakly compact in
$L^1(\R^d)$, and stronger conclusions follows for the
intersections of level sets of  $\Fcfd$,  $\Fhls$ and ${\mathcal
D}$.

The ``compactness via controlled concentration'' provided by
$\Fcfd$ and its dissipation ${\mathcal D}$ through
Theorem~\ref{both2} and  Lemma \ref{controlcenter} is the core of
our proof of Theorem~\ref{basinat}. However this is not the only
use we shall make of compactness via controlled concentration: It
is absolutely essential to our construction of properly
dissipative weak solutions.

Indeed, in many problems in which one seeks to prove an
entropy-entropy-dissipation inequality such as (\ref{enendi}),
both the entropy functional ${\mathcal H}$ and its dissipation
${\mathcal D}$ would be weakly lower semicontinuous, often due to
some convexity property.   Then, if $\{\rho_n\}_{n\in\NN}$ is a
sequence of nice or approximate solutions of the evolution
equation converging weakly to a weak solution $\rho$, one would
have
$$
{\mathcal H}[\rho(T)] \leq \lim_{n\to\infty} {\mathcal H}[\rho_n(T)]\qquad{\rm and}\qquad
\int_0^T {\mathcal D}[\rho(t)]\dd t   \leq \lim_{n\to\infty}\int_0^T {\mathcal D}[\rho_n(t)]\dd t \ ,
$$
which is very helpful if one is trying to prove something like (\ref{enendi}).

While in our case $\Fcfd$ is convex and lower semicontinuous,
${\mathcal D}$ is the {\em difference} of two non-comparable
convex functions and has no lower semicontinuity. Therefore, we
need new tools to prove (\ref{enendi}), and as we shall see, it is
once again the  compactness via controlled concentration that does
the trick.

\subsection{Further developments}

One can build on the regularity results obtained here to prove
additional regularity. Indeed, if $\rho$ is one of the solutions
we have constructed here, it is easy to prove that  for any $a>0$,
$\nabla c(x,t)$ is bounded and continuous on
$(a,\infty)\times\R^2$, only  using the continuity properties on
$\rho$ in $t$, the uniform control on first moments, and the fact
that $\rho(t)$ is uniformly bounded in both $L^1$ and $L^3$ for
all $t>a$.  Thus ``freezing'' $b := \nabla c$, $\rho$ is seen to
be a weak solution of the linear parabolic equation
$$
\frac{\partial\rho}{\partial t} = \Delta \rho - {\rm div}(b\rho)\
,
$$
with $b$ bounded and continuous. Parabolic regularity theory may
now be applied. In fact, the arguments developed in
\cite{Ali,Kowal,CC} can be applied to get $L^\infty$-bounds of the
density for all positive times based on the $L^p$-bounds, $1\leq
p<\infty$, obtained in Theorem \ref{main}. A further development
that requires new tools is to bound the rate of convergence to the
equilibrium $\varrho_\lambda$ in our convergence theorem.

An interesting problem whose solution would lead to rate
information is to characterize the stability of the GNS inequality
that we have used. That is, we know that ${\mathcal D}[\rho] = 0$
if and only if $\rho$ is a translate of $\varrho_\lambda$ for some
$\lambda > 0$, since, as we have seen, this is simply a
restatement of a sharp GNS inequality of Del Pino and Dolbeault. A
stability result for this inequality would be a result stating
that, for any $\epsilon > 0$, if  ${\mathcal D}[\rho] $ is
sufficiently small, then the distance, in some metric, from $\rho$
to some translate of some $\varrho_\lambda$, $\lambda > 0$, is no
more than $\epsilon$. It would also be useful to quantify the
qualitative stability result for the Log-HLS inequality that we
prove and use in Section 5. Work in this direction is underway.

\subsection{Other equations with a second Lyapunov functional}

The second Lyapunov functional $\Fcfd$ is more useful to us than
the primary Lyapunov functional $\Fhls$, which actually drives the
evolution, because of its convexity properties, especially its
displacement convexity, as explained in Section 3.

There is a ``canonical way'' to produce gradient flow evolution
equations that have a convex second Lyapunov functional that has
been investigated in \cite{MMS}. Indeed, both the PKS equation and
the fast diffusion equation are gradient flow systems where the
gradient is computed using the $2$-Wasserstein metric, as we
recall in Section 3. To keep things simple here, let us explain
the mechanism studied in  \cite{MMS} in the finite dimensional
Euclidean case.

Let $V$ be a smooth convex function on $\R^n$. Let $W$ be the
smooth function on $\R^n$ defined by $W(x) = |\nabla V(x)|^2$. Now
consider the evolution equation
$$
{\displaystyle \frac{{\rm d}}{{\rm d}t}x(t) = - \nabla W(x(t))}\,.
$$
Then of course, for any solution $x$,
$$
\frac{{\rm d}}{{\rm d}t}W(x(t)) = -|\nabla W(x(t))|^2 \leq 0\ ,
$$
and so $W$ is monotone decreasing along the evolution. It is the
primary Lyapunov function for this flow. Next, note that since
$\nabla W = 2[{\rm Hess}V]\,\nabla V$,
$$
\frac{{\rm d}}{{\rm d}t}V(x(t)) = -\left[\nabla V\cdot  \nabla
W\right](x(t)) = - 2\left\{\nabla V\cdot   [{\rm Hess}V]\,\nabla
V\right\}(x(t)) \leq 0\ ,
$$
since the Hessian of $V$ is positive. Thus, $V$ is a second
Lyapunov function for the gradient flow driven by $W$.

An example in   \cite{MMS} concerns a porous medium equation on
the line, which is gradient flow in the $2$-Wasserstein metric for
a certain entropy functional. With this entropy functional playing
the role of $V$, the  gradient flow equation for the functional
corresponding to $W$ is a certain fourth order equation of
thin-film type.

The fact that the entropy for the porous medium equations  is a
second Lyapunov functional for this  fourth order thin film
equation had been discovered earlier in  \cite{CT2} and exploited
as the key to understanding the long time behavior of the latter
equation. Again in this case, the second Lyapunov function is
strictly and uniformly displacement convex, while the primary
Lyapunov functional is not displacement convex at all.

In the case studied here, the second Lyapunov functional does not
arise through the mechanism studied in \cite{MMS}, or any other
evident natural mechanism, and we have no ``explanation'' of why
one should expect $\Fcfd$ to decrease along the PKS flow. However,
as explained in  \cite{MMS}, once one knows this, it is a
consequence, formally at least,  that $\Fhls$ decreases along the
fast diffusion flow. This has interesting consequences  that are
investigated in  \cite{CCL}.

The motivation for doing the computation to check the monotonicity
is twofold: First, both evolution equations have the same steady
states, which is certainly  necessary, but not at all sufficient,
for the computation to work out. Second, there are many sharp
inequalities that have negative powers of $1+|x|^2$ as their cases
of equality, so there are tools available to try to prove the
positivity of the dissipation.

\subsection{A brief outline of the rest of the paper}

The rest of the paper is organized as follows. Section 2 begins
with a brief summary of  some results concerning the
2{\em-Wasserstein metric} and gradient flows with respect to it.
In particular, we  recall a discrete variational scheme due to
Jordan, Kinderlehrer and Otto \cite{JKO} for constructing
solutions of a class of equations including both the PKS equation
and the critical fast diffusion equation.   We also recall
McCann's  \cite{McC2} notion of {\em displacement convexity}, and
explain how this should, at least formally, lead to the
entropy-entropy dissipation inequality that we seek. Making the
formal calculation rigorous will then  be reduced to rigorously
proving certain consequences of displacement convexity for
$\Fcfd$, and this will be facilitated by the ``robustness'' of
displacement convexity.

The latter half of Section 2 is more novel. As we have noted
earlier, $\Fcfd$ is formally a sum of displacement convex terms,
however, for the densities that concern us, each of the terms is
divergent. Thus, we are forced to introduce a {\em regularization}
of $\Fcfd$.  While there are many tools available to regularize
functions that are convex in the usual sense (e.g. infimal
convolution), there is no general approach to regularizing
functionals while preserving, or at least not severely damaging,
their formal displacement convexity properties. The regularization
developed in the second half of Section 2 is one of the
cornerstones of the paper.

In Section 3 we prove the controlled
concentration results that have been stated and discussed in
previous subsections.

In Section 4, we lay the ground work for the proof of
Theorem~\ref{main} on the existence of properly dissipative weak
solutions. These will be constructed using a variant of the
Jordan, Kinderlehrer and Otto \cite{JKO} scheme, which constructs
the evolutions by solving a sequence of variational problems, as
in di Giorgi's ``minimizing steps'' method.

In this method, the Euler-Lagrange equation for the variational
problem solved at each step often provides essential {\em
a-priori} regularity on the minimizing density $\rho$. Once again,
at this point in our problem, we encounter difficulties due to
potential cancelation of infinities. To resolve these, we are
forced to  regularize $\Fhls$. The discrete scheme provides a very
convenient framework in which to impose and control the
regularization: We use a different degree of regularization at
each discrete time step. Because of the regularization, we will at
least know that at each time step, $\nabla \sqrt{\rho}$ is square
integrable, but we shall have no useful quantitative bound on
$\|\nabla \sqrt{\rho}\|_2$. Still, this gives us enough regularity
to make some crucial integrations by parts, and then eventually
through the use of Theorem~\ref{both2}, we shall obtain a useful
quantitative bound on $\|\nabla \rho^{1/4}\|_2$.

In Section 5, we pass to the continuous time limit, and provide
the proofs of Theorems~\ref{main} and  \ref{basinat}. Here, the
flexibility of choosing the degree of regularization at each time
step is crucial to cope with the errors committed in the
displacement convexity of the regularized functional and to get
the right dissipation in the limit. The convergence in $L^1$ needs
a qualitative control of the error in the log-HLS inequality: if
$\Fhls[\rho]$ is small enough, the densities are closed to some
$\varrho_\mu$ in $L^1$ norm. This together with the fact that
$\Fcfd[\rho]$ is non increasing proves that $\mu=\lambda$.

\section{Displacement convexity and the PKS system}\label{metricl}

\subsection{Gradient flows in the Wasserstein metric and displacement convexity}

We recall some facts concerning the $2$-Wasserstein metric that
will be used here. We shall be brief, aiming mainly to establish
terminology and notation. For more background, see
\cite{Villani,AGS}. Let  $\PP(\RR^2)$ denote the set of
probability measures in $\RR^2$, and let $\PP_2(\RR^2)$ the subset
of probability measures with finite second moments.  Define the
functional $\W_2$ in $\PP(\RR^2) \times \PP(\RR^2)$ by
\begin{equation}\label{w2def}
\W_2^2(\mu, \nu) =\inf_{\Pi\in\Gamma} \iint_{\RR^2 \times
\RR^2} \vert x - y \vert^2 \, \dd \Pi(x, y) \ ,
\end{equation}
where $\Pi$ runs over the set $\Gamma$ of all {\em couplings} of
the probability measures $\mu$ and $\nu$; that is, the set of
joint probability measures in $\RR^2 \times \RR^2$ with first
marginal $\mu$ and second marginal $\nu$. For absolutely
continuous probability measures $f\dd x$ and $g\dd x$ we will
simply write $\W_2(f,g)$ in place of  $\W_2(f\dd x,g\dd x)$.
Clearly, $\W_2$ is finite in $\PP_2(\RR^2) \times \PP_2(\RR^2)$,
though it takes on the value $+\infty$ in certain pairs $(\mu,
\nu)\ \in \PP(\RR^2) \times \PP(\RR^2)$ -- for example if $\mu$
belongs to $ \PP_2(\RR^2)$, but $\nu$ does not. It is easy to see
that  $\W_2$ is a metric on $\PP_2(\RR^2)$; it is called the {\em
$2$-Wasserstein metric}, where the~$2$ refers to the exponent  $2$
on the distance $|x-y|$.  More generally, given any $\nu\in
\PP(\RR^2)$,  $\W_2$ is a metric on the subset of $\PP(\RR^2)$
given by $\left\{ \mu \in \PP(\RR^2)\ :\ \W_2 (\mu,\nu) < \infty\
\right\}$.

A result of Brenier \cite{B} as extended by McCann \cite{McC1},
provides effective control over the minimization problem defining
$\W_2(\mu, \nu)$.  To recall this result, let $T$ be a measurable
map $\RR^2 \to \RR^2$. We say that $T$ {\it transports} $\mu$ onto
$\nu$, if for any measurable set $B \subset \RR^2$, $\nu(B)=\mu
\circ T^{-1}(B)$. In this case we say that $\nu$ is the {\it
push-forward} of $\mu$ by $T$, $\nu=T\# \mu$. An equivalent
formulation is that $\nu=T\# \mu$ if
\begin{equation}\label{push}
\int_{\RR^2} \zeta(T(x)) \dd \mu(x) = \int_{\RR^2} \zeta(y)
\dd \nu(y) \quad \forall \zeta \in \mathcal C_b^0(\RR^2)\;.
\end{equation}
By the  Brenier-McCann Theorem \cite{B,McC1}, for any two
probability measures $\mu$ and $\nu$ on $\R^2$ not charging
Haus\-dorff dimension $1$ sets, there is an essentially unique
convex function $\varphi$ in $\R^2$ such that $\nabla\varphi\# \mu
= \nu$ and
\begin{equation}\label{wasmap}
  \W_2^2(\mu, \nu) = \int_{\RR^2} |x-\nabla \varphi(x)|^2\dd\mu(x)\;.
\end{equation}
The essential uniqueness is that if $\varphi$ and
$\widetilde\varphi$ are two such convex functions, then $\nabla
\varphi = \nabla \widetilde\varphi$ almost everywhere with respect
to $\mu$. In this paper we will be concerned with densities whose
mass is not necessarily one. If $\mu$ and $\nu$ are two positive
measures of mass $M>0$, we define $\W_2(\mu, \nu)$ in terms of the
$2$-Wasserstein distance between the probability measures $\mu/M$
and $\nu/M$ as follows:
\begin{equation}\label{w2def2}
  \W_2^2(\mu, \nu) =  M\W_2^2(\mu/M, \nu/M)\ .
\end{equation}
This normalization convention has the advantage that if
$\nabla\varphi\# (\mu/M) = (\nu/M)$, then~\eqref{wasmap} is still
valid for arbitrary $M$. Note that if~\eqref{push} holds for $\mu$
and $\nu$, it also holds if we change $\mu$ and $\nu$ by
multiplying them by a positive constant, i.e., $\nabla\varphi\#
(\mu/M) = (\nu/M)$ if and only if $\nabla\varphi\# \mu = \nu$.

In Section 5 we shall also use the $p$-Wasserstein distance,
$1\leq p <2$, especially for $p=1$, on account of a useful
description of compact sets for this metric. For two probability
measures $\mu$ and $\nu$ on $\R^2$, $p$-Wasserstein distance
$\W_p(\mu,\nu)$ is defined by \eqref{w2def} where $2$ is
substituted by $p$. For two positive measures of mass $M$, we
define $\W_p(\mu, \nu) = \sqrt{M}\,\W_p(\mu/M, \nu/M)$. This
normalization is chosen taking into account \eqref{w2def2} to
extend the standard ordering relation for the $\W_p$-metrics on
probability measures; that is, by H\"older's inequality, we have
for any $1\leq p <2$
\begin{equation}\label{wascom}
\W_p(\mu,\nu) \leq  \W_2(\mu,\nu)\ .
\end{equation}

A fundamental insight of Otto \cite{Otto96} is that the
$2$-Wasserstein metric is useful when considering any evolution
equation on  densities $\rho$ that can be written in the form
\begin{equation}\label{gengrad}
 \frac{\partial \rho}{\partial t}(t,x)  = {\rm div}\left(\rho(t,x)\nabla \left[ \frac{\delta {\cal G}[\rho(t)]}{\delta  \rho}(t,x) \right]\right)\
\end{equation}
for some functional ${\cal G}$. The prime example
of~\eqref{gengrad} considered in \cite{JKO} is the Fokker-Planck
equation for probability densities for which
\begin{equation}\label{fp}
{\cal G}[\rho] = \int_{\R^2}  \rho(x)\,\ln \rho(x)\dd x +
\frac{1}{2}\int_{\R^2} |x|^2\rho(x)\dd x\  .
\end{equation}
In \cite{Otto01} rescaled porous medium equations were also
included. Otto's insight \cite{Otto96,Otto01} is that the equation
(\ref{gengrad}) is {\em gradient flow for the functional ${\cal
G}$ with respect to the $2$-Wasserstein metric}. This is true for
a large class of equations of the form~\eqref{gengrad}, see
\cite{AGS,BCC,CG,CarUlu,
Carrillo-McCann-Villani03,Carrillo-McCann-Villani06,Villani}. The
``gradient flow in the 2-Wasserstein metric'' point of view is
useful to us here for two reasons:
\begin{itemize}
\item It provides the means for constructing well-behaved
solutions of the equation in question through the solution of a
sequence of variation problems; the  {\em
Jordan-Kinderlehrer-Otto} (JKO) {\em scheme}.

\item It provides the means for studying the rate at which
solutions of (\ref{gengrad}) converge to minimizers of the
functional ${\cal G}$, at least when the functional ${\cal G}$ has
a certain convexity property.
\end{itemize}

The convexity property referred to in the second point is McCann's
notion of {\em displacement convexity} \cite{McC2}, which turns
out to be convexity in the ``Riemannian metric'' associated to the
2-Wasserstein metric; see \cite{Otto01}. If the functional ${\cal G}$ is uniformly
displacement convex, then there are automatically a family of
functional inequalities that govern the convergence of solutions
of (\ref{gengrad}) to minimizers of ${\cal G}$. In concrete terms,
the functional ${\cal G}$ is said to be {\em displacement convex}
in case the following is true: For any two densities $\rho_0$ and
$\rho_1$ of the same mass $M$, let $ \varphi$ be  the
essentially unique convex function such that $\nabla \varphi\#
\rho_0 = \rho_1$. For $0 < t < 1$, define
\begin{equation*}
\varphi_t(x) = (1-t)\frac{|x|^2}{2} + t\varphi(x)\qquad{\rm and}\qquad \rho_t = \nabla \varphi_t \# \rho_0\ .
\end{equation*}
The {\em displacement interpolation} between $\rho_0$ and $\rho_1$
is the path of densities $t\mapsto \rho_t$, $0\le t \le 1$. Let
$\gamma$ be any real number. To say that ${\cal G}$ is {\em
$\gamma$-displacement convex} means that for all such densities
$\rho_0$ and $\rho_1$, and all $0\le t \le 1$, $ (1-t){\cal
G}[\rho_0] + t{\cal G}[\rho_1] - {\cal G}[\rho_t]\ge \gamma
t(1-t)\W_2^2(\rho_0,\rho_1)$. ${\cal G}$ is simply {\em
displacement convex} if this is true for $\gamma=0$, and ${\cal
G}$ is {\em uniformly displacement convex} is this is true for
some $\gamma>0$. Let us recall the characterization of
displacement convexity given by McCann in \cite{McC2} for
functionals of the form
\begin{equation}\label{i:func}
{\cal G}_\Phi [\rho]:= \int_{\R^d} \Phi(\rho(x))\,\dx\, ,
\end{equation}
where $\Phi:[0,\infty)\to[0,\infty)$ with  $\Phi(0) =0$. Then
McCann's Theorem says that if $s\mapsto s^d\Phi(s^{-d})$ is convex
non-increasing on $(0,+\infty)$ then the functional ${\cal
G}_\Phi$ is displacement convex, and this condition is essentially
necessary.

A much simpler result, also from \cite{McC2}, is that if $V$ is
any real valued  function on $\R^2$ such that for all $x_0,x_1\in
\R^2$ and all $0\leq t \leq 1$,  $(1-t)V(x_0) + tV(x_1) -
V((1-t)x_0+tx_1) \geq \gamma t(1-t)|x_0-x_1|^2$, then the
functional
$$
{\displaystyle V[\rho] = \int_{\R^2}V(x)\,\rho(x)\dd x}
$$
is $\gamma$-displacement convex. Using these results, one readily
checks that in the case of the Fokker-Planck equation, the
functional \eqref{fp} is indeed 2-displacement convex. The
consequent inequalities that govern the long time behavior of
solutions are {\em Gross's logarithmic Sobolev inequality} and the
{\em Talagrand inequality} for Gaussian measures. Our analysis of
long time behavior for the PKS system falls outside the scope of
previous work in this direction since the functional $\Fhls$ is
{\em not} displacement convex. The key reason that it is useful to
bring the second formal Lyapunov functional $\Fcfd$ into the
analysis of the PKS system is that it is displacement convex. In
the next section we prove the displacement convexity of $\Fcfd$,
and study its consequences.

\subsection{The critical fast diffusion equation as gradient flow of a uniformly displacement convex entropy}

The equation
\begin{equation}\label{pmfp}
\frac{\partial u}{\partial t}(t,x) = \Delta u^m(t,x) + \kappa\,
{\rm div}(x\,u(t,x))\ ,
\end{equation}
where $\kappa$ is a non-negative constant and $m>0$, is called the
{\em porous medium equation} with $\kappa=0$ and $m>1$ while for
$\kappa=0$ and $0<m<1$ is called the fast diffusion equation. When
$\kappa>0$, there is a restoring drift. In case $m=1$,~\eqref{pmfp}
is of course the heat equation  for $\kappa= 0$, and the linear
Fokker-Planck equation for $\kappa>0$.

Equation~\eqref{pmfp} can be written in the gradient flow
form
\begin{equation}\label{pmfpvf}
\frac{\partial u}{\partial t}(t,x) = {\rm div}\,\left(u(t,x)\nabla
\frac{\delta {\cal G}}{\delta u}\right)\quad \mbox{with}\quad
{\cal G}[u] = \int_{\R^2} \left( \frac{1}{m-1} u^m(x) + \kappa
\frac{|x|^2}{2}u(x)\right)\dd x\ ,
\end{equation}
which shows that the evolution equation (\ref{pmfp}) is gradient
flow for ${\cal G}$ with respect to the Wasserstein metric.

The value $m=1/2$ for this  equation in $\R^2$ is {\em critical}
in the sense that  the functional ${\cal G}$ in~\eqref{pmfpvf} is
strictly displacement convex for $m \geq  1/2$, but is not
displacement convex  for $m < 1/2$. There are many of other
``critical'' values of $m$ between $0$ and $1/2$ at which other
things happen, see \cite{CV} for instance. But since displacement
convexity plays a crucial role in our work, it is natural to refer
to the $m=1/2$ case as critical here. Indeed, by the criteria of
McCann introduced above since ${\cal G}[u]=V[u]+{\cal
G}_{\phi_m}[u]$ with $V(x)=\kappa |x|^2/2$ and
$\phi_m(s)=s^m/(m-1)$, then $V[u]$ is uniformly displacement
convex, and for $m\geq 1/2$, ${\cal G}_{\phi_m}[u]$ is
displacement convex. As might be expected, some difficulties arise
at the critical value $m=1/2$.

Since~\eqref{pmfp} is gradient flow for ${\cal G}$, one might hope
to find stable steady states by finding the minimizers $\bar u$ of
${\cal G}$. Computing the Euler-Lagrange equation we find
$m/(m-1)\bar u^{m-1} + \kappa |x|^2/2 = C $, where $C$ is a
Lagrange multiplier for the constraint $M:=\int_{\R^2}u(x)\dd x$,
which is conserved. In the case $m=1/2$ and choosing
$$
\kappa=\kappa_{M,\lambda}:=2\sqrt{\frac{\pi}{M\,\lambda}} \qquad
\mbox{we find} \qquad \bar u(x) =\frac{M}{\pi}\frac{\lambda}
{(\lambda + |x|^2)^2} = \rhohls(x)\ .
$$

One readily checks that $\bar u =  \rhohls$ is a
steady state solution to~\eqref{pmfp} with
$\kappa=\kappa_{M,\lambda}$ and so the family of stationary
solutions of the PKS system which we are investigating are also
stationary solutions of the critical fast diffusion equation for
different drifts $\kappa=\kappa_{M,\lambda}$. However, as neither
$\sqrt{\rhohls}(x)$ nor $|x|^2\rhohls$
is integrable, these functions are not in the domain of definition
of ${\cal G}$, and so are {\em not} minimizers of
\begin{equation}\label{ment2}
{\cal G}[u] = \int_{\R^2} \left( -2 \sqrt{u(x)} +
\kappa_{M,\lambda} \frac{|x|^2}{2}u(x)  \right)\dd x\ ,
\end{equation}
the $m=1/2$ version of~\eqref{pmfp} with
$\kappa=\kappa_{M,\lambda}$ as above.

The cure is a simple renormalization as introduced in
\cite{CV,LM}: Consider instead the functional $\Fcfd[u]$ defined
by
\begin{equation}\label{renorm}
u \mapsto \int_{\R^2}  -2 \left(\sqrt{u(x)}  -
\sqrt{\rhohls(x)}\right)\dd x  +
\kappa_{M,\lambda}\int_{\R^2} \frac{|x|^2}{2}\left[u(x)
-\rhohls(x)\right] \dd x\, .
\end{equation}
Then, at least  as long as $u$ has the same behavior at infinity
as does $\rhohls$, the integrals will converge. The counter terms
that we have subtracted off from our functional do not depend on
$u$, and hence they do not affect $\delta {\cal G}/\delta u$. This
is the key idea used in the improvements of rates of convergence
for the fast diffusion equation, see \cite{CV,LM}. Since
$\kappa_{M,\lambda} |x|^2 = 2/\sqrt{\rhohls} - \lambda
\,\kappa_{M,\lambda}$, the functional in~\eqref{renorm} can be
written in the following simpler form, which we take to be the
definition of the {\em critical fast diffusion entropy}:
\begin{equation*}
\Fcfd[u] := \int_{\R^2}\frac {\left(\sqrt{u} - \sqrt{
\rhohls}\right)^2} {{\sqrt{\rhohls}}
}\dd x .
\end{equation*}
It is easy to check that  for $m=1/2$,~\eqref{pmfp} can be written
in the general form \eqref{gengrad} with ${\cal G}=\Fcfd$. As
noted above the displacement convexity of $\Fcfd$  is {\em
formally} obvious from the fact that where $\sqrt{u(x)}$,
$\sqrt{\rhohls(x)}$ and $|x|^2u(x)$ are integrable, $\Fcfd[u]$
would differ from the right hand side of~\eqref{ment2} by a
constant. We provide a rigorous proof in the next subsection.


\subsection{Regularization of the  critical fast diffusion entropy}

To show that $u\mapsto \Fcfd[u]$ is displacement convex, and more
generally, to make rigorous computations involving critical fast
diffusion entropy, $\Fcfd[u]$, we introduce a regularized version
of the critical fast diffusion entropy:
\begin{defi}[Regularized fast-diffusion relative entropy functional]
For $\delta> 0$, and $u$ a density with mass $M$, define
$\Fcfdr[u]$ by
\begin{equation*}
\Fcfdr[u] = \int_{\R^2}\frac {\left(\sqrt{u+\delta } - \sqrt{
\rhohls+\delta }\right)^2} {\sqrt {\rhohls+\delta} }\dd x\,.
\end{equation*}
\end{defi}
\begin{prop}[Displacement convexity of relative entropy functionals]\label{entap}
For any density $u\in L^1_+(\RR^2)$ of mass $M$, $\delta \mapsto
\Fcfdr[u]$ is monotone increasing as $\delta$ decreases to zero,
and
\begin{equation}\label{regdeg1}
\lim_{\delta\to 0} \Fcfdr[u] =  \Fcfd[u]\ .
\end{equation}
Furthermore, let $u_0$ and $u_1$ belong to $L^1_+(\RR^2)$ of total
mass $M$ such that $\W_2(u_0,u_1) < \infty$, and let $u_t$, $0\le t
\le 1$ be their displacement interpolation.  Then for each
$\delta>0$,
\begin{equation*}
(1-t)\,\Fcfdr[u_0] + t\,\Fcfdr[u_1] - \Fcfdr[u_t] \ge t\,(1-t)\,
K_\delta(u_0,u_1)\ ,
\end{equation*}
where $K_\delta(u_0,u_1)$ satisfies
\begin{equation}\label{dico1}
\lim_{\delta\to 0}K_\delta(u_0,u_1) =
\kappa_{M,\lambda}\,\W_2^2(u_0,u_1)\ ,
\end{equation}
and $K_\delta(u_0,u_1)\geq \gamma_\delta \,\W_2^2(u_0,u_1)$ with
$\gamma_\delta < 0$. Consequently, the maps $u\mapsto \Fcfdr[u]$
are $\gamma_\delta$-displacement convex and the map $u\mapsto
\Fcfd[u]$ is strictly uniformly displacement convex:
\begin{equation*}
(1-t)\,\Fcfd[u_0] + t\,\Fcfd[u_1] - \Fcfd[u_t] \ge
\kappa_{M,\lambda}\,t\,(1-t)\,\W_2^2(u_0,u_1) \ .
\end{equation*}
\end{prop}

\noindent{\bf Proof:} As $\delta \mapsto \left(\sqrt{u+\delta } - \sqrt{ \rhohls+\delta
}\right)^2$ is non-increasing. Then, as $\delta$ decreases,
$$
\frac{\left(\sqrt{u+\delta } - \sqrt{ \rhohls+\delta }\right)^2}
{\sqrt {\rhohls+\delta} }\qquad{\rm increases\ to}\qquad
\frac{\left(\sqrt{u } - \sqrt{ \rhohls }\right)^2}
{\sqrt {\rhohls} }\ .
$$
By the monotone convergence theorem and~\eqref{regdeg1}, the
monotonicity in $\delta$ follows. Next,
\begin{equation*}
\frac{\left(\sqrt{u+\delta } - \sqrt{ \rhohls+\delta }\right)^2}
{\sqrt {\rhohls+\delta} } =\frac{u}{\sqrt {\rhohls+\delta } }
-2\left(\sqrt{u+\delta } -\sqrt\delta\right)+ \frac{\delta}{\sqrt
{\rhohls+\delta } } +  \sqrt {\rhohls+\delta} -2\,\sqrt{\delta}\ .
\end{equation*}
Where by the mean value theorem
$$
\frac{u}{\sqrt{\rhohls+\delta}} \le \frac{u}{\sqrt{\delta}} \,,\;
\sqrt{u + \delta}-\sqrt{\delta} \leq \frac{u}{2\sqrt{\delta}}
\;\mbox{ and }\;\left|\frac\delta{\sqrt{\rhohls+\delta}} + \sqrt
{\rhohls+\delta} -2\,\sqrt{\delta}\right|\le
\frac{\rhohls}{\sqrt\delta}\,.
$$
These three terms are integrable and
\begin{equation}\label{dico5}
\Fcfdr[u] =  {\rm (I)} + {\rm (II)} + {\rm const.}
\end{equation}
where
\begin{equation*}
 {\rm (I)}:= \int_{\R^2}\frac{1}{\sqrt {\rhohls+\delta } }\,u \dd x \quad\mbox{and}\quad {\rm (II)}:= 2\int_{\R^2}\left( \sqrt\delta -\sqrt{u+\delta }\right)\dd
 x\, .
\end{equation*}
The criterion~\eqref{i:func} are easily checked for $\Phi(u) = \sqrt\delta
-\sqrt{u+\delta }$, and thus ${\rm (II)}$ is displacement convex.

The term ${\rm (I)}$ in~\eqref{dico5} is unfortunately {\em not}
displacement convex in general. In fact, we will show that is
$\gamma_\delta$-displacement convex with a explicit computable
constant. In order to check the $\gamma_\delta$-displacement
convexity of the regularized functional, notice that
$({\rhohls+\delta})^{-1/2}$ is a function of
$|x|^2$. Thus, the functional (I) is of the general form:
$$
u\mapsto \int_{\R^2}V_\delta(x)\, u(x)\,\dd x
$$
with $V_\delta(x)=({\rhohls+\delta})^{-1/2}$. The characterization
of $\gamma_\delta$-displacement convexity
\cite{AGS} ensures that this is implied by $D^2
V_\delta \geq \gamma_\delta I_2$. Let us compute the hessian of
the potential $V_\delta(x)$. Define the function $f_\delta$ on
$[0,\infty)$ by
$$
f_\delta(|x|^2) := V_\delta(x) = \frac{1}{\sqrt
{\rhohls(x)+\delta } }\ .
$$
We compute
\begin{equation*}
  f_\delta'(s)=\frac{A}{\left[A  + \delta(\lambda + s)^{2}\right]^{3/2}} \quad\mbox{and} \quad
  f_\delta''(s)=-\frac{3\,A\,\delta(\lambda + s)}{[A + \delta(\lambda + s)^2]^{5/2}}\ ,
\end{equation*}
with $A=M\,\lambda/\pi$. Therefore
$$
D^2 V_\delta (x)= 2 f_\delta'(|x|^2) \delta_{ij} + 4
f_\delta''(|x|^2) (x\otimes x),
$$
and taking into account that $f_\delta''(s)\leq 0$, then
\begin{equation}\label{new}
\xi\cdot D^2 V_\delta (x)\cdot \xi^T \geq \left[2 f_\delta'(|x|^2)
+ 4 f_\delta''(|x|^2) |x|^2\right] |\xi|^2:= F_\delta(|x|^2)
|\xi|^2
\end{equation}
for all $x,\xi\in\R^2$, where the function $F_\delta$ is given by
$$
F_\delta(s)=2 f_\delta'(s) + 4 s
f_\delta''(s)=\frac{2A^2+2A\delta\lambda^2 -8A\lambda\delta s
-10A\delta s^2}{[A + \delta(\lambda + s)^2]^{5/2}}.
$$
It is obvious that the function $F_\delta$ converges point-wise to
the constant $\kappa_{M,\lambda}$ as $\delta \to 0$ in
$[0,\infty)$. Moreover, since for each $\delta>0$, the function
$F_\delta(s)\to 0$ as $s\to \infty$ and it is clear that is
negative for $s$ large enough since the denominator is positive
and the numerator has a negative dominant term, then $F_\delta$
attains its maximum and minimum in $[0,\infty)$. Then, we can
choose its minimum value as $\gamma_\delta<0$ and the
$\gamma_\delta$-displacement convexity is proved.

In order, to show the limiting uniform displacement convexity, we
need to refine our arguments. For that, we come back to the
definition of convexity. Let $\psi$ be the essentially unique
convex function such that $\nabla \psi\# u_0 = u_1$. For $0 \le t
\le 1$, define
\begin{align*}
  \eta_\delta(t) &:=\int_{\R^2}V_\delta(x)\left[(1-t)\,u_0(x) + t\,u_1(x)-u_t(x)\right]\;\dd x \\
&=\int_{\R^2}\left[(1-t)\,V_\delta(x)+
t\,V_\delta(\nabla\psi(x))-V_\delta(x + t(\nabla \psi(x)-
x))\right]\,u_0(x)\dd x \ .
\end{align*}
We seek a lower bound on $\eta_\delta$ of the form $\eta_\delta(t)
\ge t\,(1-t)\,K_\delta(u_0,u_1) $.  Since $\eta_\delta(0) =
\eta_\delta(1) = 0$, it suffices for this purpose to show that
$\eta_\delta ''(t) \ge 2\,K_\delta (u_0,u_1)$ for all $0 \le t \le
1$. By denoting $y:=\nabla \psi(x) - x$, we have
\begin{equation*}
\eta_\delta''(t) = \int_{\R^2} y\cdot D^2 V_\delta(x+ty) \cdot y^T
\,u_0(x)\,\dd x\ .
\end{equation*}
Using \eqref{new}, we
readily obtain that
$$
\eta_\delta''(t) \geq \int_{\R^2} F_\delta(|x+ty|^2) |y|^2
\,u_0(x)\,\dd x \geq K_\delta(u_0,u_1)\ ,
$$
with
$$
K_\delta(u_0,u_1) := \min_{0\leq t\leq 1} \int_{\R^2}
F_\delta(|x+ty|^2) |y|^2 \,u_0(x)\,\dd x .
$$
Now, let us observe that the function $F_\delta$ is bounded in
$[0,\infty)$ uniformly in $\delta$. For that, note that
$f'_\delta$ is decreasing and thus $f'_\delta(s)\leq f'_\delta
(0)\leq A^{-1/2}$. On the other hand, by the geometric-arithmetic
mean inequality, we get
\begin{equation}\label{passingtolimitdelta}
|sf''_\delta(s)|\leq 3\delta^{1/2} \frac{s}{A+\delta s^2} \leq
\frac{3}{2\sqrt{A}}\ .
\end{equation}
As a consequence, we get
$$
\|F_\delta\|_{L^\infty(0,\infty)} \leq \frac{8}{\sqrt{A}},
$$
and thus,
$$
|F_\delta(|x+ty|^2)| |y|^2 \,u_0(x) \leq \frac{8}{\sqrt{A}}
|\nabla \psi(x)|^2 \,u_0(x)\in L^1(\R^2)\ ,
$$
for all $0\leq t\leq 1$. Thus, the dominated convergence theorem
guarantees that
$$
\lim_{\delta \to 0} \int_{\R^2} F_\delta(|x+ty|^2) |y|^2
\,u_0(x)\,\dd x = \kappa_{M,\lambda} \int_{\R^2} |y|^2
\,u_0(x)\,\dd x = \kappa_{M,\lambda} \W_2^2(u_0,u_1)\ ,
$$
uniformly in $0\leq t\leq 1$, which together with the definition
of $K_\delta(u_0,u_1)$ implies the uniform displacement convexity
of the limiting functional $\Fcfd[u]$.  \lanbox
\medskip

Continuing with the notation of Proposition~\ref{entap}, define
the function $h_\delta$ on $[0,1]$ by $h_\delta(t) = \Fcfdr[u_t] -
K_\delta(u_0,u_1)t^2$. Then by Proposition~\ref{entap},
$(1-t)\,h_\delta(0) + t\,h_\delta(1)- h_\delta(t) \ge 0$, so that
$h_\delta$ is convex.  Therefore, for all $t\in (0,1)$,
$$
h_\delta(1) - h_\delta(0) \ge \frac{h_\delta(t) - h_\delta(0)}{t}\ .
$$
This in turn implies that
\begin{equation*}
\Fcfdr[u_1] -  \Fcfdr[u_0] \ge \limsup_{t\to 0} \frac{\Fcfdr[u_t]  - \Fcfdr[u_0] }{t} + K_\delta(u_0,u_1)\ .
\end{equation*}

To compute the $\limsup$ of the right hand side, we treat the two
non-constant  terms ${\rm (I)}$ and ${\rm (II)}$ in~\eqref{dico5}
separately. As we have noted ${\rm (II)}$ is displacement convex,
and by well known theorems on the sub-gradients of displacement
convex functions \cite[Chapter 10]{AGS}, this
part contributes
\begin{equation*}
\int_{\R^2}\frac{\nabla
u_0(x)}{2(u_0(x)+\delta)^{3/2}}\cdot(\nabla\psi(x)
-x)\,u_0(x)\,\dd x\ ,
\end{equation*}
as long as the integrand satisfies mild regularity properties; in
particular whenever $u_0$ is bounded below on every compact set by
some strictly positive number, and $\sqrt{u_0}$ has a square
integrable distributional gradient.  We shall show that both of
these conditions hold in our application.  Given that they do,
then by the Cauchy-Schwarz inequality,
\begin{align*}
\int_{\R^2}\frac{\left|\nabla
u_0\right|}{(u_0+\delta)^{3/2}}\left|\nabla\psi(x) -x\right|u_0
\,\dd x &\le \sqrt{\int_{\R^2}\frac{\left|\nabla
u_0\right|^2}{(u_0+\delta)^{3}}\,u_0\,\dd
x}\;\sqrt{\int_{\R^2}\left|\nabla\psi(x)
-x\right|^2\,u_0\,\dd x}\\
&\le  \sqrt{\frac1{3\delta}\int_{\R^2}\frac{|\nabla
u_0|^2}{u_0}\,\dd x}\,\W_2(u_0,u_1) = \frac2{\sqrt{3\delta}}
 \sqrt{\int_{\R^2}\left| \nabla
\sqrt{u_0}\right|^2\,\dd x}\,\W_2(u_0,u_1)\ .
\end{align*}

The contribution of ${\rm (I)}$ in~\eqref{dico5} can be treated by
appealing to the general results in \cite{AGS}
since this functional is $\gamma_\delta$-displacement convex, in
the notation of \cite{AGS}, and thus, this part
contributes
\begin{equation*}
2\int_{\R^2} f_\delta'(|x|^2)\,x\cdot  (\nabla\psi(x) -x)\,u_0\dd x\ .
\end{equation*}
Which is finite because $s \mapsto \sqrt s\,f'_{\delta}(s)$ is a
rational functional which tends to 0 when $s$ goes to infinity.
Thus we have:
\begin{lm}[First-order characterization of displacement convexity]\label{discobnd}
Let $u_0$ and $u_1$ be two densities of total mass $M$ such that
$\W_2^2(u_0,u_1)< \infty$, and such that $u_0$ is uniformly
bounded below on  compact subsets of $\R^2$ by a strictly positive
number, and that $\sqrt{u_0}$ has a square integrable
distributional gradient.  Let $\nabla\psi$ be the unique gradient
of a convex function $\psi$ in $\R^2$ so that $\nabla\psi \# u_0 =
u_1$. Then
\begin{equation}\label{dico8}
 \Fcfdr[u_1] - \Fcfdr[u_0] \ge\! \int_{\R^2}\!\left[ \frac{2Ax}{\left[A  \!+\! \delta(\lambda\! +\! |x|^2)^{2}\right]^{\frac32}}\! +\! \frac{\nabla u_0}{2(u_0\!+\!\delta)^{\frac32}}\right]\! (\nabla\psi(x) -x)\,u_0\dd x + K_\delta(u_0,u_1)
\end{equation}
where $K_\delta(u_0,u_1)$ is defined in {\rm Proposition~\ref{entap}},
and the integrand in~\eqref{dico8} is integrable.
\end{lm}
One might be tempted to take the limit $\delta\to 0$ at this stage
and to conclude
\begin{equation*}
 \Fcfd[u_0]  \le  \Fcfd[u_1]   -  \int_{\R^2} \left[\kappa_{M,\lambda}\,x +
 \frac{\nabla u_0}{2u_0^{3/2}}\right] \cdot (\nabla \psi(x) -x)\,u_0\,\dd x - \kappa_{M,\lambda}\, \W_2^2(u_0,u_1)\ ,
\end{equation*}
but without further information about $\nabla \psi(x) -x$, it is
not possible to do this, or to justify the convergence of the
integral. In our applications, it will be simpler to use the
specific information that we obtain on $\nabla \psi(x) -x$, then
to do some integrations by parts, {\em and then} take the limit
$\delta\to 0$.

Let us finally deduce as an application of the uniform
displacement convexity of the functional $\Fcfd[u]$, an
interesting functional inequality of Talagrand type. Actually,
generalized Log-Sobolev-type inequalities lead formally to
generalized Talagrand-type inequalities for this functional by
repeating arguments due to Otto and Villani \cite[Theorem 1,
Proposition 1]{OV} in the linear case and generalized in
\cite[Theorem 2.1]{Carrillo-McCann-Villani03}. Here,
we are able to show it in full rigor by the previous approximation
argument.

\begin{thm}[Talagrand's inequality]\label{tal}
Whenever $u\in L^1_+ (\RR^2)$ of mass $M$ with $\Fcfd[u]<\infty$,
then
\begin{equation*}
\W_2(u,\rhohls) \le
\sqrt{\frac{2\Fcfd[u]}{\kappa_{M,\lambda}}}.
\end{equation*}
\end{thm}

\noindent{\bf Proof:} Using Lemma \ref{discobnd} with
$u_0=\rhohls$ and $u_1=u$, we obtain that
$\Fcfdr[u] \ge K_\delta(\rhohls,u)$.
for all $\delta>0$, since $\Fcfdr[\rhohls]=0$ and
$$
\frac{2Ax}{\left[A  + \delta(\lambda + |x|^2)^{2}\right]^{3/2}} +
 \frac{\nabla \rhohls}{2(\rhohls+\delta)^{3/2}} = 0 .
$$
Thus, passing to the limit $\delta \to 0$ taking into account
\eqref{dico1}, the desired inequality is obtained. \lanbox

\begin{remark}[Basins of attraction]\label{basinattraction}
Theorem~{\rm\ref{tal}} tells us that the $2$-Wasserstein distance
of our initial data to the stationary state $\rhohls$ is finite
provided $\Fcfd[u_0]<\infty$. Moreover, each of the equilibrium
solutions $\rhohls$ are infinitely far apart in the $\W_2$ metric:
We can easily check that with $\varphi(x) =
\sqrt{\lambda/\mu}|x|^2/2$, one has $\nabla\varphi\#\varrho_{\mu}
= \rhohls$. Thus, the uniqueness part of Brenier-McCann Theorem
ensures
$$
\W_2^2(\varrho_{\mu} ,\rhohls ) =
\frac{1}{2}\int_{\R^2}\left| \sqrt{\frac{\lambda}{\mu}} x
-x\right|^2 \varrho_{\mu}(x) \dd x = +\infty
$$
since the equilibrium densities $\rhohls$ all have
infinite second moments. In particular,
$\Fcfd[\varrho_{\mu}]=+\infty$ for $\mu\neq\lambda$.
\end{remark}

\section{Proof of the concentration controlled  inequalities}

\subsection{Concentration control for $\Fhls$}

To prepare the way for the proof of Theorem~\ref{both}, it is
useful to give an elementary demonstration of a crude form of the
log HLS inequality, without sharp constants, but which would
nonetheless provide bounds on ${\cal E}[\rho]$ for all $M< 8\pi$.

\begin{lm}[Bounds on the entropy]\label{lemhlssub}
Let $\rho$ be a density of mass $M$ on  $\R^2$ such that $\rho\log
\rho$ and $|x| \rho$ are in $L^1(\R^2)$. Then, for any $\alpha >
1/(8\pi)$, there exists a constant $C(M,\alpha,\lambda)>0$ only
depending on $M$, $\alpha$ and $\lambda$ such that
\begin{equation*}
\frac12 (G)_+*\rho(x) \leq    \alpha \int_{\RR^2} \rho \log \rho \,\dd x  + C(M,\alpha,\lambda) +  4\,\alpha\,M\,\log\left(\int_{\R^2}\sqrt{\lambda +|x|^2}\,\rho\dd x\right)\
\end{equation*}
for al $x\in\RR^2$. For $\alpha$ and $\lambda$ fixed,
$C(M,\alpha,\lambda)$ is monotone increasing in $M$.
\end{lm}

\noindent{\bf Proof:} Recall the following Young type inequality:
For all $s,t>0$,  $st \le s\log s + e^{t-1}$. Then,  for any  $\alpha>0$, we have
\begin{equation}\label{young}
 st = \alpha\,[s\,({t}/{\alpha})] \le \alpha s\log s
+ \alpha e^{t/\alpha-1}\ .
\end{equation}
We now apply this to
$$
\frac{1}{2}(G)_+*\rho(x) = \frac{1}{4\pi}\int_{\R^2}(\log |x-y|)_-\rho(y)\dd y  =
\int_{|x-y|< 1}\frac{1}{4\pi}(-\log |x-y|)\frac{\rho(y)}{\rhohls(y)}\rhohls(y)\dd y\ ,
$$
under the integral sign with weight $\rhohls$, and
with
$$
s =  \frac{\rho(y)}{\rhohls(y)}\qquad{\rm and}\qquad t = -\frac{1}{4\pi}\log |x-y|\ .
$$
Since $\rhohls$ is bounded above by $M/(\lambda\,\pi)$, this yields
\begin{equation}\label{praa}
\frac{1}{2}(G)_+*\rho(x) \le \alpha
\int_{|y-x|<1}\left(\frac{\rho}{\rhohls}\right) \log
\left(\frac{\rho}{\rhohls}\right){\rhohls}
\dd y  + \alpha\left(\frac{M}{\lambda\pi e}\int_{|z|\le
1}\frac{1}{|z|^{1/(4\pi\alpha)}}\dd z\right)\ .
\end{equation}
The second integral on the right converges as long as $\alpha > 1/8\pi$,
in which case, doing the integral explicitly, we find
\begin{equation}\label{prbb}
 \frac{\alpha\,M}{\lambda\,\pi\, e}\int_{|z|\le 1}\frac{1}{|z|^{1/(4\pi\alpha)}}\dd z  =
\frac{M}{\lambda e}\frac{8\pi\alpha^2}{8\pi\alpha -1}\qquad\mbox{ for $8\pi \alpha > 1$}\ .
\end{equation}

To relate the first integral to ${\cal E}[\rho]$, use the fact
that $s \mapsto s \log s$ is bounded below by $-1/e$ to conclude
that
\begin{align*}
\int_{|y-x|<1}\left(\frac{\rho}{\rhohls}\right)
 \log \left(\frac{\rho}{\rhohls}\right){\rhohls} \dd y
 \le\,& \int_{\R^2}\left(\frac{\rho}{\rhohls}\right)
 \log \left(\frac{\rho}{\rhohls}\right){\rhohls} \dd y + \frac{M}{e}\\
 \le\,& \int_{\R^2}\rho\log \rho \dd y -   \int_{\R^2}\rho\log
\rhohls  \dd y +  \frac{M}{e}\, .
\end{align*}
By Jensen's inequality for the concave function $\log$ in $L^1((\rho/M) \dd x)$,
\begin{align}
 \int_{\R^2}\rho \log \rhohls \dd x
&= M\log\left(\frac{\lambda M}{\pi}\right) -4\int_{\R^2}\log\left(\sqrt{\lambda +|x|^2}\right)\,\rho\dd x\nonumber\\
&\ge M\log\left(\frac{\lambda M}{\pi}\right)
-4M\log\left(\frac{1}{M}\int_{\R^2}\sqrt{\lambda + |x|^2}\rho\dd
x\right)\ .\label{stack2}
\end{align}
Using~\eqref{prbb} and~\eqref{stack2} in~\eqref{praa}, we obtain
\begin{align*}
\frac{1}{2}(G)_+*\rho(x)  \le &\,\alpha \int_{\R^2}\rho\log \rho \dd y  +  \frac{M}{\lambda e}\frac{8\pi\alpha^2}{8\pi\alpha -1}\nonumber\\
 &+\alpha\,M\left[\frac{1}{e}  -  \log\left(\frac{\lambda}{\pi}\right) +  3\log M + 4\log\left(\int_{\R^2}\sqrt{\lambda +|x|^2}\,\rho\dd x\right)\right]\ .
\end{align*}
The result follows with
\begin{equation*}
  C(M,\alpha,\lambda):= \alpha\,M\left[\frac{1}{\lambda e}\frac{8\pi\alpha}{8\pi\alpha -1} +\frac{1}{e}  -
  \log\left(\frac{\lambda\,}{\pi}\right) + 3(\log M)_+\right]\, .
\end{equation*}
\lanbox\smallskip

To apply this, let $\rho$ be any density on $\R^2$ of mass $8\pi$
{\em or less}. For any number $R>1$, define
$$
\rho_1(x) := 1_{\{|x| \geq R\}}\rho(x) \qquad {\rm and}\qquad
\rho_2(x) := 1_{\{|x| < R\}}\rho(x)\ .
$$
Also define
$$
\widetilde \rho_1(x) := 1_{\{|x| \geq R-1\}}\rho(x) \qquad {\rm
and}\qquad  \widetilde  \rho_2(x) := 1_{\{|x| < R+1\}}\rho(x)\ .
$$
Then since the support of $(G)_+$ has unit radius,
\begin{eqnarray}\label{gosp}
\int_{\R^2} \rho\, (G)_+*\rho  \dd x &=&   \int_{\R^2} \rho\, (G)_+*\rho_1  \dd x + \int_{\R^2} \rho\, (G)_+*\rho_2  \dd x\nonumber \\
&=&   \int_{\R^2} \widetilde \rho_1\, (G)_+*\rho_1  \dd x +
\int_{\R^2} \widetilde \rho_2\, (G)_+*\rho_2  \dd x \,.
\end{eqnarray}
Now suppose it is possible to choose $R>1$ so that for some
$0<a<8\pi$,
\begin{equation}\label{cut}
 \int_{|x|>R-1} \rho \dd x = \int_{\R^2} \widetilde \rho_1 \dd x \leq 8\pi - a \qquad {\rm and}\qquad
 \int_{|x|<R+1} \rho \dd x = \int_{\R^2} \widetilde \rho_2 \dd x \leq 8\pi - a \ .
 \end{equation}
Then choosing $\alpha = (8\pi - a/2)^{-1}$, and applying the
pointwise bounds from Lemma~\ref{lemhlssub} in \eqref{gosp}, we
obtain that
\begin{equation*}
\frac12   \int_{\R^2} \rho\, (G)_+*\rho \dd x \leq \frac{16\pi
-2a}{16\pi -a}\int_{\R^2}\rho\log \rho \dd x   +
2C(8\pi,\alpha,\lambda) +
32\pi\,\alpha\,\log\left(\int_{\R^2}\sqrt{\lambda +|x|^2}\,\rho\dd
x\right)\ .
\end{equation*}
It follows from this and the obvious fact that
$$
 -\int_{\R^2} \rho\, G*\rho \dd x \ge -\int_{\R^2} \rho\,
(G)_+*\rho \dd x,
$$
that
\begin{equation}\label{gosp3}
\Fhls[\rho]  \geq    \frac{a}{16\pi -a}\int_{\R^2}\rho\log \rho
\dd x   - 2C(8\pi,\alpha,\lambda) -
32\pi\,\alpha\,\log\left(\int_{\R^2}\sqrt{\lambda +|x|^2}\,\rho\dd
x\right)\ .
\end{equation}

If the mass $M$ is {\em less} that $8\pi$, the estimate
(\ref{gosp3}) gives us an upper bound on the entropy of $\rho$ in
terms of $\Fhls[\rho]$ and the first moment of $\rho$. This would
suffice, in place of the sharp logarithmic Sobolev inequality, to
prove that no blow-up occurs or mass less than $8\pi$ as it was
done in \cite{BDP}. This cutting-in-pieces argument is similar but
simpler than the one in \cite[Lemma~3.1]{BCM} used for the bounded
second initial moment case.

If the mass is {\em equal} to $8\pi$, one needs additional
information to find an $0<a<8\pi$ and an $R>1$ for which
\eqref{cut} is true. The additional information in case the second
initial moment is bounded was given by contradicting the
convergence to a Delta Dirac, see \cite[Lemma~3.1]{BCM}. Here, we
need to localize the mass and quantify the tails by using a bound
on $\Fcfd[\rho]$. One of the bounds is easy:

\begin{lm}[Solid core]\label{soco2}
Let $\rho$ be a density of mass $M$ such that for some
$\lambda>0$, $\Fcfd[\rho] < \infty$. Then
\begin{equation}\label{new2}
\int_{\{|x| \geq 4\sqrt{\lambda} +
4(\lambda/M\pi)^{1/4}\sqrt{\Fcfd[\rho]}\}}\rho\dd x \leq
\frac{M}{2}\ .
\end{equation}
\end{lm}

\noindent{\bf Proof:} We start by reminding the bound
\begin{equation*}
  \int_{\R^2}|x|\,\rho\dd x <
  \int_{\R^2}\sqrt{\lambda +|x|^2}\,\rho\dd x \le 2\,\sqrt\lambda\,M +
  2M^{3/4}(\lambda/\pi)^{1/4}\,\sqrt{\Fcfd[\rho]}\;,
\end{equation*}
proved in Lemma~\ref{controlcenter}, here written with explicit
constants. The result is a direct consequence of this bound since:
$$
\mbox{if }\,\, \displaystyle\int_{|x|\leq r}\rho\dd x \leq M/2\, ,
\qquad \mbox{then } \,\,\displaystyle\int_{\R^2}|x|\rho\dd x \geq
rM/2\,.
$$
Choosing $r=4\sqrt{\lambda} +
4(\lambda/M\pi)^{1/4}\sqrt{\Fcfd[\rho]}$, we get a contradiction
unless \eqref{new2} is satisfied. \lanbox\smallskip

Now define
$$
R :=   4\sqrt{\lambda}
\left(1 + \sqrt{
\frac{\Fcfd[\rho]}
{\sqrt{M\pi\lambda}}}
\right) + 1\ .
$$
and then using this value of $R$ define $\rho_1$, $\rho_2$,
$\widetilde \rho_1$ and $\widetilde \rho_2$ as above. By
Lemma~\ref{soco2},
$$
\int_{\R^2}\widetilde \rho_1 \dd x \leq \frac{M}{2} \leq 4\pi\ .
$$
Thus, the left inequality in (\ref{cut}) is valid for this choice
of $R$, and any $a \leq  a_1:=4\pi$. The inequality on the right
in (\ref{cut}) requires more work to achieve. The key is the
following result on ``thick tails'':

\begin{thm}[Thick tails]\label{thick}
Let $\rho$ be a density of mass $M$ such that
$\W_2^2(\rho,\varrho_\lambda) < \infty$. Then for
 $\eta_*:= \frac15e^{-1/5}$ and any
$s>1$,
\begin{equation}\label{fira}
\int_{|x|^2\geq \lambda s^2}  \rho(x)\dd x  \geq \eta_*
e^{-\frac{4}{M \lambda}\W_2^2(\rho,\varrho_\lambda)}
\int_{|x|^2\geq \lambda s^2} \varrho_\lambda(x)\dd x\ =
\frac{M\eta_*}{1+s^2} e^{-\frac{4}{M
\lambda}\W_2^2(\rho,\varrho_\lambda)}
\end{equation}
and
\begin{equation}\label{firab}
\int_{|x|^2\geq \lambda s^2}  \rho(x)\dd x  \geq \eta_*
e^{-\frac{4}{\sqrt{ \pi M \lambda}}\Fcfd[\rho]} \int_{|x|^2\geq
\lambda s^2} \varrho_\lambda(x)\dd x\ = \frac{M\eta_*}{1+s^2}
e^{-\frac{4}{\sqrt{ \pi M \lambda}}\Fcfd[\rho]}.
\end{equation}

\end{thm}

\noindent{\bf Proof:}  Let $0<\eta_*< 1/5$ to be fixed later. Fix
any $s>1$ and define $\eta$ by
\begin{equation}\label{firb}
\eta :=  \frac{\int_{|x|^2\geq \lambda s^2}  \rho(x)\dd x }{
\int_{|x|^2\geq \lambda s^2} \varrho_\lambda(x)\dd x}\ .
\end{equation}
We may assume that $\eta <\eta_*$, or else there is nothing to
prove for this $s$.

Let $\nabla \varphi$ be the gradient of a convex function such that
$\nabla \varphi \# \varrho_\lambda = \rho$.  Define $r := \sqrt{\lambda}s$, and define
$$
A := \{ \ x\ :\ |x| \geq r\quad{\rm and}\quad |\nabla \varphi(x)|
< r\ \} \, .
$$
Then since on $A$, $|\nabla \varphi(x) - x|^2 \geq
(|x| - r)^2$,
\begin{equation*}
\W_2^2(\rho,\varrho_\lambda) = \int_{\R^2}|\nabla \varphi(x) -
x|^2\varrho_\lambda\dd x \geq
 \int_{A} (|x| - r)^2\varrho_\lambda\dd x\ .
\end{equation*}
We now claim that
\begin{equation}\label{fir2}
 \int_{A} \varrho_\lambda\dd x \geq (1-\eta) \int_{|x|\geq r}  \varrho_\lambda (x)\dd x\ .
\end{equation}
Assuming (\ref{fir2}) for the moment,
\begin{equation}\label{fir3}
\W_2^2(\rho,\varrho_\lambda)  \geq  \inf\left\{
 \int_{\widetilde A} (|x| - r)^2\varrho_\lambda\dd x\ :\  \widetilde A\subset \{ |x| \geq r\}\ ,\
  \int_{\widetilde A} \varrho_\lambda\dd x \geq (1-\eta) \int_{|x|\geq r}  \varrho_\lambda (x)\dd x\ \right\}\ .
\end{equation}
By the ``bathtub principle'', the infimum is achieved by choosing
$\widetilde A$ to be the level set of $(|x| -r)^2$ in $\{ |x| \geq
r\}$ that has the minimal mass, that is, $\widetilde A = \{ x\ :\
r \leq  |x| \leq \widetilde r\ \}$ where $\widetilde r$ is such
that
$$
\int_{r \leq |x| \leq \widetilde r} \varrho_\lambda\dd x =
(1-\eta) \int_{|x|\geq r}  \varrho_\lambda (x)\dd x\, ,
$$
which is equivalent to
$$
  \int_{|x|\geq \widetilde r}  \varrho_\lambda (x)\dd x = \eta
   \int_{|x|\geq r}  \varrho_\lambda (x)\dd x\ .
$$
By direct computation, this means $M\lambda/(\lambda +\widetilde
r^2)= \eta  M\lambda/(\lambda + r^2)$, and then since $\eta<1$,
this means that $\widetilde r^2 > r^2/\eta$. In fact, since $\eta
< 1/5$, $\widetilde r > 2r$, and so the optimal set $\widetilde A$
in (\ref{fir3}) contains $\{ x\ :\ 2r \leq |x| \leq
r/\sqrt{\eta}\}$. Therefore, combining the last estimates with
(\ref{fir3}),
\begin{equation*}
\W_2^2(\rho,\varrho_\lambda)  \geq
 \int_{2r \leq |x| \leq r/\sqrt{\eta}} (|x| - r)^2\varrho_\lambda\dd x \geq
 \frac{1}{4}  \int_{2r \leq |x| \leq r/\sqrt{\eta}} |x|^2\varrho_\lambda\dd x \ .
\end{equation*}
Now recalling that $r = \sqrt{\lambda}s$, explicitly calculating  the integral yields
$$
\W_2^2(\rho,\varrho_\lambda)  \geq  \frac{\lambda M}{4} \left[
\ln\left(\frac{1+ s^2/\eta}{1+4s^2}\right) +
\left(4-\frac1{\eta}\right) \frac{s^2}{(1+s^2/\eta)(1+4s^2)}
\right]\ .
$$
Remembering that $s>1$ and $\eta<\frac 15$,
$$
\frac{1+ s^2/\eta}{1+4s^2} \geq \frac{s^2/\eta}{5s^2} =
\frac{1}{5\eta}\qquad \mbox{and} \qquad
\left(4-\frac1{\eta}\right)\frac{s^2}{(1+s^2/\eta)(1+4s^2)} \geq
\left(4-\frac1{\eta}\right)\frac{\eta}{5(1+\eta)} \geq -\frac15.
$$
Therefore, fixing $\eta_* := e^{-1/5}/5<1/5$, we get
$$
\eta \geq \eta_* e^{-\frac{4}{\lambda
M}\W_2^2(\varrho,\varrho_\lambda) }\ .
$$
Combining this with (\ref{firb}) yields (\ref{fira}). Thus, to
prove (\ref{fira}), it suffices to prove  (\ref{fir2}). By
definition,
\begin{align*}
 \int_{A} \varrho_\lambda\dd x &=  \int_{\{|x| > r\}\cap \{|\nabla \varphi(x)| < r\}}
 \varrho_\lambda\dd x \geq \int_{|x| > r }\varrho_\lambda\dd x  -
 \int_{  |\nabla \varphi(x)| \geq r }\varrho_\lambda\dd x\\
 & = \int_{|x| > r  }\varrho_\lambda\dd x  -
 \int_{   |y| \geq r }\rho \dd y = \int_{|x| > r}\varrho_\lambda\dd x  -
\eta \int_{|x| \geq r}\varrho_\lambda \dd x \ ,
\end{align*}
and this proves (\ref{fir2}).  Finally, (\ref{firab}) follows from
(\ref{fira}) and the Talagrand-type inequality, Theorem~\ref{tal}.
\lanbox\smallskip

Now note that, by the definition of $R$, $R^2/\lambda > 1$, and
hence Theorem~\ref{thick} implies that
\begin{equation*}
\int_{|x|\geq R+1}  \rho(x)\dd x  \geq \eta_* e^{-\frac{4}{\sqrt{
\pi M \lambda}}\Fcfd[\rho]} \int_{|x|\geq R+1}
\varrho_\lambda(x)\dd x\ .
\end{equation*}
Then since
$$
\int_{|x|\geq R+1} \varrho_\lambda(x)\dd x \geq \frac{8\pi\lambda}{(R+1)^2}\ .
$$
Thus, with $R$ as above and
$$
a_2 :=  8\pi\lambda\eta_* e^{-\frac{4}{\sqrt{ \pi M
\lambda}}\Fcfd[\rho]} \left( 4\sqrt{\lambda} \left(1 + \sqrt{
\frac{\Fcfd[\rho]} {\sqrt{M\pi\lambda}}} \right) +
2\right)^{-2}<8\pi\ ,
$$
the second inequality in (\ref{cut}) is also satisfied.

\medskip

\noindent{\bf Proof of Theorem~\ref{both}:} Most of the work is
now done since \eqref{cut} is satisfied by choosing
$a=\min(a_1,a_2)$. We need a final argument to take care of the
control for the negative contribution of the entropy in terms of
the localization of the mass of the distribution known as
Carleman-type estimate.

\begin{lm}[Control on the negative part of the entropy]\label{logminus}
For any density $\rho\in L^1_+(\R^2)$, if the moment $\int_{\R^2}
m(x)\,\rho(x)\dx$ is bounded with $e^{-m(x)}\in L^1(\R^2)$ and
$m:\R^+_0 \longrightarrow \R^+_0$, then
\[
\int_{\R^2}\rho(x)\,\log_- \rho(x)\,\dd x \leq
\int_{\R^2}m(x)\rho(x)\,\dx + \frac 1e \int_{\R^2} e^{-m(x)} \dd
x\;.
\]
\end{lm}

\noindent{\bf Proof:} Let $\bar \rho:=\rho\,\chi_{\{\rho\le 1\}}$
and $\bar M=\int_{\R^2}\bar \rho(x)\dd x\le\int_{\R^2} \rho(x)\dd
x=M$. Then
\[
\int_{\R^2}\bar \rho(x)\,\left(\log\bar \rho(x)+m(x)\right)\dd
x=\int_{\R^2} [U(x)\log U(x)]\mu\dd x-\bar M \log Z
\]
where $U:=\bar \rho/\mu$, $\mu(x)= e^{-m(x)}/Z$ with
$Z=\int_{\R^2} e^{-m(x)} \dd x$. The Jensen inequality yields
\[
\int_{\R^2} [U(x)\,\log U(x)]\mu \dd x \ge \left(\int_{\R^2}
U(x)\mu \dd x\right)\;\log\left(\int_{\R^2} U(x)\mu \dd
x\right)=\bar M\,\log \bar M\;
\]
and
\begin{align*}
-\int_{\R^2}\rho(x)\,\log_- \rho(x)\,\dd x = \int_{\R^2}\bar
\rho(x)\,\log\bar \rho(x)\dd x & \ge \bar M\log
\bar M - \bar M \log Z - \int_{\R^2}m(x)\,\bar \rho(x)\dd x\\
& \ge -\frac Ze-\int_{\R^d}m(x)\, \rho(x)\dd x\;.
\end{align*}\lanbox\smallskip

\noindent To control the negative part of the entropy in
\eqref{gosp3}, use Lemma~\ref{logminus} with
$m(x)=\sqrt{\lambda+|x|^2}$ and Lemma~\ref{controlcenter}:
\begin{align*}
-\frac{a}{16\pi -a} \int_{\R^2}\rho \log_- \rho \dd x \geq \,&-
\int_{\R^2}\rho \log_- \rho \geq -\int_{\R^2}m(x)\rho(x)\,\dx -
\frac 1e \int_{\R^2} e^{-m(x)}
\dd x\\
\geq \,& - 2\sqrt\lambda M - 2M^{3/4}(\lambda/\pi)^{1/4}\,
\sqrt{\Fcfd[\rho]} - \frac 1e \int_{\R^2} e^{-m(x)} \dd x\ .
\end{align*}
This gives the final control on the positive part of the entropy
from \eqref{gosp3}:
\begin{align}
\Fhls[\rho] \ge &\, \frac{a}{16\pi -a} \int_{\R^2}\rho \log_+ \rho
\dd x - 2C(8\pi,\alpha,\lambda) -
32\pi\,\alpha\,\log\left(\int_{\R^2}\sqrt{\lambda +|x|^2}\,\rho\dd
x\right) \nonumber\\&- 2\sqrt\lambda M -
2M^{3/4}(\lambda/\pi)^{1/4}\, \sqrt{\Fcfd[\rho]} - \frac 1e
\int_{\R^2} e^{-m(x)} \dd x\ \ .\label{final}
\end{align}
Finally, we choose
$$
\gamma_1 := \frac{a}{16\pi -a}
$$
where $a$ is given just above and an explicit expression for
$C_{{\rm CCF}}$ follows from (\ref{final}). \lanbox


\subsection{Concentration control  for ${\cal D}$}

\noindent{\bf Proof of Theorem~\ref{both2}}. While in the proof of
Theorem~\ref{both}, we used a ``horizontal splitting'' of $\rho$,
here we use a ``vertical splitting'': Let $f := \rho^{1/4}$. For
$\beta>0$,  define $f_\beta:=\min\{f,\beta^{1/4}\}$ and
$h_\beta:=f-f_\beta$. We have
\begin{equation}\label{removef}
  \pi\,{\cal D}[\rho]=8\,\pi\,\int_{\R^2}|\nabla f|^2\dd x - \pi \int_{\R^2}f^6\dd x\;.
\end{equation}
Defining $A_\beta=\{x\,:\,f(x)\ge \beta^{1/4}\}=\{x\,:\,\rho(x)\ge
\beta\}$, we get
\begin{equation}\label{decoup}
  \int_{\R^2}f^6\,\dd x = \int_{\R^2\setminus A_\beta}f_\beta^6\,\dd x +  \int_{A_\beta}(h_\beta +\beta^{1/4})^6\,\dd x =
  \int_{\R^2}f_\beta^6\,\dd x - \beta^{3/2}|A_\beta| + \int_{A_\beta}(h_\beta +\beta^{1/4})^6\,\dd x\;.
\end{equation}
By the convexity of $x \mapsto x^6$, for any $\eta \in (0,1)$
\begin{equation}\label{convex6}
 \int_{A_\beta}(h_\beta +\beta^{1/4})^6 \,\dd x \le
 \frac{\beta^{3/2}}{\eta^5}|A_\beta| +\frac{1}{(1-\eta)^5}\int_{A_\beta} h_\beta^6\,\dd
 x\, .
\end{equation}
By the inequality $f_\beta^6 \le \sqrt \beta f^4$, and
plugging~\eqref{convex6} and~\eqref{decoup} into~\eqref{removef},
we obtain
\begin{equation}\label{eq:tmpboth}
  \pi\,{\cal D}[\rho] \ge 8\,\pi\,\int_{\R^2}|\nabla f|^2 \,\dd x - \frac{\pi}{(1-\eta)^5}\int_{\R^2}h_\beta^6\,\dd x  - 8\,\pi^2\sqrt \beta- \pi\,\beta^{3/2}\left(\frac{1}{\eta^5}-1\right)|A_\beta|\;.
\end{equation}
By the GNS inequality, Lemma~\ref{spgn}, applied to $h_\beta$:
\begin{equation}\label{gnsh}
 -\frac{\pi}{(1-\eta)^5} \int_{\R^2}h_\beta^6\dd x  \geq  - \int_{\R^2}|\nabla h_\beta|^2\dd x \left(\frac{1}{(1-\eta)^5}\int_{\R^2}h_\beta^4\dd x\right) \, .
\end{equation}
By definition of $f_\beta$ and $h_\beta$, $\nabla f_\beta=0$ in the support of $h_\beta$ so that
\begin{equation}\label{ffbeta}
\int_{\R^2}|\nabla f|^2\dd x=  \int_{\R^2}|\nabla f_\beta|^2\dd x
+ \int_{\R^2}|\nabla h_\beta|^2\dd x \,.
\end{equation}
Using ~\eqref{gnsh} and~\eqref{ffbeta}  in (\ref{eq:tmpboth}), we obtain
\begin{equation}\label{dosa}
  \pi\,{\cal D}[\rho] \geq 8\,\pi\,\int_{\R^2}|\nabla f_\beta |^2 \,\dd x
  + \left(8\pi -  \frac{1}{(1-\eta)^5}\int_{\R^2}h_\beta^4\dd x \right) \int_{\R^2}|\nabla h_\beta|^2\dd x
  - 8\,\pi^2\sqrt \beta- \pi\,\beta^{3/2}\left(\frac{1}{\eta^5}-1\right)|A_\beta|\ .
\end{equation}

We obtain a result of the type we seek under any conditions that
ensure the second term on the right is positive. Our first
approach uses Theorem~\ref{both}, and so requires that both
$\Fhls[\rho]$, and $\Fcfd[\rho]$, for some $\lambda>0$, be finite.

By (\ref{young}) once more and Theorem~{\ref{both}, we get for all
$\alpha>0$,
\begin{equation}\label{dosa2}
\int_{\R^2}h_\beta^4\dd x = \int_{\R^2} 1_{A_\beta}\rho \dd x \leq
\alpha \left(  \int_{\R^2}\rho \log \rho \dd x  +
e^{(1/\alpha)-1}|A_\beta|\right) \leq \alpha \frac{ \Fhls[\rho] +
{{\rm C_{\rm CCF}}}}{\gamma_1} + \alpha e^{(1/\alpha)-1}|A_\beta|\
.
\end{equation}
This sort of estimate is frequently used in large deviations
problems. Now choose $\eta =1/2$ in (\ref{dosa}), and then choose
$\alpha$ so that the first term in the right hand side in
(\ref{dosa2}) is $\pi/8$, we have, for this choice of $\eta$ and
$\alpha$,
$$\left(8\pi -  \frac{1}{(1-\eta)^5}\int_{\R^2}h_\beta^4\dd x \right)   \geq 4\pi  -  \alpha e^{(1/\alpha)-1}|A_\beta|\ .$$
By Chebychev's inequality, $|A_\beta| \leq \frac{8\pi}{\beta}$,
and so we can choose $\beta$ so that for any $\gamma\in (0,4\pi)$,
$$
\left(8\pi -  \frac{1}{(1-\eta)^5}\int_{\R^2}h_\beta^4\dd x \right) \geq \gamma\ ,
$$
as was to be shown with $C_{\rm CCD} := 248 \pi^2\sqrt \beta$.
Though we have explained how to compute $\beta$, we shall not
write down a formula.  This proves the first part of the theorem.

As for the second, note that we used the bound on $\Fhls$ only to
obtain a bound on the entropy which was used in (\ref{dosa2}).
However, if we have by other means a bound on the entropy, we can
use that in (\ref{dosa2}) in place of the bound on $\Fhls$. This
proves the final part of the theorem.
 \lanbox

\section{Analysis of the discrete time variational scheme for the critical mass PKS system}

From now on, we will assume that the mass is $8\pi$.

\subsection{The Jordan-Kinderlehrer-Otto  scheme for the critical mass PKS system}
The Jordan-Kinderlehrer-Otto (JKO) scheme for constructing
solutions to~\eqref{gengrad}, as described in
Section~\ref{metricl}, would be to fix a time step $\tau>0$, and
inductively define the sequence $\{\rho^k\}_{k\in \N}$ by setting
$\rho_0$ to be the initial density, and then for $k\ge 0$,
\begin{equation}\label{fp45}
\rho^{k+1} \in {\rm argmin}\left\{
\frac{\W_2^2(\rho,\rho^k)}{2\tau} + {\cal G}[\rho]\right\}\ .
\end{equation}

In other words, $\rho^{k+1}$ is {\em some} minimizer of the
functional $\rho \mapsto \W_2^2(\rho,\rho^k)/(2\tau) + {\cal
G}[\rho]$. Only existence of the minimizer is an issue, and not
uniqueness, although in many examples that have been investigated
a strict convexity argument furnishes the uniqueness. The key
point is existence of a minimizer, since that provides a solution
of the Euler-Lagrange equation for the minimization problem in
(\ref{fp45}). Then, as shown in \cite{JKO}, the fact that each
$\rho^{k+1}$ satisfies this Euler-Lagrange equation means that,
under certain conditions on ${\cal G}$, if one defines
$\tilde\rho_\tau$ by
$$
\tilde\rho_\tau(t,x) = \rho^k (x)\qquad {\rm for} \quad k\tau \leq
t < (k+1)\tau\ , \quad \mbox{with } \tilde\rho_\tau(0,x) = \rho_0
\, ,
$$
there is a sequence of values of $\tau$ tending to zero along
which $\tilde\rho_\tau$ tends to a solution of  (\ref{gengrad}) in
a suitable weak sense. This scheme of constructing  weak solutions
of the PKS system for $M < 8\pi$ was developed in \cite{BCC}.
However, for $M=8\pi$ we can not proceed in a very direct manner.
Our problem lies outside the scope of previous applications of the
JKO scheme, since  at the critical mass $M = 8\pi$, (\ref{bdp})
provides no upper bound on ${\cal E}[\rho]$, and hence, it is not
even clear that minimizers exist for the variational problem in
(\ref{fp45}) when ${\cal G} = \Fhls$ and $M = 8\pi$.  Our
controlled concentration inequalities could be used to solve this
problem, but other more thorny issues arise when we would try to
analyze the Euler-Lagrange equation.

To circumvent these difficulties, we introduce a regularized
functional.  In fact, for reasons that will become evident later
on, we shall even be forced to choose a different degree of
regularization  at each time step.

\subsection{Regularization of $\Fhls$}

Let
$$
\gamma(x) := \frac{1}{2\pi}e^{-|x|^2/2}
$$
be the standard Gaussian probability density in $\R^2$. Then, for
all $\epsilon>0$ define $\gamma_\epsilon(x) =
\epsilon^{-2}\gamma\left(x/\epsilon\right)$, and define the
regularized Green's function $G_\epsilon =
\gamma_\epsilon*G*\gamma_\epsilon$, where $*$ denotes convolution,
and $G(x)=-1/(2\pi)\log|x|$. The radially symmetric, $C^\infty$
probability density $\gamma_{\sqrt{t}}$ is the fundamental
solution of the heat equation satisfying for $t>0$, $x\in\R^2$:
$$
\frac{\partial}{{\partial t}}\gamma_{\sqrt{t}} = \frac12 \Delta   \gamma_{\sqrt{t}} \ .
$$
It follows that
\begin{equation}\label{deriveps}
\frac{{\partial}}{{\partial t}} G*\gamma_{\sqrt{t}}  = \frac12
\Delta  G*\gamma_{\sqrt{t}} = - \frac12 \gamma_{\sqrt{t}}\ ,\mbox{
and thus, } \,\frac{{\partial}}{{\partial \epsilon}} G_\epsilon =
- 2\epsilon \gamma_{2\epsilon}\, .
\end{equation}
The right hand side is strictly negative everywhere. From this we
deduce a useful monotonicity property:
\begin{equation}\label{moneps}
\epsilon_1 < \epsilon_2 \quad \Rightarrow \quad    G_{\epsilon_1}(x)  >  G_{\epsilon_2}(x)
\end{equation}
for all $x\in\R^2$. Let us point out that
\begin{equation}\label{deriveps2}
\frac{{\partial^2}}{{\partial \epsilon^2}} G_\epsilon = - 2
\gamma_{2\epsilon}-8\epsilon^2 \Delta \gamma_{2\epsilon}= 2
x\cdot\nabla \gamma_{2\epsilon} + 2\gamma_{2\epsilon} =
2\left(1-\frac{|x|^2}{4\epsilon^2} \right) \gamma_{2\epsilon}\, ,
\end{equation}
since $\gamma_\epsilon$ satisfies $\epsilon^2 \Delta
\gamma_{\epsilon} + \mbox{div } (x \gamma_\epsilon)=0$.

\begin{lm}[First properties of $G_\epsilon$]\label{ptyofG}
Let $G_\epsilon$ be defined as above then:
\begin{itemize}
\item[(i)] For all $x\in \R^2$, $G_\epsilon(x) \le G(x)$.

\item[(ii)] There exists $C>0$ such that for all $x\in \R^2$,
$G_\epsilon(x) \le C\epsilon^{-2}$.

\item[(iii)] For all $(x,y)\in \R^4$,
  \begin{equation*}
  G_\epsilon(x-y)\ge  -\frac{1}{4\pi}\left[4 + \log\left(e+|x|^2 \right) + \log\left(e+|y|^2 \right)\right]\;.
  \end{equation*}
\end{itemize}
\end{lm}

\noindent{\bf Proof:} (i) As $\gamma$ is radially symmetric and
subharmonic in $\R^2$ so that, by the mean value property the
first item holds.

\noindent(ii) Since $\ln_- |x|$ is locally integrable in $\R^2$,
for any $x \in \R^2$
\begin{equation*}
  G*\gamma_\epsilon(x)=\int_{|y|\le \epsilon} G(x-y)\,\gamma_\epsilon(y)\,\dd y \le \frac1{2\pi}
  \int_{|y|\le \epsilon} \log_-|x-y|\,\gamma_\epsilon(y)\,\dd y \le \frac{C}{\epsilon^2}\;.
\end{equation*}
since $\gamma$ is bounded. Thus, we get
\begin{equation*}
  G_\epsilon(w) \le \frac{C}{\epsilon^2} \int_{\R^2} \gamma_\epsilon(z)\,\dd z = \frac{C}{\epsilon^2}\;.
\end{equation*}

\noindent(iii) From the elementary inequality $|z-w| \le |z| + |w| \le 2\max\{|z|\, ,\, |w|\}$, we obtain
\begin{equation*}
\log |z-w| \le \log 2 + \log|z| + \log|w|\ .
\end{equation*}
Therefore,
\begin{equation*}
G(z-w) \ge - \frac{1}{2\pi}(2 + \log|z| + \log|w|)\ .
\end{equation*}
Integrating both sides against $\gamma_\epsilon(x-z)\,\gamma_\epsilon(y-w)$, and using Jensen's
inequality, we find
\begin{align*}
G_\epsilon(x-y)&\ge -\frac{1}{2\pi}\left[2 + \log\left(\int_{\R^2}|z|\,\gamma_\epsilon(x-z)\dd z\right)
+ \log\left(\int_{\R^2}|w|\,\gamma_\epsilon(y-w)\dd w\right)\right]\\
&\ge -\frac{1}{2\pi}\left[2 + \log\left(|x| + \int_{\R^2}|z|\,\gamma_\epsilon(z)\dd z\right)
+ \log\left(|y| + \int_{\R^2}|w|\,\gamma_\epsilon(w)\dd w\right)\right]\\
&\ge -\frac{1}{4\pi}\left[4 + \log\left(e+|x|^2 \right)+ \log\left(e+|y|^2 \right)\right]\ ,
\end{align*}
at least for $\epsilon$ small enough so that
$\int_{\R^2}|z|\gamma_\epsilon(z)\dd z$ is small enough.
\lanbox\medskip

One of the main uses that we will make of the regularization of
the self interaction functional is that it provides a regularized
density for the chemical attractant: given a mass density $\rho$,
we define the {\em regularized chemical attractant density}
$c_\epsilon$ by $c_\epsilon(x)  =    G_\epsilon *\rho(x)$.

\begin{lm}[Uniform estimate regularized chemoattractant]\label{regconc}
For all $\epsilon> 0$ and all densities $\rho$ with mass $8\pi$,
the regularized chemical attractant density $c_\epsilon =
G_\epsilon*\rho$ satisfies
\begin{equation*}
\|\nabla c_\epsilon\|_{\infty} \le \frac{4 \,C_{{\rm
HLS}}}{\epsilon} \,\|\gamma\|^2_{4/3}\
\end{equation*}
and
\begin{equation}\label{regconc2}
\||x| \nabla c_\epsilon\|_{\infty} \le 8C_{{\rm
HLS}}\,\|\gamma\|_{4/3} \||x|\gamma\|_{4/3} + 4 +\,\frac{C_{{\rm
HLS}}}{2\pi\epsilon}\,\|\gamma\|_{4/3}^2\,\| |x|\rho\|_{1} \, .
\end{equation}
\end{lm}
Here $C_{{\rm HLS}}$ denotes the constant of the sharp
Hardy-Littlewood-Sobolev (HLS) inequality \cite{L83} for the special case $p=q =4/3$:

\begin{equation}\label{spclhls}
\iint_{\R^d \times\R^d}f(x)\frac{1}{|x-y|} g(y)\dd x\dd y\le
C_{{\rm HLS}}\|f\|_{4/3}\|g\|_{4/3}\ .
\end{equation}
Though the explicit value of $C_{{\rm HLS}}$ is simple enough,
see~\cite{L83},
our bounds and their proofs will perhaps be easier to
read if leave $C_{{\rm HLS}}$ unevaluated in them, as a marker of
the use of the HLS inequality.

\smallskip\noindent{\bf Proof:} By the Young inequality, we have
$\|\nabla c_\epsilon\|_{\infty} = \|\nabla G_\epsilon *
\rho\|_{\infty} \le 8\,\pi\,\|\nabla G_\epsilon\|_\infty$. And by
the HLS inequality,
$$
\|\nabla G_\epsilon\|_\infty \le \frac{1}{2\pi}\iint_{\R^2 \times
\R^2}
\gamma_\epsilon(x-z)\,\frac{1}{|z-w|}\,\gamma_\epsilon(w-y)\dd
z\dd w \le \frac{C_{{\rm HLS}}}{2\pi}\,\|\gamma_\epsilon\|_{4/3}^2
= \frac{C_{{\rm HLS}}}{2\pi\epsilon}\,\|\gamma\|_{4/3}^2\ .
$$

Using the triangle inequality $|x|\le |x-z|+|z-w|+|w-y|+|y|$ we have
\begin{equation*}
\left||x|\nabla c_\epsilon \right| \le 2 \left|(|x|
\gamma_\epsilon)* \nabla G* \gamma_\epsilon * \rho \right| +
\left| \gamma_\epsilon * \rho \right| + \left|\gamma_\epsilon*
\nabla G* \gamma_\epsilon * |x|\rho \right| \, .
\end{equation*}
Using the Young  and HLS inequalities, see~\eqref{spclhls}, we obtain
\begin{equation*}
  \left\Vert|x|\nabla c_\epsilon \right\Vert_\infty \le 8 \,C_{{\rm HLS}}\left\Vert|x| \gamma_\epsilon\right\Vert_{4/3}\left\Vert\gamma_\epsilon\right\Vert_{4/3} + 4  + \frac{C_{{\rm HLS}}}{2\pi}\left\Vert\gamma_\epsilon\right\Vert_{4/3}^2\left\Vert|x|\rho \right\Vert_1\; ,
\end{equation*}
and the second part of the result is obtained by using $\left\Vert\gamma_\epsilon\right\Vert_{4/3}=\eps^{-1/2}\left\Vert\gamma\right\Vert_{4/3}$.\lanbox\medskip

Using the regularized Green's function $G_\epsilon$, we introduce the {\em regularized  self-interaction functional} ${\cal W}_\epsilon$:
\begin{equation*}
{\cal W}_\epsilon(\rho) =   \iint_{\R^2 \times \R^2}\rho(x) \,G_\epsilon(x-y)\,\rho(y)\dd x\dd y\ .
\end{equation*}

\begin{lm}[Continuity of the regularized interaction energy]\label{reg}
Let $\rho_1$ and $\rho_2$ be any two densities in $\R^2$
of mass $8\pi$ bounded in $L^1(\R^2,\log(e+|x|^2))$. Then, for all $0< \epsilon \le 1$,
\begin{equation}\label{stconv}
\left|{\cal W}_\epsilon[\rho_1] - {\cal W}_\epsilon[\rho_2]
\right| \le \left[\frac{3}{\pi} +
2C\epsilon^{-2}\right]\|\rho_1\|_{L^1(\R^2,\log(e+|x|^2)\dd
x)}\,\|\rho_1 -\rho_2\|_{L^1(\R^2,\log(e+|x|^2)\dd x)}
\end{equation}
Moreover, let $\{\rho_n\}_{n\in\NN}$ be a sequence of densities
all bounded in $L^1(\R^2,\log(e+|x|^2))$ uniformly in~$n$. If
$\{\rho_n\}_{n\in\NN}$ converges weakly in $L^1(\R^2)$ to $\rho$,
then for each $\epsilon> 0$,
\begin{equation*}
{\cal W}_\epsilon[\rho]\leq \liminf_{n\to\infty}\,{\cal
W}_\epsilon[\rho_n] \ .
\end{equation*}
\end{lm}

\noindent{\bf Proof:} By definition of $G_\epsilon$:
\begin{equation*}
  {\cal W}_\epsilon[\rho] = {\cal W}_0[\rho*\gamma_\eps]=\int_{\R^2} \left| \nabla G* \rho*\gamma_\eps \right|^2(z) \dd z \ge 0\ .
\end{equation*}
Using $\widetilde \rho:=\rho_1-\rho_2$ we write
\begin{equation}\label{recomp5}
{\cal W}_\epsilon[\rho_1] - {\cal W}_\epsilon[\rho_2] = {\cal
W}_\epsilon[\widetilde \rho]+
2\int_{\R^2}\rho_1*G_\epsilon(x)\,\widetilde \rho(x)\dd x \ge
2\int_{\R^2}\left[\rho*G_\epsilon\right]\widetilde \rho\dd x \ .
\end{equation}
Then
combining Lemma~\ref{ptyofG} and~\eqref{recomp5}, we obtain
\begin{align*}
{\cal W}_\epsilon[\rho_1] - {\cal W}_\epsilon[\rho_2]
\ge\,& 2\iint_{\widetilde
\rho>0}\rho(y)\,G_\epsilon(x-y)\,\widetilde \rho(x)\dd y\dd x + 2\iint_{\widetilde \rho\le0}\rho(y)\, G_\epsilon(x-y)
\,\widetilde \rho(x)\dd y\dd x \\
\ge \,& -\frac{1}{2\pi}\iint_{\widetilde
\rho>0}\rho(y)\left[4 + \log\left(e+|x|^2 \right)+
\log\left(e+|y|^2 \right)\right]\widetilde \rho(x)\dd
y\dd x\\
&+2\,\frac{C}{\epsilon^2}\iint_{\widetilde \rho\le0}
\rho(y)\widetilde \rho(x)\dd y\dd x \\
\ge\,& -\left[\frac{3}{\pi} +
2\frac{C}{\epsilon^{2}}\right]\|\rho_1\|_{L^1(\R^2,\log(e+|x|^2)\dd
x)}\,\|\widetilde \rho\|_{L^1(\R^2,\log(e+|x|^2)\dd x)}\ .
\end{align*}
Now swapping the roles of $\rho_1$ and $\rho_2$, we obtain~\eqref{stconv}.

By Lemma~\ref{ptyofG} $\{\rho_n\}_{n\in\NN}$ bounded in
$L^1(\R^2,\log(e+|x|^2)\dd x)$ uniformly in $n$ implies that
$\rho_n*G_\epsilon$ is bounded in $L^\infty(\R^2)$ uniformly in
$n$. Since $\{\rho_n\}_{n\in\NN}$ converges to $\rho$ weakly in
$L^1(\R^2)$, then $\rho_n*(\chi_R \,G_\epsilon)\to \rho*(\chi_R
\,G_\epsilon)$ point-wise for given any cut-off function $\chi_R$
with support in $B(0,R)$ and thus $\rho*G_\epsilon\in
L^\infty(\R^2)$. Therefore, applying~\eqref{recomp5} with
$\rho_1=:\rho$ and $\rho_2:=\rho_n$, we have
$$
\liminf_{n\to\infty}\left({\cal W}_\epsilon[\rho_n] - {\cal
W}_\epsilon[\rho]\right)\ge
\lim_{n\to\infty}2\,\int_{\R^2}\left[\rho*G_\epsilon\right]( \rho_n -
\rho)\dd x = 0\ ,
$$
where we have used the weak convergence on the right hand side.
\lanbox \medskip

We are now  ready to introduce our regularized free energy
functional.

\begin{defi}[Regularized free energy functional]
For all $0<\epsilon\le1$, define
\begin{equation*}
\Fhls^{\epsilon}[\rho] := \int_{\R^2}\rho(x) \log \rho (x)\dd x -
\frac{1}{2}\iint_{\R^2 \times\R^2}
\rho(x)\,G_\epsilon(x-y)\,\rho(y)\,\dd x \dd y
\end{equation*}
on the set of densities $\rho$  of mass $8\pi$ such that
$\rho\in
L^1(\R^2,\log(e+|x|^2)\dd x)$ and $\rho\log \rho$ is integrable.
\end{defi}

The following lemma is an immediate consequence of (\ref{moneps}):

\begin{lm}[Monotonicity in $\epsilon$]\label{monineps}  For all densities  $\rho$  of mass $8\pi$ such that
$\rho\in L^1(\R^2,\log(e+|x|^2)\dd x)$ and $\rho\log \rho$ is
integrable, $\epsilon \mapsto  \Fhls^{\epsilon}[\rho]$ is monotone
decreasing in $\epsilon$.
\end{lm}

Note that by Lemma~\ref{ptyofG} (ii) and (iii),
$\rho\,(G_\epsilon*\rho)$ is integrable for  $\rho\in
L^1(\R^2,\log(e+|x|^2)\dd x)$. Moreover, by Lemma~\ref{ptyofG} (i)
\begin{equation}\label{regcomp45}
\Fhls^{\epsilon}[\rho] \ge \Fhls[\rho]\ .
\end{equation}
In particular, by the sharp log HLS inequality, see Lemma~\ref{lm:loghls}
\begin{equation}\label{regcomp46}
\Fhls^{\epsilon}[\rho] \ge-C(8\pi)=8\pi(-1+\log 8)\ .
\end{equation}
By Lemma~\ref{ptyofG} (iii), we have the upper bound independent
of $\epsilon$:
\begin{equation}\label{regcomp477}
\Fhls^{\epsilon}[\rho] \le  \int_{\R^2}\rho(x) \log \rho (x)\dd x
+ 32\,\pi + 2\|\rho\|_{L^1(\R^2,\log(e+|x|^2)\dd x)}\ .
\end{equation}

\begin{lm}[Error estimate for regularized free energy]\label{erroreg}
For all $\rho \in L^1_+\cap L^{3/2}(\R^2)$ with mass $8\pi$, and
all $\epsilon < (2\sqrt{e})^{-1}$,
\begin{equation*}
\Fhls^{\epsilon}[\rho] - \Fhls[\rho] \le C_\gamma
\,\|\rho\|_{3/2}^{3/2}\,\epsilon\, .
\end{equation*}
\end{lm}

\noindent{\bf Proof:} We use H\"older's inequality and Young's
inequality for convolutions to get
$$
\Fhls^{\epsilon}[\rho] - \Fhls[\rho] \le \|\rho\|^{2}_{4/3}
\|G_{\epsilon}- G\|_{2}\ .
$$
H\"older's inequality gives $ \Fhls^{\epsilon}[\rho] - \Fhls[\rho]
\le \sqrt{8\pi}\,\|\rho\|^{3/2}_{3/2} \|G_{\epsilon}- G\|_{2} $.
Since $G\in L^1_{loc}(\RR^2)$, then $G_\epsilon (x) \to G(x)$ a.e.
in $\R^2$ as $\epsilon\to 0$. By Lemma \ref{ptyofG}, $G_\epsilon
\le G_{\tilde\epsilon}\leq G$ for $0<\tilde\epsilon<\epsilon$.
This implies that $\|G_{\epsilon}- G\|_{2}$ is non decreasing in
$\epsilon$ and it has zero limit as $\epsilon \to 0$ by monotone
convergence theorem. Moreover, the same arguments also show that
$$
\lim_{\tilde\epsilon\to 0} \|G_{\epsilon}-
G_{\tilde\epsilon}\|_{2} = \|G_{\epsilon}- G\|_{2} \, .
$$
On the other hand, using \eqref{deriveps} and \eqref{deriveps2}
for $|x|>2\epsilon$, we get for any fixed $x\in\R^2$
$$
G_{\tilde\epsilon}(x)- G_\epsilon(x) \leq
(\epsilon-\tilde\epsilon) \left|\frac{\partial}{\partial \epsilon}
G_\epsilon (x)\right|  \leq 2 \epsilon^2 \gamma_{\epsilon}(x)
$$
where second-order Taylor expansion at $\epsilon$ of the function
$G_\epsilon (x)$ for $0<\tilde\epsilon<\epsilon$ was used and the
second order term is nonpositive due to \eqref{deriveps2}. Taking
the limit $\tilde\epsilon\to 0$ and integrating, we deduce that
$$
\left(\int_{|x|>2\epsilon} (G(x)- G_\epsilon(x))^2 \dd
x\right)^{1/2} \leq 2 \epsilon\, \|\gamma\|_2.
$$
Since $G_\epsilon \le G$, we can thus directly compute
$$
\int_{|x|<2\epsilon} (G(x)- G_\epsilon(x))^2 \dd x \leq
4\int_{|z|\le 2\epsilon} |G(z)|^2\, \dd z = 16 \pi\,\epsilon^2
\left(\frac12 -\log (2\epsilon) + (\log (2\epsilon))^2\right).
$$
Finally, simple computations show that when $|\log(2\epsilon)| >
1/2$, the term in parentheses on the right in no greater than
$5|\log(2\epsilon)|^2$. Collecting all together leads to the
result with $C_\gamma$ explicitly computable. \lanbox


\subsection{Existence and first properties of the JKO scheme minimizers}
Let ${\mathcal S}$ denote the set :
\begin{equation*}
  {\mathcal S}:=\{\rho \in L^1(\R^2)\,:\,\int_{\R^2}\rho(x)\dd x =M,\;\W_2(\rho,\rhohls) < \infty,\;{\cal E}[\rho]<\infty,\;\int_{\R^2}|x| \rho(x)\dd x < \infty\}\;.
\end{equation*}
By~\eqref{regcomp46},  the functional
$$
\rho \mapsto \frac{\W_2^2(\rho,\rho_0)}{2\tau} + \Fhls^{\epsilon} [\rho],
$$
is bounded from below on $\mathcal S$. The next lemma asserts that
it has minimizers, and begins the task of their analysis. We state
this lemma for a single step since we shall be changing the value
of $\epsilon$ from step to step.
\begin{thm}[Existence of minimizers]\label{lower3}
Let $\lambda>0$, $0< \tau \le 1 $ and $0< \epsilon \le 1$. For all $\rho_0 \in {\mathcal S}$
\begin{equation}\label{fff}
{\rm arg}\min_{\rho\in {\cal S}}\left\{
\frac{\W_2^2(\rho,\rho_0)}{2\tau} + \Fhls^{\epsilon}
[\rho]\right\}
\end{equation}
is not empty, and each minimizer $\rho$ belongs to ${\mathcal S}$.
Moreover, there exists $K_1>0$ depending only on ${\cal
E}[\rho_0]$ and $\int_{\R^2}|x| \rho_0\,\dd x$ such that
$\Fhls^{\epsilon}[\rho] \le K_1$.
\end{thm}
\noindent{\bf Proof:} Let $\{\rho^{(k)}\}_{k\in \N}$ be a minimizing sequence {\it i.e.} such that
$$
\lim_{k\to\infty}\left( \frac{\W_2^2(\rho^{(k)},\rho_0)}{2\tau} +
\Fhls^{\epsilon} [\rho^{(k)}]\right) = \inf_{\rho\in{\cal S}}
\left\{ \frac{\W_2^2(\rho,\rho_0)}{2\tau} + \Fhls^{\epsilon}
[\rho] \right\}\ .
$$
By what we have noted just above, the infimum on the right hand
side is finite. The following observation is the starting point
for obtaining all of the bounds we need: Considering the trial function $\rho =
\rho_0$ itself, one sees that we may suppose
$$
\frac{\W_2^2(\rho^{(k)},\rho_0)}{2\tau} + \Fhls^{\epsilon} [\rho^{(k)}]    \le
\frac{\W_2^2(\rho_0,\rho_0)}{2\tau}  + \Fhls^{\epsilon} [\rho_0] = \Fhls^{\epsilon} [\rho_0]
$$
for all $k$. Consequently, for all $k$,
\begin{equation}\label{tel1}
\Fhls^{\epsilon} [\rho^{(k)}] \le \Fhls^{\epsilon} [\rho_0] \quad\mbox{and}\quad\W_2^2(\rho^{(k)},\rho_0) \le 2\tau\left[ \Fhls^{\epsilon} [\rho_0] - \Fhls^{\epsilon} [\rho^{(k)}]\right]\ .
\end{equation}

We first bound $\W_2(\rho^{(k)},\rhohls)$ uniformly in $k$. Since
$\rho_0\in {\mathcal S}$, \eqref{regcomp477} ensures that
$\Fhls^{\epsilon} [\rho_0] < \infty$, and provides a bound
depending only on ${\cal E}[\rho_0]$ and $\int_{\R^2}|x| \rho_0\dd
x$. Then~\eqref{regcomp46} provides a universal lower bound on
$\Fhls^{\epsilon} [\rho^{(k)}]$, and thus by~\eqref{tel1}, there
is a finite constant $K_1$ depending only on ${\cal E}[\rho_0]$,
and $\int_{\R^2}|x| \rho_0\dd x$ such that for all $k$,
\begin{equation}\label{tel2}
\Fhls^{\epsilon} [\rho^{(k)}]\leq  K_1
\quad\mbox{and}\quad
\W_2^2(\rho^{(k)},\rho_0) \leq K_1\ .
\end{equation}
In particular, by the triangle inequality, for all $k$,
$\W_2(\rho^{(k)},\rhohls) \leq \sqrt{K_1} +  \W_2(\rho_0,\rhohls) < \infty$.

We next bound the first moments of $\rho^{(k)}$ uniformly in $k$.
 Let $\nabla \varphi$ be the
optimal transportation plan  $\nabla \varphi\#\rho^{(k)} =\rhohls$.
Then since  $|x| \leq |x- \nabla \varphi(x)| + |\nabla \varphi(x)|$ for all $x$,  integrating against $\rho^{(k)}$ and using the
Cauchy-Schwarz inequality yields
$$
\int_{\R^2}|x|\,\rho^{(k)}(x)\dd x \leq
\sqrt{8\,\pi}\,\W_2(\rho^{(k)},\rhohls) +
\int_{\R^2}|x|\,\rhohls(x)\dd x
 \leq \sqrt{8\,\pi}\,[\sqrt{K_1} +  \W_2(\rho_0,\rhohls) ] + \int_{\R^2}|x|\,\rhohls(x)\dd x \ .
$$
The right hand side is finite and independent of $k$.

We next bound ${\mathcal E}[\rho^{(k)}]$.  By part {\it (ii)} of  Lemma~\ref{ptyofG},  there is a constant $C$ such that
\begin{align*}
\int_{\R^2}\rho^{(k)}(x)\log \rho^{(k)}(x)\dd x \le \Fhls^{\epsilon}[\rho^{(k)}] +\frac{C\,(8\pi)^2}{\epsilon^{2}}\le \Fhls^{\epsilon}[\rho_0]
+ \frac{C\,(8\pi)^2}{\epsilon^{2}} \le K_1 + \frac{C\,(8\pi)^2}{\epsilon^{2}}\ ,
\end{align*}
where we have used (\ref{tel2}) once more. Again the right side is finite and independent of $k$.

The last two uniform bounds show that $\{\rho^{(k)}\}_{k\in\NN}$
is uniformly integrable. Hence, by the Dunford-Pettis theorem,
there exist a weakly in $L^1$ convergent sub-sequence whose limit
we shall denote by $\rho$.

By a standard weak lower semicontinuity argument (see e.g.
\cite{AGS} for the weak lower semicontinuity of $\W_2^2$), $\rho$
satisfies each of the three bound that we have proved uniformly
for $\{\rho^{(k)}\}_{k\in\NN}$, and thus $\rho\in{\mathcal S}$.

It remains to prove that the functional $\Fhls^{\epsilon}$ is
lower semi-continuous on $L^1(\R^2)$. For the entropy part, this
is standard. For the self interaction part, this follows from
Lemma~\ref{reg}. So that
\begin{equation}
  \label{eq:newnumber}
  \Fhls^{\epsilon}[\rho] \le \liminf \Fhls^{\epsilon}[\rho^{(k)}] <K_1\;.
\end{equation}
Finally, the weak limit $\rho$
is a minimizer.\lanbox
\begin{prop}[Strict positivity of the minimizers]
Let $\rho_0$ satisfies the conditions of
{\rm Theorem~\ref{lower3}}. Then any minimizer $\rho$
of~\eqref{fff} is uniformly bounded below on compact sets,
\emph{i.e.} for all $R >0$, there exists $s>0$ such that
\begin{equation}\label{dp4}
\rho(x) \ge s \, \mbox{ almost everywhere in } \, D_R:=\{ x : |x|
\le R \}\, .
\end{equation}
Moreover, $s$ does not depend on the chosen minimizer of Problem~\eqref{fff} in case of non-uniqueness.
\end{prop}
\begin{remark}[Idea of the proof]
Let us try to quantify this simple statement: Entropy abhors a
vacuum. The functional derivative of ${\cal E}[\rho]$ is
$\log\rho$. On any set where $\rho$ is very close to zero,
$\log\rho$ is very negative, and we can lower ${\cal E}[\rho]$ by
transporting some mass from where $\rho$ is relatively large to
this spot where it is very small. This will lower the entropy by a
very large multiple of the transported mass. On the other hand, if
we do not have to transport the mass too far, the effects on
$\W_2^2(\rho,\rho_0)$ and ${\cal W}_\epsilon(\rho)$ will be
relatively small.
\end{remark}

\noindent{\bf Proof:}
For any $s>0$, let
$$
A_R(s) := \{ x\in D_R : \rho(x) \ge s \}\quad\mbox{and}\quad
C_R:=\int_{A_R(2/R^2)}\rho(x)\dd x\, .
$$
For any $s>0$, let $\alpha_R(s) := \{ x\in D_R : \rho(x)
\le s \}$, and let $|\alpha_R(s)|$ denote its Lebesgue measure.
By Theorem~\ref{lower3}, $|x|\rho(x)$ is integrable, and hence
$$
\int_{D_R}\rho(x)\dd x =\int_{\R^2}\rho(x)\dd x- \int_{|x|> R}\rho(x)\dd x \ge
8\,\pi-\int_{\R^2} \frac{|x|}{R}\rho(x)\dd x \geq 4\,\pi\ ,
$$
as long as $4\pi R \ge \int_{\R^2}|x|\rho(x)\dd x$.

If $|\alpha_R(s)| = 0$ for some $s>0$, there is nothing to prove:
$\rho$ is bounded below uniformly by $s$ on $D_R$. Therefore,
suppose that $|\alpha_R(s)| > 0$ for all $s>0$. Pick some small
positive numbers $\delta$ and $s$, and define a new density
$\widetilde\rho$ by transporting a mass
$\delta\,C_R\,|\alpha_R(s)|$ from $A_R(2/R^2)$ to $\alpha_R(s)$,
distributing it {\em uniformly} there, which raises the density
there by $\delta\,C_R$. In formulas, choose $s< 2/R^2$ to have
$\alpha_R(s)\cap A_R(2/R^2)=\emptyset$, and define a new density
$\widetilde \rho$ by
\begin{equation*}
\widetilde \rho(x) =
\begin{cases}
 (1-\delta\,|\alpha_R(s)|)\,\rho(x) & x\in A_R(2/R^2)\;,\\
\rho(x) + \delta\,C_R &x \in \alpha_R(s)\;,\\
\rho(x)& {\rm otherwise}\, .
 \end{cases}
\end{equation*}
In order to ensure positivity, we have to impose
$\delta|\alpha_R(s)|\leq \delta \,\pi\,R^2 \leq 1/2$. In this way,
it is easy to check that $\widetilde \rho$ is a density.

Note that $\|\widetilde \rho  - \rho\|_1 \le 2\,\delta
\,|\alpha_R(s)|C_R\,$, and since all the modifications take place
on $D_R$,
$$
\|\widetilde \rho  - \rho\|_{L^1(\R^2,\log(e+|x|^2)\dd x)} \le
\log(e+R^2)\,2\,\delta\, |\alpha_R(s)|\,C_R\ .
$$
It now follows from the bounds on $\rho$ derived Theorem~\ref{lower3}  and from  Lemma~\ref{reg} that there is
a constant $K$ depending only on $R$, $\epsilon$, ${\cal
E}[\rho_0]$ and $\Fcfd[\rho_0]$ such that
\begin{equation}\label{tri18}
{\cal W}_\epsilon[\widetilde\rho] \le  {\cal W}_\epsilon[\rho] +
\delta \,|\alpha_R(s)|\,K\ .
\end{equation}

Using Taylor's expansion of $x \mapsto x\log x$, that $\log x$ is
increasing and assuming $s\leq \delta \,C_R$, we obtain
\begin{align}
\int_{\R^2}[\widetilde \rho \log \widetilde  \rho -\rho \log  \rho] \dd x &=
\int_{A_R(2/R^2)}\left\{(1-\delta\,|\alpha_R(s)|)\,\rho\,\log\left[(1-\delta\,|\alpha_R(s)|)\,\rho\right] - \,\rho\log\rho\right\}\dd x\nonumber\\
&\quad+ \int_{\alpha_R(s)}\left[(\rho+ \delta\,C_R) \log\left(\rho+ \delta\,C_R\right) -  \rho \log  \rho\right]\dd x\nonumber\\
&\le  -
\delta\,|\alpha_R(s)|\int_{A_R(2/R^2)}\!\!\!\rho\,\log\left[(1-\delta\,|\alpha_R(s)|)\,\rho\right]\dd
x + \delta\,|\alpha_R(s)|\,C_R  \log(2\delta\,C_R)\nonumber\\
&\le  \delta\,|\alpha_R(s)|\, C_R \left[-
\log\left(\frac1{R^2}\right) + \log(2\delta\,C_R)\right]\ ,
\label{tri17}
\end{align}
where $\delta\,|\alpha_R(s)| \leq 1/2$ and $x\in A_R(2/R^2)$ were
used in the last estimate.

To estimate the difference $\W_2^2(\widetilde\rho,\rho_0) -
\W_2^2(\rho,\rho_0)$, let $\Pi$ denote the optimal coupling of
$\rho$ and $\rho_0$, and use it to define a non-optimal coupling
$\widetilde \Pi$ of $\widetilde\rho$ and $\rho_0$. To do this, let
$\mu$ be the measure supported on $A_R(2/R^2)$ with density
$\rho$, and hence total mass $C_R$. Let $\nu$ be the $\dd x$--uniform distribution on $\alpha_R(s)$ with
total mass $C_R$. Let $\nabla\psi$ be the optimal transportation
plan with $\nabla \psi \# \mu = \nu$, and define the map $T:\R^2
\to \R^2$ by
$$
T(x) =
\begin{cases}
\nabla \psi(x) & x\in A_R(2/R^2)\;,\\
 x & {\rm otherwise}\, .
 \end{cases}
$$
Then $\widetilde \Pi$, given by $\widetilde \Pi = (1 -
\delta|\alpha_R(s)|) \Pi + \delta|\alpha_R(s)|(T\otimes {\rm
Id})\# \Pi$ is a coupling of $\widetilde\rho$ and $\rho_0$, and
hence
\begin{align*}
 \W_2^2(\widetilde\rho,\rho_0) &\le \iint_{\R^2\times\R^2} |x-y|^2 \dd \widetilde \Pi(x,y) \\
 &= (1 - \delta|\alpha_R(s)|)\,\W_2^2(\rho,\rho_0) + \delta|\alpha_R(s)|
 \iint_{\R^2\times\R^2} |T(x)-y|^2 \dd \widetilde \Pi(x,y)\, .
\end{align*}
Then, since $|T(x) - y|^2 \le 2|T(x) - x|^2 + 2|x-y|^2$, and
$|T(x) - x| \le 2R$, since all of the transportation induced by
$\nabla\psi$ takes place inside $D_R$, it follows that
$$
\W_2^2(\widetilde\rho,\rho_0) \le
(1+\delta|\alpha_R(s)|)\,\W_2^2(\rho,\rho_0) +
\delta|\alpha_R(s)|(8\pi)^2\, 8R^2\ .
$$
By  bounds on $\rho$ derived in the proof of Theorem~\ref{lower3}, there is a constant
$\tilde K$ depending only on ${\cal E}[\rho_0]$ and
$\Fcfd[\rho_0]$ such that  $\W_2^2(\rho,\rho_0) \le \tilde K\tau$.
Finally then, there is a constant depending only on $R$, $\tau$,
${\cal E}[\rho_0]$ and $\Fcfd[\rho_0]$ such that
\begin{equation}\label{tri19}
\W_2^2(\widetilde \rho,\rho_0) \le  \W_2^2(\rho,\rho_0) +
\delta|\alpha_R(s)| K\ .
\end{equation}
Combining~\eqref{tri17},~\eqref{tri18} and~\eqref{tri19} yields
\begin{align*}
\frac{\W_2^2(\widetilde \rho,\rho_0)}{2\tau} + \Fhls^{\epsilon}
[\rho] \le  \frac{\W_2^2(\rho,\rho_0)}{2\tau} + \Fhls^{\epsilon}
[\rho] + \delta\,|\alpha_R(s)|\, C_R \left[-
\log\left(\frac1{R^2}\right) + \log(2\delta C_R) + K'\right]\ ,
\end{align*}
with a given constant $K'$. If $|\alpha_R(s)|>0$ for all $s>0$,
then choosing $\delta$ small enough such that
$$
- \log\left(\frac1{R^2}\right) + \log(2\delta\,C_R) + K' < 0
$$
contradicts the optimality of $\rho$. For instance, choosing $s_R=
\delta/C_R$, the procedure described above can be carried out, and
we conclude that $\rho$ is bounded below by $s_R$ on $D_R$. This
proves~\eqref{dp4}. \lanbox

\

We now continue the analysis of the minimizers $\rho$
begun in Theorem~\ref{lower3}. We obtained $\rho\in {\mathcal S}$
and the lower bound~\eqref{dp4} directly from the variational principle,
but to proceed, we need
the Euler-Lagrange equation for the variational problem~\eqref{fff}.

By the Brenier-McCann Theorem, there is a a lower
semi-continuous  convex function $\varphi$ in $\R^2$ such that
$\nabla \varphi \#  \rho = \rho_0$, and $\nabla \varphi$ is
uniquely determined on the support of $\rho$, which is all
$\R^2$ by~\eqref{dp4} .  The
Euler-Lagrange equation for~\eqref{fff} relates $\rho$, $\rho_0$ and
$\nabla\varphi$:

\begin{lm}[Euler-Lagrange equation]\label{lm:eulerlagrange}
Let $\rho_0$ satisfy the conditions of Theorem~{\rm\ref{lower3}}
and $\rho$ be any minimizer for the variational problem
in~\eqref{fff}, and let $\nabla\varphi$ be the unique gradient of
a lower semi-continuous convex function such that  $\nabla \varphi
\#  \rho = \rho_0$. Then the distributional gradient of $\rho$
satisfies
\begin{equation}\label{disgrad}
-\nabla \rho + \rho \nabla c_\epsilon = \frac{{\rm id}-\nabla
\varphi }{\tau} \rho
\end{equation}
where $c_\epsilon  = G_\epsilon*\rho$. In particular, since $c_\epsilon$ is
differentiable everywhere, and $\varphi$ is differentiable almost
everywhere, $\rho$ is differentiable almost everywhere.
\end{lm}

The proof follows exactly the original procedure in \cite{JKO},
see also \cite[Theorem 3.4]{BCC}, and we skip it here for the sake
of brevity. The interested reader can see the details in the
preprint version of this paper.

\begin{lm}[Qualitative regularity estimates]\label{feul}
Let $\rho_0$ satisfy the conditions  of Theorem~{\rm\ref{lower3}},
and let $\rho$ be any minimizer for the variational problem
in~\eqref{fff}. Then $\sqrt{\rho}$  has a square integrable
distributional gradient, and for any $1 < p < \infty$, $\rho^p$ is
integrable.
\end{lm}\medskip

\noindent{\bf Proof:} By the positivity of $\rho$, see Proposition~\ref{lower3}, we can divide both sides
of~\eqref{disgrad} by $\sqrt{\rho}$, to obtain
$$
2\nabla \sqrt{\rho} = \left(\nabla c_\epsilon - \frac{x-\nabla \varphi}{\tau}\right)\sqrt{\rho}\ ,
$$
where $\nabla\varphi$ is such that $\nabla \varphi \#  \rho =
\rho_0$. By the triangle inequality,
\begin{align}\label{nscrummy2}
2\|\nabla \sqrt{\rho}\|_2 &\le  \left(\int_{\R^2}\left|\nabla c_\epsilon(x)\right|^2\rho(x)\dd x\right)^{1/2}
+\frac1{\tau} \left(\int_{\R^2}|x-\nabla \varphi(x)|^2\rho(x)\dd x\right)^{1/2}\nonumber\\
&\le  \left(\int_{\R^2}\left|\nabla
c_\epsilon(x)\right|^2\rho(x)\dd x\right)^{1/2} + \frac1{\tau}
\W_2(\rho,\rho_0)\ .
\end{align}
By Lemma~\ref{regconc}, $\|\nabla c_\epsilon\|_\infty$ is
uniformly bounded, and so is the first term of~\eqref{nscrummy2}.
This proves that $\sqrt{\rho}$ has a square integrable
distributional gradient. The integrability of $\rho^p$ is then a
consequence of  the following classical version of the GNS
inequality valid for functions on $\R^2$ with $p \in [2,\infty)$
\begin{equation*}
  \int_{\R^2}|v|^{p}\dd x \le D_p\left[\int_{\R^2}|\nabla v|^{2}\dd x\right]^{p/2-1}\int_{\R^2}|v|^{2}\dd x
\end{equation*}
applied to $v=\sqrt \rho$.\lanbox

\begin{remark}
Although the bounds in Lemma~{\rm \ref{feul}} are not
quantitative, and would certainly be~$\epsilon$ dependent if we
were to extract quantitative bounds, we shall use them only to
justify certain integrations by parts, and otherwise show that
substraction of infinities does not invalidate computations that
follow. Thus, these qualitative estimates are all we require
concerning $\nabla \sqrt{\rho}$ and $\rho$. However, they are
absolutely crucial for their purpose, and their necessity is the
main reason we have had to introduce the regularized Green's
function $G_\epsilon$, and along with it, the regularized chemical
attractant. Without the regularization, we would only know that
$2\nabla \sqrt{\rho} - \nabla c\,\sqrt{\rho}$ was square
integrable -- but the possible cancelation effects would not allow
us to conclude that  $\nabla \sqrt{\rho}$ was square integrable.
\end{remark}


\subsection{A discrete form of the entropy-entropy dissipation inequality}

Our main goal in this subsection is to prove a discrete version of
the entropy-entropy dissipation inequality (\ref{enendi}).  The
key idea is to use the $\kappa_\lambda$ displacement convexity of
$\Fcfd$ and the ``above the tangent'' inequality for convex
functions as follows:  For given initial density $\rho_0$, let
$\rho$ be any minimizer for the variational problem
in~\eqref{fff}.  Let $u_t$, $0\leq t \leq 1$ denote the
displacement interpolation between $\rho$ and $\rho_0$ starting at
$\rho$ and ending at $\rho_0$. Then $u_0 = \rho$ and $u_1 =
\rho_0$. Since $\Fcfd$ is displacement convex, the ``above the
tangent'' inequality for convex functions says that
$$
\Fcfd[\rho]  + \frac{{\rm d}}{{\rm d}t}\Fcfd[u_t]\bigg|_{t=0} +
\kappa_\lambda\,\W_2^2(\rho,\rho_0)
 \leq \Fcfd[\rho_0]\ .
$$
A formal computation of the second term on the left would give,
for $\epsilon = 0$,
$$
\frac{{\rm d}}{{\rm d}t}\Fcfd[u_t]\bigg|_{t=0}  = {\mathcal
D}[\rho]\ .
$$
Indeed, assuming
Lemma~\ref{discobnd} holds for $\delta=\epsilon=0$ applied to
$u_0=\rho$ and $u_1=\rho_0$, we get
\begin{equation*}
   \Fcfd[\rho] \le \Fcfd[\rho_0] - \frac{1}{2}\int_{\R^2}\left[  \kappa_\lambda x
   +\frac{\nabla \rho}{\rho^{3/2}}\right]\cdot  (\nabla\varphi(x) -x)\,\rho\dd x  - \kappa_\lambda\,\W^2_2(\rho,\rho_0)\ .
\end{equation*}
Using~\eqref{disgrad}, {\it i.e.
$(\nabla\varphi(x)-x)\rho = \tau \left(\nabla
\rho -\rho\nabla c\right)$} and expanding, we can rewrite this as
\begin{align*}
\Fcfd[\rho] \le \, &\Fcfd[\rho_0] - \frac\tau2 \left[ \int_{\R^2}
\frac{\left|\nabla \rho\right|^2}{\rho^{3/2}}\dd x
-\int_{\R^2}\frac{\nabla c \cdot\nabla \rho}{\sqrt \rho}\dd x
-\kappa_\lambda\int_{\R^2} x \cdot\nabla c\,\rho\dd x
+\kappa_\lambda\int_{\R^2} x \cdot  \nabla \rho\dd x \right] \\&
\qquad \qquad- \kappa_\lambda\,\W_2^2(\rho,\rho_0)
\\
:= &\Fcfd[\rho_0] - \frac\tau2 \left[ {\rm (I)} + {\rm (II)} +{\rm
(III)} + {\rm (IV)} \right] - \kappa_\lambda\,\W_2^2(\rho,\rho_0)\ .
\end{align*}
Using $-\Delta c=\rho$ we have
\begin{equation*}
  {\rm (II)}=-2 \int_{\R^2}\nabla c \cdot\nabla \sqrt\rho\dd x = -2 \int_{\R^2}\rho^{3/2}\dd x\ .
\end{equation*}
Using the symmetrization argument we obtain
\begin{equation*}
  {\rm (III)}=\kappa_\lambda\frac1{4\pi}\left(\int_{\R^2} \rho\dd x\right)^2=16\pi\kappa_\lambda\ .
\end{equation*}
And by integration by parts ${\rm (IV)}=-16\pi\kappa_\lambda$,
resulting into $\Fcfd[\rho]  \le  \Fcfd[\rho_0] - \tau {\cal
D}[\rho]  - \kappa_\lambda\,\W_2^2(\rho,\rho_0)$.

However, to do the calculation in a rigorous manner we must take into account that
$\epsilon >0$, and we must use the regularized entropy functional $\Fcfdr$.   Before proceding with this, we point out that
no such estimate can be given for $\Fhls$ since this functional is not displacement convex.

\begin{lm}[Convexity estimates at the regularized level]\label{displconvexitydelta}
Let $\rho_0$ satisfy the conditions  of Theorem~{\rm
\ref{lower3}}, and let $\rho$ be any minimizer for the variational
problem in~\eqref{fff}, then
\begin{align*}
\Fcfdr[\rho] \le & \, \Fcfdr[\rho_0]-\frac{\tau}{2}\int_{\R^2}
\frac{\left|\nabla \rho\right|^2}{(\rho + \delta)^{3/2}} \dd
x+\tau \int_{\R^2}\rho^{3/2}\dd
x + 16\,\pi\,\sqrt{\delta}\tau+ 16\pi\sqrt{\delta} (J_\gamma + 1 + \hat C_\epsilon)\tau \nonumber\\
&-16\pi\kappa_\lambda \tau+ 2C_\epsilon \,
\|\,|2f_\delta'-\kappa_\lambda|\,(1+|x|)\rho \|_1 \tau+
\frac{16\pi\,J_\gamma}{\sqrt{2\lambda}}\tau + 2\tau \int_{\R^2}
\nabla\cdot[x\,f_\delta'(|x|^2)]\rho(x)\dd x\nonumber\\&
-K_\delta(\rho_0,\rho)\ ,
\end{align*}
and
\begin{align*}
\Fcfdr[\rho] \le &\, \Fcfdr[\rho_0] -\frac{\tau}{2}\int_{\R^2}
\frac{\left|\nabla \rho\right|^2}{(\rho + \delta)^{3/2}} \dd
x+\tau \int_{\R^2}\rho^{3/2}\dd
x + 16\,\pi\,\sqrt{\delta}\tau+ 16\pi\sqrt{\delta} (J_\gamma + 1 + \hat C_\epsilon)\tau \nonumber\\
&-16\pi\kappa_\lambda \tau+ 2C_\epsilon \,
\|\,|2f_\delta'-\kappa_\lambda|\,(1+|x|)\rho \|_1 \tau+
\frac{32\pi\,C_{{\rm HLS}}}{\sqrt{2\lambda}} \, \sqrt\epsilon\,
\||x|\,\gamma\|_{4/3} \,\| \rho \|_{4/3} \tau\nonumber\\& + 2\tau
\int_{\R^2} \nabla\cdot[x\,f_\delta'(|x|^2)]\rho(x)\dd
x-K_\delta(\rho_0,\rho)\ ,
\end{align*}
where $K_\delta$ is defined in {\rm Proposition~\ref{entap}} and
the constants $J_\gamma$, $\hat C_\epsilon$, and $C_\epsilon$ are
explicit constants.
\end{lm}\medskip

\noindent{\bf Proof:} This is an elaborate calculation in which a
number of integrations by parts operations must be carefully
examined for boundary behavior. It is relegated to the Appendix.
\lanbox

\

As a consequence of this lemma, letting $\delta$ go to 0, we obtain the
following  result concerning the dissipation of
$\Fcfd$ in one discrete time step.

\begin{cl}[Convexity estimates]\label{cordis}
Let $\rho_0$ satisfy the conditions of Theo\-rem~{\rm
\ref{lower3}}. If $\rho$ is any minimizer for the variational
problem in~\eqref{fff} then
\begin{equation}\label{eq:witouteps}
\Fcfd[\rho] \le \Fcfd[\rho_0]-\tau\,{\cal D}[\rho] +\tau\,\AA\,\|\gamma\|_{4/3} - \kappa_\lambda
\W_2^2(\rho,\rho_0)\ ,
\end{equation}
and
\begin{equation}\label{eq:witheps}
\Fcfd[\rho] \leq \Fcfd[\rho_0]-\tau \,{\cal D}[\rho] +
\tau\sqrt\epsilon\,\AA\,\| \rho \|_{4/3} -
\kappa_\lambda \W_2^2(\rho,\rho_0)\ .
\end{equation}
where $\AA:=32\pi(2\lambda)^{-1/2} C_{{\rm HLS}}\,
\||x|\,\gamma\|_{4/3}$.
\end{cl}

\noindent{\bf Proof:} Let us first observe that
$$
\int_{\R^2}\nabla\cdot[x\,f_\delta'(|x|^2)] \rho(x)\dd x =
\int_{\R^2}\left[2f_\delta'(|x|^2)+2|x|^2f_\delta''(|x|^2)\right]
\rho(x)\dd x .
$$
Let us recall from the proof of Proposition~\ref{entap} that
$2f_\delta'(s)\nearrow \kappa_\lambda$ and $f_\delta''(s)\to 0$ as
$\delta\to 0$ for all $s\geq 0$. Moreover, we have that
$sf_\delta''(s)$ is a bounded function uniformly in $\delta$ from
\eqref{passingtolimitdelta}. These properties together with the
dominated convergence theorem leads easily to
$$
\|\,|2f_\delta'-\kappa_\lambda|\,(1+|x|)\rho \|_1\to 0 \qquad
\mbox{and} \qquad \int_{\R^2}\nabla\cdot[x\,f_\delta'(|x|^2)]
\rho(x)\dd x \to 8\kappa_\lambda\pi
$$
as $\delta \to 0$, since $(1+|x|)\rho\in L^1(\R^2)$. By monotone
convergence theorem, we obtain
$$
\lim_{\delta\to 0} \int_{\R^2} \frac{\left|\nabla
\rho\right|^2}{(\rho + \delta)^{3/2}} \dd x = \int_{\R^2}
\frac{\left|\nabla \rho\right|^2}{\rho^{3/2}}\dd x\,.
$$

Putting together all these facts and Proposition~\ref{entap}, we
can pass to the limit as $\delta\to 0$ in Lemma
\ref{displconvexitydelta} to get the desired estimates
\eqref{eq:witheps} and \eqref{eq:witouteps}. \lanbox

\subsection{One-step estimates}

Neither of the one step dissipation estimates that we have so far,
namely \eqref{eq:witheps} and \eqref{eq:witouteps}, are exactly
what we need. The problem is the term $\tau A\|\gamma\|_{4/3}$ in
the first of these, and the term $\tau A\|\rho\|_{4/3}$ in the
second of these. These terms might be large compared to the other
terms so that these estimates might even give only ``negative
dissipation''.

In the first main result of this subsection, we use one and then
the other of these inequalities in combination with the controlled
concentration inequality of Theorem~\ref{both2} to produce the
kind of dissipation estimate that we really want.  In the second
main result, we show that $L^p$ norms of the densities are
essentially propagated along each step of the discrete variational
scheme. Again, Theorem~\ref{both2} plays a crucial role in both
proofs.

\begin{thm}[One-step theorem]\label{thmstep}
Let $\rho_0$ satisfy the conditions  of Theorem~{\rm
\ref{lower3}}, $\rho$ be any minimizer for the variational
problem in~\eqref{fff}, and choose any $\crho$ such that
\begin{equation}\label{resini}
\Fcfd[\rho_0]  < \crho \, .
\end{equation}
Define $Q_0 > 0$, $\tau_0^{\star}>0$ by
\begin{equation}\label{taude}
Q_0:=\crho-  \Fcfd[\rho_0]\quad \mbox{and} \quad
  \tau_0^\star := \min\left\{\frac{Q_0}{2\,\AA\|\gamma\|_{4/3}}\ ,\ 1\right\}\ ,
\end{equation}
where $\AA$ is the constant given in {\rm Corollary~\ref{cordis}}. Finally, given $Q_0$ and $0<\tau\le\tau_0^\star$, and also any positive
integer $\ell$, let $\epsilon_\ell$ be given by
\begin{equation*}
 \tau^{1/3}\sqrt{\epsilon_\ell}\left[8\,\pi^{1/3}\,\AA\, \gamma_2^{-2/3} \left(\pi \,\crho+\tau_0^\star C_{{\rm CCD}} \right)^{2/3}\right] = \frac{Q_0}{4}\tau^2 2^{-\ell}\ .
\end{equation*}
Then for all $\tau\le \tau_0^\star$ and all $\epsilon =\epsilon_\ell$, $\rho$ satisfies
\begin{equation}\label{boundforboth}
\Fhls[\rho] < +\infty\;, \quad \Fcfd[\rho]  < \crho
\end{equation}
and
\begin{equation}\label{taude3}
\Fcfd[\rho] -\Fcfd[\rho_0] \le -\,\tau{\mathcal D}[\rho] +
\frac{Q_0}{4} \tau^22^{-\ell} -  \kappa_\lambda \W_2^2(\rho_0,\rho)\ .
\end{equation}
\end{thm}

\noindent{\bf Proof:}
By \eqref{eq:witouteps}, our choice of $\tau$ and $Q_0$ in~\eqref{taude} implies that
\begin{equation}\label{fonda}
  \Fcfd[\rho] \le \Fcfd[\rho_0] -\tau{\mathcal D}[\rho]+ \frac{Q_0}{2} =
  \crho - Q_0 -\tau{\mathcal D}[\rho]+ \frac{Q_0}{2}\le \crho-\tau{\mathcal D}[\rho] \ .
\end{equation}
On one hand, the GNS
inequality, see Lemma~\ref{spgn}, implies ${\cal D}[\rho] \ge 0$ so that~\eqref{fonda} implies that $\rho$ also satisfies~\eqref{resini}.
On the other hand, since $\Fcfd[\rho]$ cannot be negative it implies
\begin{equation*}
{\mathcal D}[\rho] \le \frac{\crho}{\tau}\ .
\end{equation*}
Moreover, by the monotonicity of $\eps \mapsto \Fhls^\eps[\rho]$,
Lemma \ref{monineps}, and~\eqref{eq:newnumber}, $\Fhls[\rho] \le
\Fhls^\eps[\rho] < +\infty$. We can thus apply the concentration
controlled inequality, Theorem~\ref{both2} which implies
$$
 \int_{\R^2} \left|\nabla \rho^{1/4}\right|^2 \dd x \le
 \frac{1}{\gamma_2}\left[\pi{\mathcal D}[\rho] + C_{{\rm CCD}}\right] \le
 \frac{1}{\tau} \frac{1}{\gamma_2}\left[\pi \crho + \tau_0^\star C_{{\rm CCD}}\right] \ .
$$
By the GNS inequality of Lemma
\ref{spgn}, we have
\begin{equation}\label{l32}
\int_{\R^2}\rho^{3/2} \dd x \le 8 \int_{\R^2} \left|\nabla
\rho^{1/4}\right|^2 \dd x \le  \frac{1}{\tau}
\frac8{\gamma_2}\left[\pi \crho+\tau_0^\star C_{{\rm CCD}} \right]
:= \frac{C_3}{\tau} \ .
\end{equation}
Next, by H\"older's inequality,
\begin{equation*}
\int_{\R^2}\rho^{4/3} \dd x = \int_{\R^2}\rho^{1/3} \rho \dd x \le
(8\pi)^{1/3}\left(\int_{\R^2}\rho^{3/2} \dd x\right)^{2/3}\le
(8\pi)^{1/3}\left( \frac{C_3}{\tau} \right)^{2/3}\ .
\end{equation*}
Now using this bound in~\eqref{eq:witheps}, we obtain
\begin{equation*}
\Fcfd[\rho] -\Fcfd[\rho_0] \le -\tau{\mathcal D}[\rho]
+\tau^{1/3}\sqrt{\epsilon}\,\left[\AA\, (8\pi)^{1/3} C_3^{2/3}\right]
-\kappa_\lambda \W_2^2(\rho_0,\rho).
\end{equation*}
We thus obtain the stated result by choosing $\epsilon=\epsilon_\ell$
for any positive integer $\ell$. \lanbox

\begin{lm}[Propagation of the $L^p$-norm]\label{propagate}
Let $\rho_0$ satisfy the conditions of Theorem~{\rm \ref{lower3}}.
Assume additionally that $\rho_0\in L^p(\R^2)$, $2\leq p< \infty$,
and let $\rho$ be any minimizer for the variational problem
in~\eqref{fff}, then there exists $K_0>0$ which only depends on
${\displaystyle \int_{\R^2}\rho|\log\rho|\dd x}$ such that for all
$K\ge K_0$
\begin{equation*}
 \int_{\R^2} (\rho-K)_+^p \dd x \le \int_{\R^2} (\rho_0-K)_+^p\dd x +\tau  A_1 + \tau A_2 {\mathcal
 D}[\rho]\;\,,
\end{equation*}
where $A_1$ and $A_2$ are universal positive constants depending
on $K$.
\end{lm}

\noindent{\bf Proof:} The displacement convexity of the functional
$$
{\displaystyle \rho\mapsto
\int_{\R^2} (\rho-K)_+^p \dd x}
$$
with $2\leq p <\infty$ and $K>0$, is easy to check using McCann's
criterion~\eqref{i:func}. The Euler-Lagrange equation
of the variational scheme $\left(x-\nabla \varphi\right) \rho =
-{\tau}\nabla \rho + {\tau}\rho \nabla c_\epsilon$ together with
the standard first-order displacement convexity characterization
\cite{Villani,AGS} imply
\begin{align}
\int_{\R^2} (\rho-K)_+^p \dd x &- \int_{\R^2} (\rho_0-K)_+^p\dd x  \le
-p\,\int_{\R^2}  \nabla\left[(\rho-K)_+^{p-1}\right] \left(\nabla \varphi-x\right)\rho\dd x\notag\\
\le\, &-(p-1) p \tau\int_{\R^2} \left|\nabla(\rho-K)_+\right|^2(\rho-K)_+^{p-2}\dd x\notag\\
& +(p-1) \tau\int_{\R^2} \nabla\left[(\rho-K)_+^{p}\right]\nabla c_\epsilon \dd x\notag +p\tau K\int_{\R^2} \nabla\left[(\rho-K)_+^{p-1}\right] \nabla c_\epsilon \dd x\notag\\
\le\, &-\frac{4(p-1)}{p}\,\tau\int_{\R^2} \left|\nabla\left[(\rho-K)_+^{p/2}\right]\right|^2\dd x + (p-1) \tau\int_{\R^2} (\rho-K)_+^{p}(-\Delta c_\epsilon) \dd x\notag\\
&+p\tau K\int_{\R^2} (\rho-K)_+^{p-1}(-\Delta c_\epsilon) \dd x :=
\tau (I_1 + I_2 + I_3)\;. \label{estt2}
\end{align}
The last two integration by parts have to be justified for any
given $\epsilon$ working as in the proof of
Lemma~\ref{displconvexitydelta} in the Appendix. Integrating by
parts on the ball of radius $R$, we obtain for any $k \in
\{p,p-1\}$
\begin{align*}
  \int_{|x|\le R} \nabla \left[(\rho-K)_+^{k}\right]\nabla c_\epsilon \dd x &= \int_{|x|\le R} (\rho-K)_+^{k}(-\Delta c_\epsilon) \dd x + \int_{|x|= R} (\rho-K)_+^{k}\nabla c_\epsilon \cdot n \dd \sigma\\
&\le  \int_{|x|\le R} (\rho-K)_+^{k}(-\Delta c_\epsilon) \dd x +
\|\nabla c_\epsilon\|_\infty\int_{|x|= R} \rho^{k}  \dd \sigma\,.
\end{align*}
It is enough to show by dominated convergence theorem that there
exists a sequence of radii $\{R_j\}_{j\in\N}$ such that the
boundary terms tend to zero as $j\to\infty$. Due to
Lemma~\ref{feul} with $p\geq 2$, for any given natural $N>1$, we
can write that
\begin{equation*}
 \sum_{N=1}^\infty \int_{N-1}^N\int_{|x|= r} \rho^k \dd \sigma \dd r = \int_{\R^2} \rho^k \dd x
 <\infty\;,\qquad \mbox{implying that }\quad
  \lim_{N \to \infty}\int_{N-1}^N\int_{|x|=r} \rho^k \dd \sigma \dd r
  =0\, ,
\end{equation*}
for $k\in\{p,p-1\}$, and the two integration by parts for any
given $\epsilon$ are justified.

We now estimate $I_2$ and $I_3$, showing in particular that they
are finite. Starting with $I_2$, using $-\Delta c_\epsilon=
\rho_\epsilon$ where $\rho_\epsilon:=
\gamma_\epsilon*\rho*\gamma_\epsilon$, so that by H\"older's
inequality  and Young's inequality for convolutions, obtain
\begin{equation}\label{fip1}
\int_{\R^2} (\rho-K)_+^{p}(-\Delta c_\epsilon) \dd x = \int_{\R^2} (\rho-K)_+^{p}\rho_\epsilon  \dd x
\leq \|(\rho- K)_+\|_{p+1}^p\|\rho_\epsilon\|_{p+1} \leq \|(\rho- K)_+\|_{p+1}^p\|\rho \|_{p+1}
\end{equation}
Likewise for $I_3$, we use the fact that on the support of
$(\rho-K)_+$, $K \leq \rho$. Therefore
$$
K\int_{\R^2} (\rho-K)_+^{p-1}(-\Delta c_\epsilon) \dd x =
K\int_{\R^2} (\rho-K)_+^{p-1}\rho_\epsilon  \dd x \leq \int_{\R^2}
(\rho-K)_+^{p-1}\rho_\epsilon \,\rho  \dd x\ .
$$
Therefore, by H\"older's inequality  and Young's inequality
for convolutions,
\begin{equation}\label{fip2}
K\int_{\R^2} (\rho-K)_+^{p-1}(-\Delta c_\epsilon) \dd x  \leq \|(\rho- K)_+\|_{p+1}^{p-1}\|\rho \|_{p+1}^{2}\ .
\end{equation}
Applying the arithmetic-geometric mean inequality to the right side of (\ref{fip1}),  we have that for any $\nu>0$,
$$
\|(\rho- K)_+\|_{p+1}^p\|\rho \|_{p+1} \leq \frac{p}{p+1}\nu^{-(p+1)/p}\|(\rho- K)_+\|_{p+1}^{p+1} +
\frac{1}{p+1}\nu^{p+1}\|\rho \|_{p+1}^{p+1}\ .
$$
Making a similar estimate for the right hand side of (\ref{fip2}), and combining results, we have that
\begin{equation}\label{fip4}
I_2+I_3 \leq F_1(\nu) \|(\rho- K)_+\|_{p+1}^{p+1} + F_2(\nu)\|\rho \|_{p+1}^{p+1}
\end{equation}
where $F_1(\nu)$ is a positive linear combination of negative powers of $\nu$, and $F_2(\nu)$
 is a positive linear combination of positive powers of $\nu$.

By Lemma~\ref{feul}, $\rho^p$ is integrable for any $1\leq p <
\infty$, and so the right side of (\ref{fip4}) is finite. Then by
(\ref{estt2}),
$$
{\displaystyle \int_{\R^2}
\left|\nabla\left[(\rho-K)_+^{p/2}\right]\right|^2\dd x <
\infty}\, .
$$
From here we show that
\begin{equation}\label{fip3}
 \int_{\R^2} \left|\nabla\left(\rho^{p/2}\right)\right|^2\dd x \leq  \int_{\R^2} \left|\nabla\left[(\rho-K)_+^{p/2}\right]\right|^2\dd x
 +  \frac{16 K^{(2p-1)/2}}{\gamma_2}\left[\pi{\mathcal D}[\rho] + C_{{\rm CCD}}\right] \ .
 \end{equation}
 Indeed,
 $$ \int_{\R^2} \left|\nabla\left(\rho^{p/2}\right)\right|^2\dd x =  \int_{\R^2} \left|\nabla\left[(\rho-K)_+^{p/2}\right]\right|^2\dd x -
\int_{\{\rho<K\}} \left|\nabla\left(\rho^{p/2}\right)\right|^2\dd x\ ,$$
and
 \begin{equation*}
\int_{\{\rho<K\}} \left|\nabla\rho^{p/2}\right|^2\dd x \leq
16 K^{(2p-1)/2}\int_{\R^2} \left|\nabla\rho^{1/4}\right|^2\dd
x\leq
 \frac{16 K^{(2p-1)/2}}{\gamma_2}\left[\pi{\mathcal D}[\rho] + C_{{\rm CCD}}\right] \,
 ,
\end{equation*}
where we applied again the concentration controlled inequality,
Theorem~\ref{both2}, using~\eqref{boundforboth}.

Following an idea of J\"ager and Luckhaus \cite{JL}, we use the GNS inequality
\begin{equation}\label{eq:otherGNS2}
\int_{\R^2}v^{p+1}\dd x\le D_p \left(\int_{\R^2}|\nabla
v^{p/2}|^2\dd x\right) \left(\int_{\R^2}v\dd x\right)\,,
\end{equation}
which is a consequence of the Sobolev embedding inequality
$\|u\|_{L^2 (\R^2)} \leq C \|\nabla u\|_{L^1 (\R^2)}$ applied to
$u=v^{(p+1)/2}$ and Cauchy-Schwarz inequality since $\nabla u =
\frac{p+1}{p} v^{1/2} \nabla v^{p/2}$. Applying
\eqref{eq:otherGNS2} to $v=(\rho-K)_+$, we get
\begin{equation*}
\int_{\R^2}(\rho-K)_+^{p+1}\dd x\le M(K)
\int_{\R^2}\left|\nabla\left[(\rho-K)_+^{p/2}\right]\right|^2\dd x
\qquad \mbox{where } M(K):=\int_{\R^2} (\rho-K)_+ \dd x\;.
\end{equation*}
Then (\ref{fip4}) becomes
\begin{equation*}
I_2+I_3 \leq F_1(\nu) M(K)\,
\int_{\R^2}\left|\nabla\left[(\rho-K)_+^{p/2}\right]\right|^2\dd x
+ F_2(\nu)\,8\pi \int_{\R^2} \left|\nabla\rho^{p/2}\right|^2\dd x
\, .
\end{equation*}

We finally work with $I_1$ to estimate it using \eqref{fip3} as
$$
\frac{p}{2(p-1)} I_1 \leq -\int_{\R^2}
\left|\nabla\left[(\rho-K)_+^{p/2}\right]\right|^2\dd x
-\int_{\R^2} \left|\nabla \rho^{p/2}\right|^2\dd x + \frac{16
K^{(2p-1)/2}}{\gamma_2}\left[\pi{\mathcal D}[\rho] + C_{{\rm
CCD}}\right] \, .
$$

Now choose  $\nu_0>0$ small enough such that $8\pi
F_2(\nu_0)<2(p-1)/p$, and then $K_0 < \infty$  large enough such
that $M(K)F_1(\nu_0)<2(p-1)/p$. This choice of $K_0$ only depends
on $\nu_0$ and the bound on  $\int_{\R^2}\rho|\log \rho|\dd x$
since
\begin{equation*}
M(K) = \int_{\R^2}(\rho-K)_+\dd x\le\int_{\rho>K}\rho\dd x\le \frac
1{\log K}\int_{\rho>K}\rho \log \rho\dd x\le\frac 1{\log
K}\int_{\R^2} \rho \log^+ \rho \dd x\;.
\end{equation*}
We find
$$
\int_{\R^2} (\rho-K)_+^p \dd x - \int_{\R^2} (\rho_0-K)_+^p\dd x
\le \tau\frac{32 K^{(2p-1)/2}(p-1)}{p\gamma_2}\left[\pi{\mathcal
D}[\rho] + C_{{\rm CCD}}\right]\, ,
$$
for all $K\geq K_0$. The desired result follows with
$$
A_1= \frac{32 K^{(2p-1)/2}(p-1)}{p\gamma_2}C_{{\rm CCD}} \qquad
\mbox{and} \qquad A_2=\frac{32 K^{(2p-1)/2}(p-1)}{p\gamma_2}\pi\,
.
$$
\lanbox


\section{Proof of the main results}


\subsection{Approximate solutions}
We now combine the single step operations described in the
previous section to inductively define infinite sequences
$\{\rho_\tau^k\}_{k\in\NN}$ giving a discrete-time approximation
to the PKS evolution \eqref{eq:sp}. {\em For the rest of this
section}, fix any $\lambda>0$, and any density  $\rho_0$ in $\R^2$
with total mass $8\pi$ with $\rho_0\log \rho_0$ integrable, and
such that there exists $\crho$ with $\Fcfd[\rho_0] < \crho$. It
then follows from Lemma~\ref{controlcenter} that  $|x|\rho_0$ is
integrable, and from the Talagrand inequality, Theorem~\ref{tal},
that ${\rm W}_2(\rho_0,\rhohls)< \infty$. Thus, $\rho_0$ satisfies
the conditions of Theorem~\ref{lower3} on the existence of
minimizers for our single step variational problem.

Fixing an arbitrarily small parameter $\tau>0$,  we now
inductively define the sequence of densities
$\{\rho_\tau^k\}_{k\in\NN}$ with $\rho_\tau^0 := \rho_0$ by
solving the sequence of variational problems
\begin{equation}\label{vark}
\rho_\tau^k\in {\rm arg}\min_{\rho\in {\cal S}}\left\{
\frac{\W_2^2(\rho,\rho_\tau^{k-1})}{2\tau} + \Fhls^{\epsilon_k}
[\rho]\right\}
\end{equation}
for a sequence of regularization parameters
$\{\epsilon_k\}_{k\in\NN}$ to be specified now. By
Theorem~\ref{lower3}, the sequence $\{\rho_\tau^k\}_{k\in\NN}$ is
well defined no matter how we choose $0 < \tau < 1$ and
$\{\epsilon_k\}_{k\in\NN}$. We moreover define
\begin{equation}\label{poinc}
Q_k := \crho - \Fcfd[\rho_\tau^k] > 0 \qquad \mbox{for each }
k \, .
\end{equation}

\begin{lm}[Good step sizes]
Let $\{\rho_\tau^k\}_{k\in\NN}$ be the sequence of minimizers
defined inductively using \eqref{vark} starting from $\rho_\tau =
\rho_0$. With $Q_k$ defined as in \eqref{poinc},   let $\AA$ be
the constant given in {\rm Corollary~\ref{cordis}}, and let
$\Lambda$ be defined by
$$
\Lambda := \prod_{m=1}^\infty
\left(1-\frac{2^{-m}}{4} \right) \ ,
$$
and note that $1 > \Lambda>0$. Choose any $\tau>0$ satisfying
\begin{equation}\label{thresh}
\tau <  \min\left\{\frac{\Lambda
Q_0}{2\,\AA\|\gamma\|_{4/3}},1\right\} := \tau^\star\ ,
\end{equation}
and define $\epsilon_k$ by
\begin{equation}\label{taude2b}
 \tau^{1/3}\sqrt{\epsilon_k}\left[8\,\pi^{1/3}\,\AA\, \gamma_2^{-2/3} \left(\pi \,\crho+
 C_{{\rm CCD}} \right)^{2/3}\right] = \frac{Q_0}{4}\tau^2 2^{-k} \, .
\end{equation}
Then for all $k$, $Q_k > \Lambda Q_0> 0$. In particular,
\begin{equation*}
  \Fhls[\rho_\tau^k] < +\infty \quad \mbox{and}\quad \Fcfd[\rho_\tau^k] <
  \crho\,.
\end{equation*}
Note that for some constant $Z$,  $\epsilon_k :=
Z\,\tau^{10/3}\,4^{-k} $.
\end{lm}

\noindent{\bf Proof:}  We shall show by induction that for each
positive integer $j$
\begin{equation}\label{indst}
Q_j \geq \prod_{m=1}^j\left(1 - \frac{2^{-m}}{4}\right) Q_0\ ,
\end{equation}
which is somewhat more than we need since the right hand side is larger than
$\Lambda Q_0$.

We now make the inductive hypothesis that for some positive
integer $k$, (\ref{indst}) is true for all positive integers $j <
k$. Since $\Lambda < 1$, we may apply Theorem~\ref{thmstep} with
$\rho_\tau^{k-1}$ in place of $\rho_0$, and $\rho_\tau^{k}$ in
place of $\rho$ and with $\tau$ and $\epsilon_k$ specified as
above. Then the conclusion (\ref{taude3}) can be simplified and
rewritten as
\begin{equation}\label{poinc2}
\Fcfd[\rho_\tau^k]  \leq  \Fcfd[\rho_\tau^{k-1}] +\tau^2 \frac{Q_{k-1}}{4}2^{-k}\ .
\end{equation}
Since $\tau<1$, this means that
$$
Q_k \geq Q_{k-1}\left(1 - \frac{2^{-k}}{4}\right)\ .
$$
By the inductive hypothesis, we  obtain (\ref{indst}) for $j=k$.
The proof that (\ref{indst}) is valid for $j=1$ is a direct
application of  Theorem~\ref{thmstep}, in the same way,  since
$\Lambda < 1$. \lanbox

\subsection{The passage to continuous time}

{\em Throughout the rest of this section, we assume that $0 < \tau
< \tau^\star$, where $\tau^\star$ is defined in \eqref{thresh},
and that $\epsilon_k$ is defined by \eqref{taude2b}, and then that
$\{\rho_\tau^k\}_{k\in\NN}$ is a corresponding sequence of
minimizers of \eqref{vark}.}

We now interpolate between the terms of the sequence
$\{\rho_\tau^k\}_{k\in\NN}$ to produce a function from
$[0,\infty)$ to $L^1(\R^2)$ that we shall show to be, for
sufficiently small $\tau$, an approximate solution of the PKS
system.  For technical reasons, we  shall need two distinct, but
closely related, interpolations.
\begin{itemize}
\item {\it The Lipschitz interpolation:} For each positive integer
$k$, let $\nabla\varphi^k$ be the optimal transportation plan with
$\nabla\varphi^k \#\rho_\tau^k =  \rho_\tau^{k-1}$.  Then for
$(k-1)\tau \leq t \leq k\tau$ we define
$$
\rho_\tau(t)=\left( \frac{t-(k-1)\tau}{\tau} \id +\frac{k\tau
-t}{\tau}\nabla \varphi^k\right) \#\rho^{k}_\tau \ .
$$

\item {\it The piecewise constant interpolation:} For each $t$ and
each positive integer $k$ with $(k-1)\tau \leq t < k\tau$ we
define $\tilde \rho_\tau(t) = \rho_\tau^{k-1}$, with $\tilde
\rho_\tau(0) = \rho_0$.
\end{itemize}
For displacement convex functionals of $\rho$, such as
$\Fcfd[\rho]$, ${\mathcal E}[\rho]$, or the absolute first moment,
any uniform bounds on the functional along the sequence
$\{\rho_\tau^k\}_{k\in\NN}$ extend to $\rho(t)$ for all $t$, since
if ${\mathcal G}$ is such a functional, then for $(k-1)\tau < t <
k\tau$,
$$
{\mathcal G}[\rho_\tau(t)] \leq  \frac{t-(k-1)\tau}{\tau}
{\mathcal G}[\rho_\tau^{k-1}] + \frac{k\tau-t}{\tau} {\mathcal
G}[\rho_\tau^{k}]\ .
$$
Of course it is evident that for {\em any} sort of functional
${\mathcal G}[\rho]$, displacement convex or not,  a uniform bound
on ${\mathcal G}[\rho]$ along the sequence
$\{\rho_\tau^k\}_{k\in\NN}$ extends to $\tilde\rho(t)$ for all
$t$. Some of the functionals with which we work, such as
${\mathcal D}[\rho]$, are not displacement convex, and this is the
reason we need the second interpolation.

The uniform equicontinuity properties that we  prove next  explain
the utility of the first interpolation, and also why we can use
the two different interpolations at once. Since $\rho^{k}_\tau$ is
a minimizer for~\eqref{vark}, using $\rho^{k-1}_\tau$ as trial
function yields
\begin{equation*}
 \Fhls^{\epsilon_{k}}[\rho^{k}_\tau] + \frac{1}{2\,\tau} \W_2^2(\rho^k_\tau,\rho^{k-1}_\tau) \le \Fhls^{\epsilon_{k}}[\rho^{k-1}_\tau]\ ,
\end{equation*}
and hence,
\begin{equation}\label{notel}
\W_2^2(\rho^k_\tau,\rho^{k-1}_\tau) \leq 2\tau \left[ \Fhls^{\epsilon_{k}}[\rho^{k-1}_\tau] -
\Fhls^{\epsilon_{k}}[\rho^{k}_\tau]\right]\ ,
\end{equation}
and by the monotonicity of $\eps \mapsto \Fhls^\eps[\rho] $ see Lemma~\ref{monineps}
\begin{equation}\label{nophi}
  \Fhls[\rho_\tau^k] \le \Fhls^\eps[\rho_\tau^k] < +\infty\;.
\end{equation}

In standard applications of the JKO scheme, in which the
functional in the variational problem does not change from step to
step, one would sum both sides in (\ref{notel}) over a range of
values of $k$, and then the sum of the terms on the right would
telescope. This is not so in our case. However, for small
$\epsilon$, $\Fhls^\epsilon \approx \Fhls$ and we recover the
telescoping sum in a useful approximate sense. The precise version
of $\Fhls^\epsilon \approx \Fhls$ follows from (\ref{regcomp45}),
Lemma~\ref{erroreg}, and  \eqref{taude2b}, which says that
$\epsilon_{k}=Z \tau^{10/3}4^{-k}$ to get
\begin{align}
\Fhls^{\epsilon_{k}}[\rho^{k-1}_\tau]-\Fhls[\rho^{k-1}_\tau] \leq
&\,C_\gamma \, Z
\|\rho^{k-1}_\tau\|_{3/2}^{3/2}\,\tau^{10/3}4^{-k}\leq \tilde{Z}
\|\rho^{k-1}_\tau\|_{3/2}^{3/2}\,\tau^{3}2^{-k} \label{tech99}
\end{align}
for $\tau<\tilde\tau^\star:=\min(\tau^\star,(2Z\sqrt{e})^{-3/10})$
with $\tilde{Z}:=C_\gamma \, Z$ according to Lemma~\ref{erroreg}.
We thus deduce
\begin{equation*}
\W_2^2(\rho^k_\tau,\rho^{k-1}_\tau) \leq 2\tau \left(
\Fhls[\rho^{k-1}_\tau] - \Fhls[\rho^{k}_\tau]\right) +
2\tilde{Z}\|\rho^{k-1}_\tau\|_{3/2}^{3/2}\,\tau^{4} 2^{-k}\, .
\end{equation*}
Using~\eqref{nophi} and~\eqref{l32} as in the proof of Theorem
\ref{thmstep} where the concentration control inequality
\eqref{both2} is crucial, we deduce that
$\|\rho^{k-1}_\tau\|^{3/2}_ {3/2}\leq C_3/\tau$ to conclude that
\begin{equation}\label{notel2}
\W_2^2(\rho^k_\tau,\rho^{k-1}_\tau) \leq 2\tau \left(
\Fhls[\rho^{k-1}_\tau] - \Fhls[\rho^{k}_\tau]\right) +
2\tilde{Z}C_3\tau^{3} 2^{-k} \, .
\end{equation}
We are almost in a position to obtain a crucial {\it a-priori}
H\"older continuity estimate, but there is still one more
consequence of our step dependent regularization to deal with:  If
for each $k$ we had been using the functional $\Fhls$ instead of
the functional $\Fhls^{\epsilon_k}$, it would be immediate that
$k\mapsto \Fhls[\rho_\tau^k]$ would be decreasing. Since by the
Log-HLS inequality, $\Fhls$ is bounded below, this would give an
immediate upper bound on the sum of the right hand side of
(\ref{notel2}) over any range of $k$.

However, we have used our freedom to choose the sequence
$\{\epsilon_k\}_{k\in\NN}$ of regularization parameters to
converge to zero as rapidly as we may require, and hence easily
obtain:

\begin{lm}[Uniform bounds on the free energy $\Fhls$]\label{almdec}
There are positive constants $\bar{\mathcal F}_0$, $\bar{\mathcal
F}_1$ depending only on the initial data and the regularization
mollifier $\gamma$ such that for each $\tau < \tilde\tau^\star$
and each $k\in\NN$,
$$
\Fhls[\rho_\tau^k]\le \bar{\mathcal F}_0 + \bar{\mathcal F}_1
\tau^{2}.
$$
\end{lm}

\noindent{\bf Proof:}  Directly from the variational problem
(\ref{vark}) we have $\Fhls^{\epsilon_k}[\rho_\tau^k] \leq
\Fhls^{\epsilon_k}[\rho_\tau^{k-1}]$. Then, as above from
\eqref{regcomp45}, \eqref{tech99}, and \eqref{l32} we get
\begin{equation}\label{monotonefpks}
\Fhls^{\epsilon_k}[\rho_\tau^{k-1}] \leq
\Fhls^{\epsilon_{k-1}}[\rho_\tau^{k-1}] + \tilde{Z} C_3 \tau^{2}
2^{-k} \ .
\end{equation}
This means that the free energy $\Fhls$ is almost decreasing along
$\{\rho_\tau^k\}_{k\in\NN}$. A telescoping sum argument yields
$\Fhls^{\epsilon_k}[\rho_\tau^k] \leq \Fhls^{\epsilon_0}[\rho_0] +
\tilde{Z} C_3 \tau^{2}$, and then one more application of
(\ref{regcomp45}) gives
$$
\Fhls(\rho_\tau^k)\leq \Fhls^{\epsilon_0}[\rho_0] + \tilde{Z} C_3
\tau^{2} \leq {\mathcal E}[\rho_0] + 32\,\pi +
2\|\rho_0\|_{L^1(\R^2,\log(e+|x|^2)\dd x)} + \tilde{Z} C_3
\tau^{2}:= \bar{\mathcal F}_0 + \bar{\mathcal F}_1 \tau^{2}\, ,
$$
where \eqref{regcomp477} was used. \lanbox

\

We are now ready to prove the H\"older continuity estimate.

\begin{lm}[H\"older continuity]\label{hocon}
There is a positive constant $\bar{\mathcal F}_2$ depending only
on the initial data and the regularization mollifier $\gamma$ such
that for each $\tau < \tilde\tau^\star$ and each $k\in\NN$, such
that for all $t> s \geq 0$,
$$
 \W_2(\rho_\tau(t),\rho_\tau(s) ) \leq \bar{\mathcal F}_2\left(t-s\right)^{1/2}\ .
$$
\end{lm}

\noindent{\bf Proof:} Let $j$ be such that $(j-1)\tau \leq s \leq
j\tau$ and let $\ell$ be such that $\ell\tau \leq t \leq
(\ell+1)\tau$. By the geodesic property of McCann's displacement
interpolation,
\begin{equation*}
\W_2(\rho_\tau(s),\rho_\tau^j) = \frac{j\tau -s}{\tau}
\W_2(\rho_\tau^{j-1},\rho_\tau^j) \qquad{\rm and}\qquad
 \W_2(\rho_\tau(t),\rho_\tau^{\ell} ) = \frac{t-\ell\tau}{\tau}\W_2(\rho_\tau^{\ell},\rho_\tau^{\ell+1})
\ .
\end{equation*}
By Lemma~\ref{almdec} and the Log-HLS inequality
\eqref{regcomp46}, $\Fhls[\rho^{k}_\tau] -  \Fhls[\rho^{k-1}_\tau]
\leq \bar{\mathcal F}_0 + \bar{\mathcal F}_1 \tau^{2}  -8\pi(\log
8 -1)$, and thus, by plugging into \eqref{notel2}, we get
\begin{align}
\W_2^2(\rho^k_\tau,\rho^{k-1}_\tau) &\leq 2\tau \left[
\bar{\mathcal F}_0 + \bar{\mathcal F}_1 \tau^{2}  -8\pi(\log 8
-1)\right] + 2\tilde{Z}C_3\tau^{3} 2^{-k} \nonumber \\
&\leq \tau \left[ 2\bar{\mathcal F}_0 + 2\bar{\mathcal F}_1  -
16\pi(\log 8 -1) + 2\tilde{Z}C_3\right]:=\bar{\mathcal F}_2^2 \tau
\, ,\label{tech98}
\end{align}
since $\tau < 1$ and $k\in\NN$. Therefore, we deduce
$$
\W_2(\rho_\tau(s),\rho_\tau^j) \leq \frac{j\tau
-s}{\tau}\bar{\mathcal F}_2\sqrt{\tau}\qquad \mbox{and} \qquad
\W_2(\rho_\tau(t),\rho_\tau^{\ell} ) \leq
\frac{t-\ell\tau}{\tau}\bar{\mathcal F}_2\sqrt{\tau} \,.
$$
Adding these two estimates and using the concavity of square root,
\begin{equation}\label{endpei2}
\W_2(\rho_\tau(s),\rho_\tau^j) + \W_2(\rho_\tau(t),\rho_\tau^\ell)
\leq \bar{\mathcal F}_2\sqrt{ (t-s) - (\ell- j)\tau }\ .
\end{equation}
Next, by the triangle inequality, the Cauchy-Schwartz inequality,
and \eqref{notel2} and proceeding as in \eqref{tech98}, we finally
conclude
\begin{align}
\W_2(\rho_\tau^\ell,\rho_\tau^{j-1} ) \leq& \, \sum_{k= j}^\ell
\W_2(\rho^k_\tau,\rho^{k-1}_\tau) \leq
  (\ell - j)^{1/2}\left(\sum_{k= j}^\ell \left(2\tau \left[ \Fhls[\rho^{k-1}_\tau] -
\Fhls[\rho^{k}_\tau]\right] + 2\tilde{Z}C_3\tau^{3} 2^{-k}\right) \right)^{1/2}\nonumber\\
\leq &[(\ell - j)\tau]^{1/2}\left(2\left[ \Fhls[\rho^{j-1}_\tau] -
\Fhls[\rho^{k}_\tau]\right] +
2\tilde{Z}C_3\tau^{2}\right)^{1/2}\leq \bar{\mathcal F}_2\,[(\ell
- j)\tau]^{1/2} .\nonumber
\end{align}
Adding this to the estimate in (\ref{endpei2}), and using the
subadditivity of the square root concludes the proof. \lanbox

\subsection{Weak compactness}

In this subsection, we will show the compactness of the sequence
of interpolating curves. We cannot proceed as usually done, for
instance in \cite{AGS}, since we want to show that the limiting
curves are not only measures but rather densities for each time
and also since our densities have infinite second moment. The idea
is to show a compactness in metrics $W_p$, $1\leq p <2$, less
strong than $W_2$ and then, pass to the limit $p\nearrow 2$ the
H\"older continuity property.

The next lemma shows a compactness property of the sets
$\{\rho_\tau(t)\ :\ 0 < \tau < \tilde\tau^\star\}$ for each fixed
$t$.

\begin{lm}[Uniform integrability at fixed $t$]\label{ufix}
There is a finite and computable constant $\bar{\mathcal F}_3$
depending only on $\rho_0$ and for any fixed $1\leq p <2$ so that
for all $\tau<\tilde\tau^\star$
$$
{\mathcal E}[\rho_\tau(t)] = \int_{\R^2}\rho_\tau(t,x)\log
\rho_\tau(t,x) \dd x \leq \bar{\mathcal F}_3 \qquad \mbox{and}
\qquad \int_{\R^2}|x|^p \rho_\tau(t,x) \dd x \leq \bar{\mathcal
F}_3\, .
$$
\end{lm}

\noindent{\bf Proof:} By  the uniform control that we have on
$\Fcfd$, the moment bound is immediate from
Lemma~\ref{controlcenter}. By the uniform control of $\Fhls$ and
$\Fcfd$ in~\eqref{boundforboth}, Lemma~\ref{almdec}, and by the
first concentration control Theorem \ref{both} we conclude that
\begin{equation}
  \label{uniflog}
  \gamma_1\int_{\R^2}\rho_\tau^k \log_+ \rho_\tau^k \dd x  \le
 \bar{\mathcal F}_0 + \bar{\mathcal F}_1
\tau^{2}+C_{{\rm CCF}}\,
\end{equation}
where $0<\gamma_1\le 1$, uniformly in $k$. Finally, using the
bound of the absolute first moment together with \eqref{uniflog},
we conclude that ${\mathcal E}[\rho_\tau^k]$ is bounded uniformly
in $k$. Then, by the displacement convexity of  ${\mathcal E}$,
this bound extends to $\rho_\tau(t)$ for all $t>0$, as explained
at the beginning of this subsection. \lanbox

\

It follows immediately from \eqref{wascom} that Lemma~\ref{hocon}
remains true if $\W_2$ there is replaced by any of the weaker
metrics $\W_p$, $1\leq p < 2$. The following characterization of
the convergence in $\W_p$ metrics in \cite[Chapter 9]{Villani}:
convergence of the absolute $p$-moment plus the weak-* convergence
as measures of a sequence of densities $\{\rho_n\}_{n\in\NN}$
towards $\rho$ is equivalent to $W_p$ convergence; implies the
following compactness result.

\begin{lm}[Compactness for the $\W_p$ metric]\label{w1com}
For any $M>0$, let $K$ be a subset of the set of densities
$\rho$ of mass $M$ that is uniformly integrable, and such that
$\{|x|^p\rho(x)\ :\ \rho\in K\}$ is also uniformly integrable.
Suppose also that $K$ is closed in the $L^1$-weak topology. Then
$K$ is compact in the $W_p$ metric.
\end{lm}

\noindent{\bf Proof:} Let $\{\rho_n\}_{n\in\NN}$ be any sequence
in $K$. Since $K$ is uniformly integrable and weakly closed in
$L^1$, the Dunford-Pettis Theorem provides us with a $\rho\in K$
and a subsequence $\{\rho_{n_k}\}_{k\in\NN}$ such that
$\lim_{k\to\infty}\rho_{n_k} = \rho$ weakly in $L^1$ and thus,
weakly-* as measures. It is trivial to check that weak-$L^1$
convergence plus the uniform integrability of $\{|x|^p\rho(x)\ :\
\rho\in K\}$ implies that
$$
\lim_{k\to\infty} \int_{\R^2}|x|^p \rho_{n_k}(x) \dd x =
\int_{\R^2}|x|^p \rho(x) \dd x \, .
$$
The  characterization of $\W_p$ convergence mentioned above then  implies
$\lim_{k\to\infty}\W_p(\rho_{n_k},\rho) = 0$. \lanbox

\begin{thm}[Convergence as $\tau\to 0$]\label{contauthm}
Given $T>0$ and any $1\leq p<2$, define $({\mathcal M}_T^\delta,
\W_p)$ to be the metric space in which ${\mathcal M}_T^\delta$ is
the set of densities on $\R^2$ satisfying ${\mathcal E}[\rho] \leq
\bar{\mathcal F}_3$ and
$$
\int_{\R^2}|x|^{p+\delta}\rho(x) \dd x \leq \bar{\mathcal F}_3(1+
T^{(p+\delta)/2})\, ,
$$
with $p<p+\delta<2$. Then there is a function $\rho$ on
$[0,\infty)$ with values in the set of densities of mass $8\pi$
such that for all $T>0$, the restriction of $\rho$ to $[0,T]$ is
continuous in $({\mathcal M}_T^\delta, \W_p)$, and there is a
sequence $\{\tau^{(n)}\}_{n\in\NN}$ such that for all $T>0$,
\begin{equation}\label{araz}
\lim_{n\to\infty} \left[ \max_{0 \leq t \leq
T}\W_p(\rho_{\tau^{(n)}}(t), \rho(t))\right] = \lim_{n\to\infty}
\left[ \max_{0 \leq t \leq T}\W_p(\tilde\rho_{\tau^{(n)}}(t),
\rho(t))\right] = 0\ .
\end{equation}
Moreover the sequence $\{\tau^{(n)}\}_{n\in\NN}$ can be chosen
independently of $p$, i.e., such that the convergence property
\eqref{araz} holds for all $1\leq p < 2$. Furthermore, as a
consequence for all $t> s \geq 0$ and all $1\leq p \leq 2$:
\begin{equation}\label{cont}
 \W_p(\rho(t),\rho(s) ) \leq \bar{\mathcal F}_2\left(t-s\right)^{1/2}\ .
\end{equation}
\end{thm}

\noindent{\bf Proof:}  For each $T>0$, $({\mathcal M}_T^\delta,
\W_p)$ is a compact metric space as a consequence of
Lemma~\ref{w1com}. By Lemma~\ref{ufix}, for each $t \leq T$, and
each $0< \tau<\tilde\tau_\star$, the restriction of $\rho_\tau$ to
$[0,T]$ takes values in $({\mathcal M}_T^\delta, \W_p)$. Next, by
(\ref{wascom}) and Lemma~\ref{hocon}, the set of these functions
for $0 < \tau < \tilde\tau_\star$ is uniformly equicontinuous into
$({\mathcal M}_T^\delta, \W_p)$.  Thus by the the Arzela-Ascoli
Theorem, we can select a uniformly convergent sequence.

Now a simple diagonal sequence argument concludes the construction
of $\rho$ and proof of (\ref{araz}). Apply the above for $T=1$ to
get the initial sequence. Now take $T=2$ and choose a subsequence
of the first sequence, and so forth. For the
piecewise-interpolation sequence, note that $\W_2(\tilde
\rho_\tau(t),\rho_\tau(t))\! \leq \! C \sqrt\tau$ by standard
properties of displacement interpolation and \eqref{tech98} in
Lemma~\ref{hocon}. Therefore, the limits $\tau^{(n)}\to 0$ of both
time interpolations are the same. Note that this argument does not
invoke any uniqueness of the limit.

Another simple diagonal sequence argument shows that the sequence
can be made independent of $p$. Take a sequence of increasing
exponents $\{p_n\}_{n\in\NN}\nearrow 2$ and
$\{\delta_n\}_{n\in\NN}\searrow 0$ with $1\leq p_n<p_n+\delta_n<2$
and apply the diagonal sequence argument to the constructed
sequences for each $p_n$. Also, take into account that the
sequence of distances $\W_p(\rho_{\tau^{(n)}}(t), \rho(t))$ is
increasing in $p$.

The last part of the claim \eqref{cont} follows directly from
\eqref{araz} and Lemma~\ref{hocon} for all $1\leq p<2$. Since the
constant $\bar{\mathcal F}_2$ obtained in Lemma~\ref{hocon} does
not depend on $1\leq p<2$, then we conclude \eqref{cont} for
$p=2$.\lanbox

\subsection{$L^p$ regularity}

Our goal in this section is to prove:

\begin{thm}[$L^p$-regularity]\label{lpreg}
For each finite $a>0$ and $p>1$, there is a finite and computable
constant $C_p$ depending only on  $a$, $p$ and $\rho_0$ such that
whenever $\tau < a$,
$$
\int_{\R^2} \rho_\tau^p(t,x)\dd x \leq \tilde C_p \qquad \mbox{for
all } t\geq a \, .
$$
\end{thm}

\noindent{\bf Proof:} From (\ref{taude3}), we deduce that for each
$m$,  ${\displaystyle \Fcfd[\rho_\tau^m] + {\mathcal
D}[\rho_\tau^m]  \leq \Fcfd[\rho_\tau^{m-1}] +
\frac{Q_m}{4}2^{-m}\tau^2}$, proceeding in the same way that we
did in deducing (\ref{poinc2}), except this time we do not discard
the dissipation term. Let  $n\geq k$ be positive integers. Since
$Q_m \le \crho$ for all $m$, summing from $m=n-k$ to $n$ yields
\begin{equation}\label{eq:discretedisener3}
\Fcfd[\rho_\tau^{(n)}] + \sum_{m=n-k}^n \tau {\mathcal
D}[\rho_\tau^m] \leq \Fcfd[\rho_\tau^{n-k-1}] +
\frac{\crho}{4}\tau^2 \, .
\end{equation}
Then since $0 \leq \Fcfd[\rho_\tau^k] \le \crho$ for all $k$,
using $\tau < 1$ and dividing by $k\tau$, we may simplify this to
\begin{equation}\label{sumbnd}
 \frac1k \sum_{m=n-k}^n {\mathcal D}[\rho_\tau^m] \leq \frac{2\crho}{k\tau}\ .
\end{equation}
We now choose $k$ to be the greatest integer less than or equal to
$a/\tau$, and of course suppose that $n>k$. Since $k\tau \leq a <
(k+1)\tau $, and $k\geq 1$, $a/2 \leq k\tau$, and then the fact
that averages dominate minima yields the conclusion that for some
positive integer $m$ with $\tau\leq a$,
$$
{\mathcal D}[\rho_\tau^m] \leq \frac{4\crho}{a}\, .
$$
Then since $\Fcfd[\rho_\tau^m] \le \crho$ and $\Fhls[\rho_\tau^m]<+\infty $ we have from
Theorem~\ref{both2} that
$$
\int_{\R^2}|\nabla (\rho_\tau^m)^{1/4}|^2\dd x \leq \frac{4\pi\crho}{a \gamma _2} + \frac{C_{{\rm CCD}}}{\gamma_2} \ .
$$
Recall the GNS inequality
$$
\|f\|_q^q \le B_q\|\nabla f\|_2^{q-4}\| f\|_4^4\ ,
$$
valid for locally integrable functions $f$ in $\R^2$ and $q> 4$. Applying this with
$q=4p$ and $f=  (\rho_\tau^m)^{1/4}$, we obtain
\begin{equation}\label{meev}
 \|\rho_\tau^m\|_p^{p} \leq \left[B_{4p}\left( \frac{4\pi \crho}{a \gamma _2} +
\frac{C_{{\rm CCD}}}{\gamma_2}\right)^{2(p-1)}8\pi\right]^{1/4}\ .
\end{equation}
Thus we have an {\em a-priori} bound on  $\|\rho_\tau^m\|_p^p$ for
{\em some} $m$ with $n-k \leq m \leq n$. We now apply
Lemma~\ref{propagate}, in which the constant $A_1$ and $A_2$ in
(\ref{bell}) are defined, to conclude that
\begin{equation}\label{bell}
\int_{\R^2}(\rho_\tau^{(n)} - K)_+^p\dd x  \leq
\int_{\R^2}(\rho_\tau^m - K)_+^p\dd x + A_1k\tau +A_2\sum_{j=m}^n
\tau{\mathcal D}[\rho_\tau^j]\ .
\end{equation}

We bound the first term on the right in~\eqref{bell} using
(\ref{meev}), the second using the fact that $k\tau \leq a$, and
the third using (\ref{sumbnd}).  The result,
$$
\int_{\R^2}(\rho_\tau^{(n)} - K)_+^p\dd x \leq \left[B_{4p}\left(
\frac{4\pi \crho}{a \gamma _2} + \frac{C_{{\rm
CCD}}}{\gamma_2}\right)^{2(p-1)}8\pi\right]^{1/4} + A_1a + 2A_2
\crho
$$
uniformly for all $n$ such that $n\tau\geq a$. Note that the bound
depends only on $a$, $p$. Since $\|\rho_\tau^{(n)}\|_p \leq
\|(\rho_\tau^{(n)} - K)_+\|_p + K^{(p-1)/p}(8\pi)^{1/p}$, we have
the same type of bound on $\|\rho_\tau^{(n)}\|_p$, uniformly  for
all $n$ such that $n\tau\geq a$.  By the displacement convexity of
$$
{\displaystyle \rho\mapsto \int_{\R^2}\rho^p(x)\dd x}
$$
for $p>1$, this bound immediately extends to $\rho_\tau(t)$
for all $t\geq a$. \lanbox

\subsection{Verification that $\rho = \lim_{\tau\to 0}\rho_\tau$ is a solution of the PKS system }

Let ${\tau^{(n)}}$, $\rho_{\tau^{(n)}}$ and $\rho$ be given as in
Theorem~\ref{contauthm}. Our goal in this subsection is to prove
that $\rho$ is a weak solution of the PKS system as given in item
(1.5.2) of the Definition~\ref{prdisso}.

\begin{lm}[$\rho$ is a weak solution of the PKS system]\label{iswkso}
Let ${\tau^{(n)}}$, $\rho_{\tau^{(n)}}$ and $\rho$ be given as in
Theorem~{\rm\ref{contauthm}}. Then for all smooth and compactly
supported test functions $\zeta$ and all $t_2 > t_1 \geq 0$,
\begin{align*}
\int_{\RR^2} \!\zeta (x) \left[ \rho(t_2,x) -
\rho(t_1,x)\right]\dx\! = &-\frac{1}{4\pi}
\int_{t_1}^{t_2}\!\!\iint_{\R^2 \times \RR^2}
\!\!\!\!\!\!\rho(s,x)\,\rho(s,y)\,\frac{(x-y)\cdot(\nabla
\zeta(x)- \nabla \zeta(y))}{|x-y|^2}\dd y\dd
x\nonumber\\
&+\int_{t_1}^{t_2}\!\!\int_{\RR^2}\!
\Delta\zeta(x)\,\rho(s,x)\dx\ds\, ,
\end{align*}
\end{lm}

\noindent{\bf Proof:} In order to prove Lemma~\ref{iswkso}, we
first remind the reader an analog for the functions $\rho_\tau$:

\begin{lm}[Approximate weak solutions of the PKS system]
For $0 < \tau < 1$, define $\epsilon(t,\tau)=\epsilon_k$ for
$t\in((k-1)\tau,k\tau]$, and all integers $k\geq1$. Then for all
smooth and compactly supported test functions $\zeta$ and all $t_2
> t_1 \geq 0$,
\begin{align*}
\int_{\RR^2} \!\!\!\zeta (x) \left[ \rho_\tau(t_2,x) -
\rho_\tau(t_1,x)\right]\!\!\dx =&\, \frac12
\int_{t_1}^{t_2}\!\!\iint_{\R^2 \times \RR^2}
\!\!\!\!\rho_\tau(s,x)\,\rho_\tau(s,y)\,\nabla
G_{\epsilon(t,\tau)}(x-y) \cdot
\left(\nabla\zeta(x)-\nabla\zeta(y)\right)\!\dd y\dd
x\nonumber\\
&+\int_{t_1}^{t_2}\!\!\int_{\RR^2}\!
\Delta\zeta(x)\,\rho_\tau(s,x)\dx\ds + O(\tau^{1/2})\, .
\end{align*}
\end{lm}

The proof of the previous Lemma follows the same lines as in
\cite[Theorem 3.4]{BCC} and we skip it here for the sake of
conciseness. The interested reader can see its proof in the
preprint version of this paper.

As $\tau\to 0$ along the sequence $\{\tau^{(n)}\}_{n\in\NN}$,
$\lim_{n\to\infty}\W_1[\rho_{\tau^{(n)}}(t),\rho(t)] = 0$
uniformly on $[0,T]$ for any finite $T$. Hence by the
Kantorovich-Rubinstein Theorem,
\begin{equation*}
\lim_{n\to\infty}
\int_{\RR^2}\Delta\zeta(x)\,\rho_{\tau^{(n)}}(t,x)\dd x =
\int_{\R^2}\Delta\zeta(x)\,\rho(t,x)\dd x \ ,
\end{equation*}
uniformly on $[0,T]$. The interaction term can be rewritten as
\begin{align}\label{interbeforelim}
\int_{\R^2}\,\left[\nabla G_{\epsilon(t,{\tau^{(n)}})}\right. &\ast \left.\rho_{\tau^{(n)}}(t)\right](x)  \cdot \nabla\zeta(x)\,\rho_{\tau^{(n)}}(t,x)\,\dd x \notag\\
= &\, \int_{\R^2}\,\left[\nabla G \ast
\gamma_{\epsilon(t,{\tau^{(n)}})} \ast
\rho_{\tau^{(n)}}(t)\right](x) \cdot
\left[\gamma_{\epsilon(t,{\tau^{(n)}})}\ast
(\rho_{\tau^{(n)}}(t)\,\nabla\zeta)\right](x)\,\dd x\notag\\
= &\, -\frac{1}{4\pi}\iint_{\R^2 \times
\RR^2}(\gamma_\epsilon*\rho_{\tau^{(n)}})(t,x)\,(\gamma_\epsilon*\rho_{\tau^{(n)}})(t,y)\,\frac{(x-y)\cdot(\nabla
\zeta(x)- \nabla \zeta(y))}{|x-y|^2} \dd x\dd y \notag\\
&+ \int_{\R^2}\,\left[\nabla G \ast
\gamma_{\epsilon(t,{\tau^{(n)}})} \ast
\rho_{\tau^{(n)}}(t)\right](x) \cdot
\left[\gamma_{\epsilon(t,{\tau^{(n)}})}\ast
(\rho_{\tau^{(n)}}(t)\,\nabla\zeta)-(\gamma_{\epsilon(t,{\tau^{(n)}})}
\ast
\rho_{\tau^{(n)}}(t))\nabla\zeta \right](x)\,\dd x \notag\\
:= &\, I_1 + I_2.
\end{align}
As $\{\rho_{\tau^{(n)}}(t)\}_{n\in\NN}$ converges weakly in
$L^1(\R^2)$ towards $\rho(t)$ as $n\to\infty$, so does
$\{\gamma_{\epsilon(t,{\tau^{(n)}})}*\rho_{\tau^{(n)}}(t)\}_{n\in\NN}$.
We then deduce that
$\{(\gamma_{\epsilon(t,{\tau^{(n)}})}*\rho_{\tau^{(n)}}(t))
\otimes
(\gamma_{\epsilon(t,{\tau^{(n)}})}*\rho_{\tau^{(n)}}(t))\}_{n\in\NN}$
converges weakly in $L^1(\R^2\times\R^2)$ towards $\rho(t) \otimes
\rho(t)$ when $n\to\infty$, see \cite[Lemma~2.3]{BCC}. As a
consequence we can pass to the limit in the first term in the
right-hand-side of~\eqref{interbeforelim} to obtain
$$
 \lim_{n\to\infty} I_1 =-\frac{1}{4\pi}\iint_{\R^2 \times \RR^2}\rho(t,x)\,\rho(t,y)\,\frac{(x-y)\cdot(\nabla \zeta(x)- \nabla \zeta(y))}{|x-y|^2} \dd x\dd y\ .
$$
We must now show that  $I_2$ disappears in
the limit.  We can estimate $I_2$ using
\begin{align*}
|\gamma_{\epsilon(t,{\tau^{(n)}})}\ast
(\rho_{\tau^{(n)}}(t)\,\nabla\zeta)-&(\gamma_{\epsilon(t,{\tau^{(n)}})}
\ast \rho_{\tau^{(n)}}(t))\nabla\zeta| (x) \le \int_{\R^2}
\gamma_{\epsilon(t,{\tau^{(n)}})}(x-y)\,
|\nabla\zeta(y)-\nabla\zeta(x)|
\,\rho_{\tau^{(n)}}(t,y)\,\dd y \\
&\le C_\zeta \int_{\R^2} \gamma_{\epsilon(t,{\tau^{(n)}})}(x-y)
\,|x-y|\, \rho_{\tau^{(n)}}(t,y)\,\dd y = C_\zeta\,
((\gamma_{\epsilon(t,{\tau^{(n)}})}|x|)\ast
\rho_{\tau^{(n)}}(t))(x) \, .
\end{align*}
By the HLS  inequality, \eqref{spclhls},
$$
|I_2| \leq \frac{C_{\rm HLS}}{2\pi}
\|\gamma_{\epsilon(t,{\tau^{(n)}})} \ast
\rho_{\tau^{(n)}}(t)\|_{4/3} \,
\|\gamma_{\epsilon(t,{\tau^{(n)}})}\ast
(\rho_{\tau^{(n)}}(t)\,\nabla\zeta)-(\gamma_{\epsilon(t,{\tau^{(n)}})}
\ast \rho_{\tau^{(n)}}(t))\nabla\zeta\|_{4/3} .
$$
Then by similar arguments similar to those used to prove
Lemma~\ref{displconvexitydelta}, we get
$$
|I_2|  \le 4\,C_{\rm HLS}\,C_\zeta \|\rho_{\tau^{(n)}}(t)\|_{4/3}
\, \|\gamma |x|\|_{4/3} \sqrt{\epsilon(t,{\tau^{(n)}})} .
$$
In case $t_1>0$, estimating $\|\rho_{\tau^{(n)}}(t)\|_{4/3}$ using
Theorem~\ref{lpreg}, we obtain the result. If $t_1=0$ we can use
instead~\eqref{l32} and H\"older's inequality to obtain
$$
|I_2| \leq 4\,C_{\rm HLS}\,C_\zeta\,(8\pi)^{1/4}\,\sqrt{2\,C_3}\,
\|\gamma |x|\|_{4/3}
\sqrt{\frac{\epsilon(t,\tau^{(n)})}{\tau}}=O((\tau^{(n)})^{7/6})\,
,
$$
where \eqref{taude2b} was used. \lanbox

\subsection{Strong Compactness}

At this point we have shown that the limit $\rho =
\lim_{n\to\infty}\rho_{\tau^{(n)}}$ posesses the properties
(1.5.1) and (1.5.2) in Definition~\ref{prdisso} of properly
dissipative weak solutions. In this subsection, we show that
(1.5.3) is also satisfied. This will complete the proof of the
existence of properly dissipative solutions in Theorem~\ref{main}.
Choosing $n=k$ in (\ref{eq:discretedisener3}) we obtain, for all
$k\in\N$ and all $\tau=\tau^{(n)}$ that
\begin{equation*}
\Fcfd[\rho_\tau^k] + \sum_{m=1}^k \tau{\mathcal D}[\rho_\tau^m]
\leq \Fcfd[\rho_0] + \frac{\crho}{4}\tau^2 \, .
\end{equation*}
Thus, using the discrete time interpolation $\tilde \rho_\tau$, we
have that for any  $T>0$ and the positive integer $N$ such that
$N\tau \leq T \leq (N+1)T$,
\begin{equation}\label{stpt8}
\Fcfd[\tilde \rho_\tau(T)] + \int_0^{N\tau} {\mathcal D}[\tilde
\rho_\tau(t)] \dd t \leq \Fcfd[\rho_0]+ \frac{\crho}{4}\tau^2\ .
\end{equation}
We emphasize that the use of the piecewise constant interpolation
is essential at this point since the functional ${\mathcal
D}[\rho]$ is not displacement convex.

Note that the $L^p$ bounds deduced in Theorem~\ref{lpreg} apply to
$\tilde \rho_\tau$ as well as to $\rho_\tau$. To make full use of
these bounds, we choose any fixed $a>0$, and then for all $\tau <
a$, we weaken the bound in (\ref{stpt8}) by increasing the lower
limit of integration in $t$ to $a$.  Also writing $b:= N\tau$,
this yields
\begin{equation*}
\Fcfd[\tilde \rho_\tau(T)] + 8\int_a^b \int_{\R^2}|\nabla \tilde
\rho_\tau^{1/4}(t,x)|^2 \dd t -   \int_a^b \int_{\R^2} \tilde
\rho_\tau^{3/2}(t,x) \dd x \dd t \leq \Fcfd[\rho_0]+
\frac{\crho}{4}\tau^2\ .
\end{equation*}
It is legitimate to express ${\mathcal D}[\tilde \rho_\tau]$ as
the difference of two integrals since Theorem~\ref{lpreg} tells us
the $\tilde \rho_\tau^{3/2}$ is integrable over $[a,T]\times
\R^2$. We now show that passing to a further subsequence of
$\{\tau^{(n)}\}_{n\in\NN}$, we may arrange that  for all $0 < a <
b < \infty$, along this subsequence,
\begin{equation}\label{stpt10b}
\lim_{n\to\infty}\int_a^b \int_{\R^2} |\tilde
\rho_{\tau^{(n)}}(t,x) - \rho(t,x)|^{3/2} \dd x \dd t = 0
\end{equation}
and
\begin{equation*}
\lim_{n\to\infty}\tilde \rho_{\tau^{(n)}}(t,x)= \rho(t,x)
\end{equation*}
for almost every $(t,x) \in [a,b]\times\R^2$.  The following
strong compactness theorem, leading to the existence of almost
everywhere convergent subsequences, is the key:

\begin{thm}[Strong compactness for $\tilde\rho_\tau^{2}(t,x)$]\label{strocom}
Let $0 < a < b < \infty $ be given. The family of functions
$$
\{\ 1_{[a,b]}(t)\tilde\rho_\tau(t,x)\ :\ \tau < \tilde \tau^\star\
\}
$$
is precompact in the strong topology in $L^2(\R^2\times (a,b))$.
\end{thm}

\noindent{\bf Proof:}  By the Kolmogorov Compactness  Theorem
\cite{HH,KOL}, also known as the Frechet-Kolmogorov theorem
\cite[Corollary IV.26]{Brezis}, a set ${\mathcal G}$ of functions
$g\in L^2((a,b)\times\R^2)$ is strongly precompact if and only if:
\begin{itemize}
\item[(K1)] ${\mathcal G}$ is uniformly bounded in
$L^2((a,b)\times\R^2)$.

\item[(K2)] For each $\epsilon>0$, there is an $R>0$ so that for
all $g\in {\mathcal G}$,
\begin{equation*}
\int_a^b \int_{|x|>R} |g(t,x)|^2 \dd x \dd t < \epsilon \,.
\end{equation*}

\item[(K3)] For each $\epsilon>0$, there is a $\delta > 0$ so that
for all  $g\in {\mathcal G}$, whenever $y\in \R^2$ satisfy $|y|
\leq \delta$,
\begin{equation}\label{kol2}
\int_{\R^2\times (a,b)} |g(t,x-y) - g(t,x)|^2\dd x\dd t   \leq
\epsilon\ ,
\end{equation}
and whenever $c\in \R$, and $[\tilde a,\tilde b]\subset (a,b)$
with $(\tilde a -c,\tilde b -c)\subset (a,b)$
\begin{equation}\label{kol3}
\int_{\R^2\times (\tilde a,\tilde b)} |g(t-c,x) - g(t,x)|^2\dd
x\dd t \leq \epsilon\ .
\end{equation}
\end{itemize}

{\it Step 1.- Uniform integrability for
$\tilde\rho_\tau^{2}(t,x)$:} First, note that for each $\tau$,
$\tilde\rho_\tau^{2} =
\tilde\rho_\tau^{1/2}\tilde\rho_\tau^{3/2}$. Therefore, by the
Cauchy-Schwartz inequality
$$
\int_{\R^2}|x|^{1/2}\tilde\rho_\tau^{2}(t,x)\dd x \leq
\left(\int_{\R^2}|x|\tilde\rho_\tau(t,x)\dd x\right)^{1/2}
\left(\int_{\R^2}\tilde\rho_\tau^{3}(t,x)\dd x\right)^{1/2}\  .
$$
Recall the first moment is controlled by $\Fcfd$ in
Lemma~\ref{controlcenter}, and that $\Fcfd[\tilde\rho_\tau(t)]
\le \crho$. Also, Theorem~\ref{lpreg} give us a bound on
$\|\tilde\rho_\tau(t)\|_{3}$ uniformly in $t\in [a,b]$ for all
sufficiently small $\tau$. Thus, there is a constant $C$ depending
only on $a$ and $b$ so that for all sufficiently small $\tau$,
$$
\int_{a}^b\int_{\R^2}|x|^{1/2}\tilde\rho_\tau^{2}(t,x)\dd x \dd
t\leq C\ .
$$
Even more simply, by Theorem~\ref{lpreg} we have a constant $C$
depending only on $a$ and $b$  so that for all sufficiently small
$\tau$, such that
$$
\int_{a}^b\int_{\R^2}\tilde\rho_\tau^{2}(t,x)\dd x \dd t\leq
C\qquad{\rm and}\qquad
\int_{a}^b\int_{\R^2}\tilde\rho_\tau^{3}(t,x)\dd x \dd t\leq C\
.
$$
The uniform integrability is an immediate consequence of these
estimates, giving the first two conditions (K1) and (K2) of the
Kolmogorov compactness theorem.

\

{\it Step 2.- Spatial translations \eqref{kol2}:} Writing
$\tilde\rho_\tau(t,x) = f^4(t,x)$, we have
$$
\tilde\rho_\tau (t,x-y) - \tilde\rho_\tau (t,x)  = -4|y|\cdot
\int_0^1 f^3(t,x-sy)\nabla f(t,x-sb )\dd s\ ,
$$
and then by Minkowski's inequality, and then H\"older's inequality
\begin{eqnarray}
\int_{\R^2} |\tilde\rho_\tau (t,x-y) - \tilde\rho_\tau (t,x)|^2\dd
x
&\leq& 4|y| \left(\int_{\R^2}f^6(t,x)\dd x\right)^{1/2} \left(\int_{\R^2}|\nabla f(t,x)|^2\dd x\right)^{1/2}\nonumber\\
&=& 4|y| \left(\int_{\R^2}\tilde\rho_\tau ^{3/2}(t,x)\dd
x\right)^{1/2} \left(\int_{\R^2}|\nabla\tilde\rho_\tau
^{1/4}|^2\dd x\right)^{1/2}\, .\nonumber
\end{eqnarray}
By our uniform $L^p$ bounds,  there is a constant $C$ independent
of $\tau$ so that
$$
\int_{\R^2} |\tilde\rho_\tau (t,x-y) - \tilde\rho_\tau (t,x)|^2\dd
x  \leq |y|C \left(\int_{\R^2}|\nabla\tilde\rho_\tau ^{1/4}|^2\dd
x\right)^{1/2}\ ,
$$
and hence
$$
\int_{a}^b\int_{\R^2} |\tilde\rho_\tau (t,x-y) - \tilde\rho_\tau
(t,x)|^2\dd x \dd t \leq C|y|b^{1/2}
\left(\int_{a}^b\int_{\R^2}|\nabla\tilde\rho_\tau ^{1/4}|^2\dd
x\dd t \right)^{1/2}\ .
$$
Since the integral on the right hand side is bounded uniformly in
$\tau$, this gives us \eqref{kol2}.

\

{\it Step 3.- Temporal translations \eqref{kol3}:} We do not have
any estimates on time derivatives, so we cannot obtain the bound
on temporal translations in such a simple manner as we have for
the spatial translations. What we do have from Lemma~\ref{hocon}
is a finite constant $C$ so that
\begin{equation}\label{inter55}
\W_2(\rho(t-c,\cdot),  \rho(t,\cdot)) \le C \sqrt{c}
\end{equation}
holds uniformly in $\tau \leq \tilde \tau^\star$ and in $t\geq c$.

We now use an interpolation argument based on an idea of Otto, see
\cite[Subsection 3.5]{Otto01}, to combine this with the spatial
regularity provided by the square integrability of $\nabla
\tilde\rho_\tau^{1/4}(t,x)$ in $\R^2\times \R$. Our task would be
very much simpler if we had, for almost every $t$, a uniform bound
on $\int_{\R^2}|\nabla \tilde\rho_\tau^{1/4}(t,x)|^2\dd x$, but
don't for any fixed $t$.

The interpolation bound we need is provided by the following
theorem whose proof will be provided after we conclude the proof
of Theorem~\ref{strocom}:

\begin{thm}[Interpolation bound]\label{inerpbnd}
Let $\sigma_0$ and $\sigma_1$ be two densities of mass $M$ in
$\R^2$ such that for some $p>2$, $\|\sigma_0\|_{p+1}^{p+1}\ ,\
\|\sigma_1\|_{p+1}^{p+1} \leq K$, Suppose also that
$\sigma_0^{1/4}$ and  $\sigma_1^{1/4}$ have square integrable
distributional gradients.  Then
\begin{eqnarray}
\|\sigma_0 - \sigma_1\|_2^2 &\leq& \left(  \|\nabla
(\sigma_0)^{1/4}\|_2 + \|\nabla (\sigma_1)^{1/4}\|_2\right)
(2^{5/2} + 2^{9/2}K)
 (\W_2(\sigma_0, \sigma_1))^{(4p-3)/(4p+2)}\nonumber\\
 &+& 16 M^{(p-1)/p}K^{(p+2)/2p} (\W_2(\sigma_0, \sigma_1))^{(p-1)/(2p+1)}\ .\nonumber
 \end{eqnarray}
\end{thm}

We now apply Theorem~\ref{inerpbnd} with
$\tilde\rho_\tau(t,\cdot)$ in place of $\sigma_0$ and
$\tilde\rho_\tau(t-c,\cdot)$ in place of $\sigma_1$.  We have a
uniform bound on $K$ in this case, for any $p<\infty$, and we also
have the bound (\ref{inter55}), so we obtain a finite constant $C$
so that
\begin{multline}
\int_{\R^2}|\tilde\rho_\tau(t,\cdot) - \tilde\rho_\tau(t-c,\cdot)|^2\dd x \leq\\
C\left(\left[  \left(\int_{\R^2}\nabla
\tilde\rho_\tau^{1/4}(t,\cdot)|^2\dd x\right)^{1/2}  +
\left(\int_{\R^2}\nabla \tilde\rho_\tau^{1/4}(t-c,\cdot)|^2\dd
x\right)^{1/2} \right]
 |c|^{\frac{4p-3}{8p+4}} +   |c|^{\frac{p-1}{4p+2}}\right)\ .\nonumber
\end{multline}
Now integrating both sides over $[\tilde a,\tilde b]$, and using
the Cauchy-Schwartz inequality, we obtain
$$
\int_{\R^2\times [\tilde a,\tilde b]}|\tilde\rho_\tau(t,\cdot) -
\tilde\rho_\tau(t-c,\cdot)|^2\dd x \dd t \leq C\left(AT^{1/2}
|c|^{\frac{4p-3}{8p+4}} +   T|c|^{\frac{p-1}{4p+2}}\right)\ ,
$$
where
$$
A :=  \left(\int_{t_0}^T\int_{\R^2}\nabla
\tilde\rho_\tau^{1/4}(t,\cdot)|^2\dd x\dd t\right)^{1/2}  +
\left(\int_{t_0}^T\int_{\R^2}\nabla
\tilde\rho_\tau^{1/4}(t-c,\cdot)|^2\dd x\dd t\right)^{1/2} \ .
$$
Our results so far give us a bound on $A$ that is uniform in $\tau
\leq \tilde \tau^\star$, and thus, choosing $p=2$, the proof of
\eqref{kol3} is complete. \lanbox

\

It remains to prove Theorem~\ref{inerpbnd}. Before beginning the
proof itself, we explain the argument in \cite[Subsection
3.5]{Otto01} that is the basis of the proof. Let $\sigma_0$ and
$\sigma_1$ be two  uniformly bounded densities of mass $M$ in
$\R^2$ .  Let $\dd\Pi$ be the optimal coupling of $\sigma_0$ and
$\sigma_1$.   For $0 < s < 1$, define $\sigma_s$ to be the
displacement interpolant between $\sigma_0$ and $\sigma_1$. That
is, for any bounded continuous function $\varphi$ in $\R^2$,
$$
\int_{\R^2}\varphi(z)\sigma_s(z)\dd z =
\int_{\R^2\times\R^2}\varphi((1-s)x + sy)\dd \Pi(x,y)\ .
$$
As is well known (see for example \cite[Lemma 3]{OLab}), we have
\begin{equation}\label{inter1}
\|\sigma_s\|_\infty \leq \max\{ \ \|\sigma_0\|_\infty\ ,\
\|\sigma_1\|_\infty \ \}\ .
\end{equation}
Suppose furthermore that both $\sigma_0$ and $\sigma_1$ have a
square integrable  gradient. Then we have Otto's interpolation
estimate:
\begin{equation}\label{inter2}
\|\sigma_0 - \sigma_1\|_2  \leq \left( \max\{ \
\|\sigma_0\|_\infty\ ,\ \|\sigma_1\|_\infty \ \}\right)^{1/2}
(\|\nabla \sigma_0\|_2  +  \|\nabla
\sigma_1\|_2)\,W_2(\sigma_0,\sigma_1)\ .
\end{equation}
To see this, note that
\begin{eqnarray}
\|\sigma_0 - \sigma_1\|_2^2 &=& \int_{\R^2\times\R^2} ([\sigma_0(x) - \sigma_1(x)] -
[\sigma_0(y)- \sigma_1(y)])\dd \Pi (x,y)\nonumber\\
&=& \int_0^1  \int_{\R^2\times\R^2}  [\nabla \sigma_1(x+s(y-x))-
\nabla \sigma_0(x+s(y-x))]\cdot (x-y)\dd \Pi (x,y)\dd s\ .
\nonumber
\end{eqnarray}
Now apply the Cauchy-Schwartz inequality. By (\ref{inter1}), we
have that for each $s$,
$$
\int_{\R^2\times\R^2}  |\nabla \sigma_0(x+s(y-x))|^2\dd \Pi (x,y)
= \int_{\R^2}  |\nabla \sigma_0(x)|^2 \sigma_s(x) \dd z
 \leq \left( \max\{ \ \|\sigma_0\|_\infty\ ,\
\|\sigma_1\|_\infty \ \}\right)\|\nabla \sigma_0\|_2^2\ ,
$$
and likewise for $\sigma_1$. Since $\int_{\R^2\times\R^2}|x-y|^2
\dd\Pi(x,y) = \W_2^2(\sigma_0,\sigma_1)$, we obtain
(\ref{inter2}).

We would like to apply this sort of argument with
$\tilde\rho_\tau(t,\cdot)$ in place of $\sigma_0$ and
$\tilde\rho_\tau(t-c,\cdot)$ in place of $\sigma_1$. We cannot do
this directly since $\tilde\rho_\tau$ is not bounded. However, by
the results of the previous section, for $t_0>0$, we have an upper
bound on $\|\tilde\rho_\tau(t,\cdot)\|_p$ for all finite $p$ that
is uniform in $t\geq t_0$.   The next lemma allows us to
approximate $\tilde\rho_\tau(t,\cdot)$ and
$\tilde\rho_\tau(t-c,\cdot)$ by uniformly bounded densities
without significantly increasing the $2$-Wasserstein distance
between them.

\begin{lm}[Approximation by bounded densities]\label{boundaplm}
Let $\sigma$ be any density of mass $M$ on $\R^2$. Suppose that
for some $p > 1$, $\sigma$ belongs to $L^{p+1}(\R^2)$ with
$\|\sigma\|_{p+1}^{p+1} \leq K$. Then for all $\lambda>0$, there
exists density of mass $M$, $\widetilde \sigma$ such that
$\|\widetilde \sigma\|_\infty  \leq 2\lambda$,
\begin{equation}\label{inter5a}
\W_2(\widetilde \sigma,\sigma)   \leq 2K\lambda^{-p-1/2}\ ,
\end{equation}
\begin{equation}\label{inter6}
\|\sigma- \widetilde \sigma\|_2 \leq 2K^{1/2}\lambda^{-(p-1)/2}\ ,
\end{equation}
and finally, $\| \widetilde \sigma\|_2 \leq 2\|\sigma\|_2$.
\end{lm}

\noindent{\bf Proof:}   Let $E_\lambda = \{ \ x\ :\ \sigma(x) >
\lambda\ \}$. Then by Chebychev's inequality,
\begin{equation}\label{inter21}
\int_{E_\lambda} \sigma(x)\dd x \leq
\frac{1}{\lambda^p}\int_{\R^2}\sigma^{p+1}(x)\dd x \leq
\frac{K}{\lambda^p}\ .
\end{equation}
Now for any $h>0$ and any integers $m$ and $n$, define $C_{m,n}$
to be the square
$$
C_{m,n} = \{ x=(x_1,x_2)\in\R^2\ :\ mh \leq x_1 < (m+1)h\ ,\  \ nh
\leq x_2 < (n+1)h\ \}\ .
$$
Define $M_{m,n}$ by
$$
M_{m,n} := \int_{E_\lambda\cap C_{m,n}}\sigma(x)\dd x\ .
$$
Next define
$$
\sigma_{m,n} := 1_{E_\lambda\cap C_{m,n}}\sigma\qquad {\rm
and}\qquad \nu_{m,n} :=  \frac{M_{m,n}}{h^2}1_{C_{m,n}}\ .
$$
There are both densities of mass $M_{m,n}$ supported in $C_{m,n}$.
Let $T_{m,n}$ be any map from $C_{m,n}$  to $C_{m,n}$ such that
$T_{m,n}\#  \sigma_{m,n}  = \nu_{m,n}$. For instance, one may use
the one associated to the optimal coupling. Define $T:\R^2\to
\R^2$ by
$$
T(x) = \begin{cases} \displaystyle\sum_{m,n\in
\Z}1_{C_{m,n}}(x)T_{m,n}(x) & \quad x\in E_\lambda\\ \qquad x &
\quad x\notin E_\lambda\ \end{cases}\ .
$$
Notice that only one term in the sum is non-zero. By construction,
\begin{equation*}
T\#\sigma = (1- 1_{E_\lambda})\sigma +  \sum_{m,n\in
\Z}\frac{M_{m,n}}{h^2}1_{C_{m,n}}\ .
\end{equation*}

By (\ref{inter21}), $M_{m,n} \leq K\lambda^{-p}$.   We now specify
$h := K^{1/2}\lambda^{-(p+1)/2}$. Then $M_{m,n}h^{-2}\leq
\lambda$, and hence with this choice of $h$,   $T\#\sigma$ is
uniformly  bounded above by $2\lambda$. We now define $\widetilde
\sigma :=  T\#\sigma$. It remains to verify that
$\widetilde\sigma$ has all of  the properties claimed in the
lemma.

First of all, notice that $|T(x) - x|^2 = 0$ on the complement of
$E_\lambda$, and on $E_\lambda$ is is bounded by $2h^2$, since
each $T_{m,n}$ maps $C_{m,n}$ into itself. Therefore,
$$
\W_2^2(\sigma,\widetilde\sigma)  \leq \int_{\R^2}|T(x) - x|^2\sigma(x)\dd z \leq 2h^2 \int_{E_\lambda}\sigma(x)\dd x \leq
2K\lambda^{-(p+1)}K\lambda^{-p}\ .
$$
This proves (\ref{inter5a}). Next, by (\ref{inter21}) and Jensen's
inequality,
\begin{equation}\label{inter31a}
\sigma - \widetilde\sigma = 1_{E_\lambda}\sigma  -   \sum_{m,n\in
\Z}\frac{M_{m,n}}{h^2}1_{C_{m,n}}  \qquad \mbox{with} \qquad
\left\Vert \sum_{m,n\in
\Z}\frac{M_{m,n}}{h^2}1_{C_{m,n}}\right\Vert_2  \leq \|
1_{E_\lambda}\sigma\|_2\ .
\end{equation}
Furthermore, by Holder's inequality with the dual indices $p$ and
$p'$,
\begin{eqnarray}\label{inter36}
\int_{\R^2}\sigma^21_{E_\lambda} \dd x =  \int_{\R^2}\sigma^{2-
1/p'}\sigma^{1/p'}1_{E_\lambda} \dd x
&\leq& \left(\int_{\R^2}\sigma^{p+1}\dd x\right)^{1/p}\left(\int_{E_\lambda}\sigma\dd x\right)^{1/p'}\nonumber\\
&\leq& \left(K\right)^{1/p}\left(K\lambda^{-p}\right)^{1/p'}\leq
K\lambda^{1-p}\ .
\end{eqnarray}
This proves (\ref{inter6}). For the final part, note that by
(\ref{inter31a}), $\|\widetilde \sigma\|_2 \leq
\|(1-1_{E_\lambda})\sigma\|_2 + \|1_{E_\lambda}\sigma \|_2  \leq
2\|\sigma\|_2$. This completes the proof of the lemma. \lanbox

\

\noindent{\bf Proof of Theorem~\ref{inerpbnd}:}  For any
$\lambda>0$, let $\widetilde \sigma_0$ and $\widetilde \sigma_1$
be the bounded approximation of $\sigma_0$ and $\sigma_1$ provided
by Lemma~\ref{boundaplm}. Also, define
$$
\underline{\sigma_0}(x)  := \min\{\sigma_1(x)\ ,\, \lambda \
\}\quad{\rm and}\quad \underline{\sigma_1}(x)  :=
\min\{\sigma_1(x)\ ,\, \lambda \ \}\ .
$$
Starting from the identity $\|\sigma_0 - \sigma_1\|_2^2 = \langle
\sigma_0 - \sigma_1,\sigma_0\rangle_{L^2} - \langle  \sigma_0 -
\sigma_1,\sigma_1\rangle_{L^2}$, and adding and subtracting
repeatedly, and using the Cauchy-Schwartz inequality, we obtain
\begin{eqnarray}\label{inter27}
\|\sigma_0 - \sigma_1\|_2^2 &\leq& \langle \underline{\sigma_0} -
\underline{\sigma_1},\widetilde \sigma_0\rangle_{L^2} -
\langle \underline{\sigma_0} - \underline{ \sigma_1},\widetilde\sigma_1\rangle_{L^2}\nonumber\\
&+&(\|\sigma_0\|_2+ \|\sigma_1\|_2)(\|\sigma_0 - \widetilde\sigma_0\|_2 + \|\sigma_1 - \widetilde\sigma_1\|_2)\nonumber\\
&+&(\|\widetilde\sigma_0\|_2+ \|\widetilde\sigma_1\|_2)\left[
(\|\sigma_0 - \underline{\sigma_0}\|_2 + \|\sigma_1
-\underline{\sigma_1}\|_2)\right]\, .
\end{eqnarray}
The heart of the matter is the estimation of $\langle
\underline{\sigma_0} - \underline{\sigma_1},\widetilde
\sigma_0\rangle_{L^2} - \langle \underline{\sigma_0} - \underline{
\sigma_1},\widetilde\sigma_1\rangle_{L^2}$. Let $\widetilde \Pi$
be the optimal coupling of $\widetilde \sigma_0$ and  $\widetilde
\sigma_1$, and for $0< s<1$,  let $\widetilde \sigma_s$ be
displacement interpolant between $\widetilde\sigma_0$ and
$\widetilde\sigma_1$. By what has been explained above,
$\widetilde \sigma_s$ is uniformly bounded above by $2\lambda$.
Then, by the Cauchy-Schwarz inequality once more,
\begin{eqnarray}
\langle \underline{\sigma_0} - \underline{\sigma_1},\widetilde
\sigma_0\rangle_{L^2} - \langle \underline{\sigma_0} - \underline{
\sigma_1},\widetilde\sigma_1\rangle_{L^2} &= &
\int_{\R^2\times\R^2}( [\underline{\sigma_0} - \underline{ \sigma_1}](x) - [\underline{\sigma_0} - \underline{ \sigma_1}](y)) \dd \widetilde \Pi(x,y)\nonumber\\
&=& \int_0^1\int_{\R^2\times\R^2}( [\nabla \underline{\sigma_0} - \nabla \underline{ \sigma_1}]((1-s)y+sx) \cdot(x-y)\dd \widetilde \Pi(x,y)\nonumber\\
&\leq& (2\lambda)^{1/2}\left(  \|\nabla \underline{\sigma_0}\|_2 +
\|\nabla \underline{ \sigma_1}\|_2\right) \W_2(\widetilde\sigma_0,
\widetilde\sigma_1)\, . \nonumber
\end{eqnarray}
Next, for $j=0,1$, let $f_j := (\sigma_j)^{1/4}$. Then  since
$\nabla  \underline{\sigma_j}(x) =0$ unless $f_j(x) \leq
\lambda^{1/4}$, $\|\nabla \underline{\sigma_j}\|_2 \leq
4\lambda^{3/4}\|\nabla f_j\|_2$. Thus,
\begin{equation}\label{inter31}
\langle \underline{\sigma_0} - \underline{\sigma_1},\widetilde
\sigma_0\rangle_{L^2} - \langle \underline{\sigma_0} - \underline{
\sigma_1},\widetilde\sigma_1\rangle_{L^2} \leq
2^{5/2}\lambda^{5/4}\left(  \|\nabla (\sigma_0)^{1/4}\|_2 +
\|\nabla (\sigma_1)^{1/4}\|_2\right) \W_2(\widetilde\sigma_0,
\widetilde\sigma_1) \ .
\end{equation}
By the triangle inequality and Lemma~\ref{boundaplm},
\begin{eqnarray}
\W_2(\widetilde\sigma_0, \widetilde\sigma_1) &\leq&
\W_2(\widetilde\sigma_0, \sigma_0)+
\W_2(\sigma_0, \sigma_1) + \W_2(\sigma_1, \widetilde\sigma_1)\nonumber\\
&\leq& 2K\lambda^{-p-1/2} + \W_2(\sigma_0, \sigma_1) +
2K\lambda^{-p-1/2}\ . \nonumber
\end{eqnarray}
Combining this with (\ref{inter31}), we obtain
\begin{eqnarray}\label{inter32}
\langle \underline{\sigma_0} - \underline{\sigma_1},\widetilde
\sigma_0\rangle_{L^2} - \langle \underline{\sigma_0} - \underline{
\sigma_1},\widetilde\sigma_1\rangle_{L^2} &\leq&
2^{5/2}\lambda^{5/4}\left(  \|\nabla (\sigma_0)^{1/4}\|_2 + \|\nabla (\sigma_1)^{1/4}\|_2\right)   \W_2(\sigma_0, \sigma_1) \nonumber\\
&+& 2^{5/2}\lambda^{5/4}\left(  \|\nabla (\sigma_0)^{1/4}\|_2 +
\|\nabla (\sigma_1)^{1/4}\|_2\right)  4K\lambda^{-p-1/2}\
.\nonumber
\end{eqnarray}

At this point we specify
\begin{equation}\label{inter54}
\lambda :=   (\W_2(\sigma_0, \sigma_1))^{-2/(2p+1)}\ .
\end{equation}
Wit this choice, we have
$$
\langle \underline{\sigma_0} - \underline{\sigma_1},\widetilde
\sigma_0\rangle_{L^2} - \langle \underline{\sigma_0} - \underline{
\sigma_1},\widetilde\sigma_1\rangle_{L^2} \leq \left(  \|\nabla
(\sigma_0)^{1/4}\|_2 + \|\nabla (\sigma_1)^{1/4}\|_2\right)
(2^{5/2} + 2^{9/2}K)
 (\W_2(\sigma_0, \sigma_1))^{(4p-3)/(4p+2)}\ .
$$

Finally, it remains to bound the terms in the last two lines of
(\ref{inter27}). Note that the same  argument using Chebychev's
inequality in (\ref{inter21}) and then H\"older's inequality in
(\ref{inter36}) yields $\|\sigma_j - \underline{\sigma_j}\|_2 \leq
K^{1/2}\lambda^{(1-p)/2}$ for $j= 0,1$. Note also that by
H\"older's inequality once more, $\|\sigma_j\|_2^2 \leq
M^{(p-1)/p}K^{1/p}$. Therefore, the sum of the terms in the last
two lines of (\ref{inter27}) is bounded above by   $16
M^{(p-1)/p}K^{(p+2)/2p}\lambda^{-(p-1)/2}$.  With the value of
$\lambda$ specified in (\ref{inter54}), the contribution of these
terms is
$$16 M^{(p-1)/p}K^{(p+2)/2p} (\W_2(\sigma_0, \sigma_1))^{(p-1)/(2p+1)}\ .$$
Combining results, the proof is complete. \lanbox

\

\begin{lm}[Further subsequence]\label{fisuse}
There is a subsequence of the sequence $\{\tau^{(n)}\}_{n\in\NN}$,
denoted with the same index, such that along this subsequence,
\eqref{stpt10b} is valid for each $0 < a < b < \infty$ and
\begin{equation*}
\lim_{n\to\infty}\tilde \rho_{\tau^{(n)}}(t,x) =   \rho(t,x)
\qquad{\rm almost\  everywhere\ on}\quad (0,T)\times\R^2\
\end{equation*}
for any $T>0$.
\end{lm}

\noindent{\bf Proof:} Let us consider any integer $N\geq N_o$ with
$N_o T>1$. Applying Theorem \ref{strocom} we get a subsequence of
$\{\tau^{(n)}\}_{n\in\NN}$ along which
$\{\rho_{\tau^{(n)}}\}_{n\in\NN}$ converges to $\rho$ strongly in
$L^{2}([1/{N_o},T]\times\R^2)$ and almost surely on
$[1/{N_o},T]\times\R^2$. Next, for $N=N_o+1$, choose a subsequence
of $\{\tau^{(n)}\}$ of the previous subsequence along which
$\{\rho_{\tau^{(n)}}\}_{n\in\NN}$ converges to $\rho$ strongly in
$L^{2}([1/(N_o+1),T]\times\R^2)$ and almost surely on
$[1/(N_o+1),T]\times\R^2$, and so forth. We finish by an obvious
diagonal sequence argument. \lanbox

\

\subsection{Entropy dissipation}

With the strong convergence results obtained in the previous
subsection, we may now establish the entropy-entropy dissipation
inequality:

\begin{thm}[Entropy-entropy dissipation]\label{enendipr}
For each $T>0$ the weak solution $\rho$ of the PKS system that we
have constructed for the initial data $\rho_0$ satisfies
\begin{equation}\label{endis}
\Fhls[\rho(T)]<+\infty\quad\mbox{and}\quad \Fcfd[\rho(T)] +
\int_0^T {\mathcal D}[ \rho(t)] \dd t \leq \Fcfd[\rho_0]\ .
\end{equation}
\end{thm}

\noindent{\bf Proof:} First remind that $\Fhls$ is lower
semi-continuous function for the weak-$L^1$ convergence,
see~\cite[Lemma~3.1]{BCC}. So that the bound in~\eqref{endis} is a
direct consequence of~\eqref{stpt10b}. Also by~\eqref{stpt10b}, it
suffices to show that
\begin{equation*}
\Fcfd[\rho(T)] \le \liminf_{n\to\infty}\Fcfd[\tilde
\rho_{\tau^{(n)}}(T)]\ ,
\end{equation*}
and that for any $0 < a < b < T < \infty$,
\begin{equation}\label{lscen2}
\int_a^b\!\!\int_{\R^2}
\left|\nabla\left(\rho^{1/4}\right)\right|^2 \dd x \dd t \le
\liminf_{n\to\infty} \int_a^b\!\!\int_{\R^2} \left|\nabla
\left(\tilde \rho_{\tau^{(n)}}^{1/4}\right)\right|^2 \dd x \dd t
\, ,
\end{equation}
for a suitable sequence $\{\tau^{(n)}\}_{n\in\NN}$, since the rest
easily follows by a monotone convergence argument for taking $a$
to $0$ and $b$ to $T$.

The first of these follows from the fact that $\Fcfd[\rho]$ is a
lower semi-continuous function on $L^1$ with respect to the $\W_1$
metric just by using the expression of $\Fcfd[\rho]$ in
\eqref{Hlam2}.

To see the second, we again make use of the almost everywhere
convergence proved in the previous subsection. Let $f_n :=
\tilde\rho_{\tau^{(n)}}^{1/4}$ and $f : =\rho^{1/4}$. Then the
sequence of functions $\{f_n\}_{n\in\NN} \to f$ in $L^4\cap L^6
((a,T)\times\RR^2)$ from Lemma~\ref{fisuse}. From \eqref{stpt8},
we have that the sequence $\{\nabla f_n\}_{n\in\NN}$ is bounded in
$L^2 ((a,T)\times\RR^2)$, therefore it has a weakly convergent
subsequence denoted with the same index such that $\{\nabla
f_n\}_{n\in\NN}\rightharpoonup \sigma$ weakly in $L^2
((a,T)\times\RR^2)$. Due to the strong convergence of the sequence
$\{f_n\}_{n\in\NN} \to f$ in $L^4\cap L^6 ((a,T)\times\RR^2)$, it
is simple to identify the weak limit as $\sigma=\nabla f$. By
standard properties of $L^2$-weak convergence, we deduce that
$$
\int_a^b\!\!\int_{\R^2} |\nabla f|^2 \dd x \,\dd t \leq
\liminf_{n\to\infty} \int_a^b\!\!\int_{\R^2} |\nabla f_n|^2 \dd x
\,\dd t
$$
which shows \eqref{lscen2}.\lanbox

\

\noindent{\bf Proof of Theorems~\ref{main}:} As noted above,
Theorem~\ref{enendipr} provides the final step in the construction
of the properly dissipative weak solutions. Theorems~\ref{lpreg}
and \ref{enendipr} provide the additional regularity properties
(1.6.1) and (1.6.2). It remains to prove (1.6.3), the dissipation
of $\Fhls$.

We now show that $\Fhls[\rho(t)]\leq \Fhls[\rho(s)]$ for all
$0\leq s < t$. Take $\tau$ to be any element of the sequence
$\{\tau^{(n)}\}_{n\in\NN}$ whose corresponding approximated
solutions $\{\tilde\rho_{\tau^{(n)}}\}_{n\in\NN}$ converges to the
constructed properly dissipative weak solution $\rho$. Let $j$ be
such that $(j-1)\tau \leq s \leq j\tau$ and let $\ell$ be such
that $\ell\tau \leq t \leq (\ell+1)\tau$. Using
\eqref{monotonefpks} in Lemma~\ref{almdec} and \eqref{regcomp45},
we deduce
\begin{equation*}
\Fhls[\tilde\rho_\tau (t)]\leq
\Fhls^{\epsilon_\ell}[\rho_\tau^{\ell}] \leq
\Fhls^{\epsilon_{\ell-1}}[\rho_\tau^{\ell-1}] + \tilde{Z} C_3
\tau^{2} 2^{-\ell} \leq \Fhls^{\epsilon_{j-1}}[\rho_\tau^{j-1}] +
\tilde{Z} C_3 \tau^{2} \sum_{k=j}^\ell 2^{-k} \leq
\Fhls^{\epsilon_{j-1}}[\rho_\tau^{j-1}] + \tilde{Z} C_3 \tau^{2}\
.
\end{equation*}
Using \eqref{tech99} and the Lemma~\ref{lpreg}, we can control the
error term in the right-hand side by
$$
\Fhls[\tilde\rho_\tau (t)]\leq \Fhls[\tilde\rho_\tau (s)] +
\tilde{Z} \|\tilde\rho_\tau (s)\|_{3/2}^{3/2}\,\tau^{3}2^{-j+1} +
\tilde{Z} C_3 \tau^{2}\leq \Fhls[\tilde\rho_\tau (s)] + 2\tilde{Z}
C_3\tau^{2}\,.
$$
Finally, the {\it a-priori} bounds uniform in $\tau$ due to Lemmas
\ref{contauthm} and \ref{lpreg} together with Lemma~\ref{qualst1}
allow us to pass to the limit $\tau \to 0$ leading to our claim.
\lanbox


\subsection{Large time asymptotics}

We start by identifying the large time asymptotics of the
solutions in a time average sense.

\begin{lm}[Time-averaged strong convergence]\label{tast}
Let $\rho$ be the properly dissipative weak solution of the PKS
system that we have constructed. Then
\begin{equation}\label{longti}
\lim_{T
\to\infty}\left(\int_T^{T+1}\!\!\!\!\int_{\R^2}|\rho(t,x) -
\rhohls(x)|\dd x \dd t\right) = 0\ .
\end{equation}
\end{lm}

\noindent{\bf Proof:} This follows by a standard entropy
dissipation argument. Let $\{t_n\}_{n\in\N}\nearrow+\infty$ be an
increasing diverging sequence of times and consider
$\sigma_n(t,x)=\rho(t+t_n,x)$, for $0\leq t\leq 1$. By using the
entropy dissipation inequality~\eqref{endis} which is true for all
$T>0$, we deduce that
\begin{equation}\label{entrodislim0}
\int_0^\infty \!\!{\mathcal D}[\rho(t)] \dd t <\infty \, ,\qquad
\mbox{and thus},\quad \lim_{n\to\infty}\int_0^1 {\mathcal
D}[\sigma_n(t)] \dd t = 0 \, .
\end{equation}
Now, again using the uniform in time bounds~\eqref{endis} for the
solution, $\Fhls[\sigma_n(t)]<+\infty$, and $\Fcfd[\sigma_n(t)]\le
\Fcfd[\rho_0]\le \crho$, the concentration control inequality in
Lemma \ref{controlcenter} and the Gagliardo-Nirenberg-Sobolev
inequality in Lemma~\ref{gns}, we deduce
\begin{equation}\label{finalcontrollimit}
\frac18 \int_0^1\!\! \int_{\R^2} \sigma_n^{3/2} \dd x \dd t\leq
\int_0^1 \!\!\int_{\R^2} \left|\nabla [\sigma_n^{1/4}]\right|^2\dd
x \dd t \le \frac{\pi}{\gamma_2} \int_0^1 {\mathcal D}[\sigma_n]
\dd t + \frac{C_{{\rm CCD}}}{\gamma_2} \leq \frac{\pi}{\gamma_2}
\crho + \frac{C_{{\rm CCD}}}{\gamma_2}\, .
\end{equation}
Moreover, due to Theorem~\ref{lpreg}, we deduce
\begin{equation}\label{finalcontrollimit2}
\sup_{n\in\NN}\,\sup_{0\leq t\leq 1} \int_{\R^2} \sigma_n^p(t,x)
\dd x \leq C_p\, ,
\end{equation}
for all $1\leq p<\infty$. Note that the sequence
$\{\sigma_n\}_{n\in\NN}$ satisfies the equicontinuity property
\eqref{cont} in Theorem~\ref{contauthm}.

Summarizing, the sequence $\{\sigma_n\}_{n\in\NN}$ has the same
properties \eqref{finalcontrollimit} and
\eqref{finalcontrollimit2} as the sequence of approximate
solutions we used in previous sections to construct the solution
in Theorem~\ref{main}. Proceeding as in Subsections 5.3-5.6, we
deduce the existence of a subsequence, denoted with the same
index, such that $\{\sigma_n\}_{n\in\N}$ converges towards
$\rho_\infty$ with the same convergence properties as in previous
Subsections 5.3-5.6. Here, $\rho_\infty$ is a weak solution of
\eqref{eq:sp} on the time interval $(0,1)$ in the sense of (1.5.2)
 in Definition~\ref{prdisso}. In particular,
$\{\sigma_n\}_{n\in\N}$ converges to $\rho_\infty$ in the metric
space $({\mathcal M}_1^{1/2},\W_1)$, with the notation of
Theorem~\ref{contauthm}, giving
\begin{equation}\label{conv1}
\lim_{n\to\infty} \left[ \max_{0 \leq t \leq 1}\W_1(\sigma_n(t),
\rho_\infty(t))\right] = 0\ .
\end{equation}
Moreover, repeating the arguments in Theorems~\ref{enendipr} and
\ref{fisuse}, we get
\begin{equation}\label{aprioricfdf3}
 \sup_{0\leq t\leq 1}\Fcfd[\rho_\infty(t)]
\le \crho \qquad \mbox{and}\qquad \liminf_{n\to\infty}
\int_0^1\!\!\int_{\R^2} \left|\nabla (\sigma_n)^{1/4}\right|^2 \dd
x \dd t \ge \int_0^1\!\!\int_{\R^2} \left|\nabla
\rho_\infty^{1/4}\right|^2 \dd x\dd t\, ,
\end{equation}
and
$$
 \lim_{n\to\infty}\int_0^1\!\!\int_{\R^2} \sigma_n^{3/2}(t,x)\dd
x \,\dd t = \int_0^1\!\!\int_{\R^2}\rho_\infty^{3/2}(t,x)\dd x
\,\dd t\, .
$$
Furthermore, Theorem~\ref{fisuse} implies the almost everywhere
convergence in $(0,1)\times\RR^2$ of $\{\sigma_n\}_{n\in\N}$
towards $\rho_\infty$, that together with \eqref{conv1} implies
that
\begin{equation}\label{conv2}
\lim_{n\to\infty}\left(
\int_0^1\!\!\int_{\RR^2}|\sigma_n(t,x)-\rho_\infty(t,x)|\,\dd x
\dd t \right) = 0\ .
\end{equation}
Now, let us identify the limit function $\rho_\infty$, passing to
the limit using \eqref{entrodislim0}, we obtain
\begin{equation*}
\int_0^1 \!\!\int_{\R^2} \left(8\left|\nabla \rho_\infty^{1/4}
\right|^2 - \rho_\infty^{3/2} \right)\dd x\dd t =0\;,
\end{equation*}
which means that $\rho_\infty(t)$ is a minimizer to the
Gagliardo-Nirenberg-Sobolev inequality for all $t\in(0,1)$, see
Lemma~\ref{spgn}, and thus that there exists $\bar\lambda(t)$ such
that $\rho_\infty(t) = \varrho_{\bar\lambda(t)}$ where $\rhohls$
is the family of the minimizers of the Gagliardo-Nirenberg-Sobolev
inequality, see Lemma~\ref{spgn}. Due to \eqref{aprioricfdf3} then
$\Fcfd[\varrho_{\bar\lambda(t)}] < \infty$, we conclude that
$\bar\lambda(t)=\lambda$ since $\Fcfd[\varrho_{\mu}] = +\infty$
for $\mu\neq\lambda$. Therefore, $\rho_\infty(t) =
\varrho_{\lambda}$ that together with \eqref{conv2} implies
\eqref{longti}. \lanbox

\

We now will take advantage of the other Lyapunov functional, we
shall prove that $\lim_{t\to\infty}\Fhls[\rho(t)] =
\Fhls[\rhohls]$. In doing this, we shall make essential use of the
monotonicity of $\Fhls[\rho(t)]$. Let us introduce for any $C>0$
and $\delta>0$ the set
$$
{\mathcal S}_{C,\delta}:=\left\{ \rho\in L^1_+(\RR^2) \, :\,
\int_{\RR^2} \rho(x) \dd x =8\pi\, , \, \int_{\R^2}|x|^\delta
\rho(x)\dd x \leq C \, \mbox{ and }
\int_{\R^2}\rho^{1+\delta}(x)\dd x \leq C \right\}\,.
$$
\begin{lm}[Convergence  for $\Fhls$]\label{qualst1} Given any
sequence $\{\rho_n\}_{n\in\NN}$ in ${\mathcal S}_{C,\delta}$ there
is a $\rho\in {\mathcal S}_{C,\delta}$ and a subsequence
$\{\rho_{n_k}\}_{k\in\NN}$ such that
$$
\lim_{k\to\infty}\|\rho_{n_k} - \rho\|_1 =0 \qquad{\rm and}\qquad \lim_{k\to\infty}\Fhls[\rho_{n_k}] = \Fhls[\rho]\ .
$$
\end{lm}
\noindent{\bf Proof:} Choose any $0 <\delta' < \min(\delta,1)$ so
that $2\delta' /(1-\delta') \leq \delta$. By uniform integrability
arguments such as we have made above, see Step 1 of
Theorem~\ref{strocom}, we can find a subsequence (denoted with the
same index) along which $\{\rho_n\}_{n\in\NN}$ is weakly
convergent in $L^{1+\delta'}(\R^2)$ and along which
$\{\rho_n^{1+\delta'}\}_{n\in\NN}$ is weakly convergent in $L^1$.
It follows as in Subsection 5.6 that $\{\rho_n\}_{n\in\NN}$ is
strongly convergent in $L^{1+\delta'}(\R^2)$, and passing to a
further subsequence, we may suppose it is also almost everywhere
convergent, and strongly convergent in $L^1(\R^2)$. Let $\rho$
denote the limit. By Fatou's Lemma, $\rho\in {\mathcal
S}_{C,\delta}$.

Since for $t\geq 1$, $t\log t \leq (t^{ \delta} -1)/\delta'$, we
have  for $\rho \geq 1$, $\rho\log \rho \leq
(1/\delta')\rho^{1+\delta'}$ and for $\rho < 1$, $\rho\log(1/\rho)
\leq (1/\delta')\rho^{1-\delta'}$.  Since for $\epsilon =
2\delta'$, we get
\begin{align*}
\int_{\R^2}\rho^{1-\delta'} (1+|x|^2)^\epsilon \dd x &= \int_{\R^2}\rho^{1-\delta'}(1+|x|^2)^{2\epsilon} (1+|x|^2)^{-\epsilon}\dd x\nonumber\\
&\leq \left(\int_{\R^2}\rho  (1+|x|^2)^{2\epsilon/(1-\delta')}\dd
x\right)^{1-\delta'} \left(\int_{\R^2} (1+|x|^2)^{-2}\dd
x\right)^{\delta'}\ ,
\end{align*}
our choice of $\delta'$ gives the uniform integrability of
$\{\rho_{n}^{1-\delta'}\}_{n\in\NN}$. Then, by what we have said
above, $\{\rho_{n}\log \rho_{n}\}_{n\in\NN}$ is uniformly
integrable, and hence
$$
\lim_{n\to\infty}\int_{\R^2}\rho_n \log \rho_n \dd x = \int_{\R^2}\rho \log \rho \dd x\ .
$$
The convergence of the positive part of the interaction potential
is straightforward, due to the uniform bound of $\{|x|^\delta
\rho_n\}_{n\in\NN}$ in $L^1(\RR^2)$ and a dominated convergence
argument. Concerning the negative part, it follows by Young's
inequality for convolutions using the convergence of
$\{\rho_n\}_{n\in\NN}\to\rho$ in $L^{1+\delta'}(\R^2)$ and the
fact that $\log_-|x| \in L^p (\RR^2)$ for all $1\leq p<\infty$.
\lanbox

\

\begin{lm}[Qualitative stability for $\Fhls$]\label{qualst2}
For any $\epsilon> 0$ and $C>0$, there exists
$\delta(\epsilon,C)>0$ so that if $\rho\in {\mathcal
S}_{C,\delta}$, then
$$
\Fhls[\rho] \leq 8\pi(-1+\log(8\pi)) + \delta(\epsilon,C)\quad
\Rightarrow\quad
  \|\rho - \varrho_\mu\|_1\leq \epsilon\qquad{\rm for\ some}\ \mu > 0\ ,
$$
and for any $R>0$, there exists $\delta(\epsilon,C,R)>0$ such that
$$
\Fhls[\rho] \leq 8\pi(-1+\log(8\pi)) + \delta(\epsilon,C,R)\quad
\Rightarrow\quad
  \left( \int_{\{|x|\leq R\}}  |(\sqrt{\rho} - \sqrt{\varrho_\mu})(x)|^2
  \dd x\right)^{1/2} \leq \epsilon\qquad{\rm for\ some}\ \mu > 0\ .
$$
\end{lm}

\noindent{\bf Proof:} Given $C,R>0$ fixed, suppose not. Then for
some $\epsilon>0$, there is a sequence $\{\rho_n\}_{n\in\NN}$ in
${\mathcal S}_{C,\delta}$ such that
$\lim_{n\to\infty}\Fhls[\rho_n] = 8\pi(-1+\log(8\pi))$ but
$$
\inf_{n,\mu} \|\rho_n - \varrho_\mu\|_1 \geq \epsilon\ .
$$
However, by Lemma~\ref{qualst1}, there is a subsequence, still
indexed by $n$, converging strongly in $L^1(\RR^2)$ to $\rho\in
{\mathcal S}_{C,\delta}$, such that
$$
8\pi(-1+\log(8\pi)) = \lim_{n\to\infty}\Fhls[\rho_n] =
\Fhls[\rho]\ .
$$
By the cases of equality in the Log-HLS inequality, $\rho =
\rho_\mu$ for some $\mu$. This is a contradiction. The second part
is proved the same way, using the uniform integrability of the
$\sqrt{\rho}$, $\rho\in {\mathcal S}_{C,\delta}$ on $\{|x|\leq
R\}$. \lanbox

\

\noindent{\bf Proof of Theorems~\ref{basinat}:} Recall that
$\Fhls[\rho(t)]\leq \Fhls[\rho(s)]$ for all $0\leq s
< t$.
We now apply this monotonicity to  improve  our large time asymptotic result.

By
\eqref{longti} in Lemma~\ref{tast}, there is a sequence of times
$\{t_n\}_{n\in\NN}\nearrow\infty$ such that $
\lim_{n\to\infty}\|\rho(t_n) - \rhohls\|_1 = 0$. By our regularity
results in Lemmas \ref{contauthm} and \ref{lpreg},
$\{\rho(t_n)\}_{n\in\NN} \subset {\mathcal S}_{C,\delta}$ for some
$0< C,\delta< \infty$. Then by Lemma~\ref{qualst1}, there is a
subsequence, still indexed by $n$, such that
$\lim_{n\to\infty}\Fhls[\rho(t_n)] = \Fhls[\rhohls]$. By the
monotonicity of $\Fhls[\rho(t)]$ it follows that
$$
 \lim_{t\to\infty}\Fhls[\rho(t)] = \Fhls[\rhohls] = 8\pi(-1 + \log(8\pi))\ .
$$
Then by Lemma~\ref{qualst2} it follows that given $R>0$ there
exists $\mu>0$ such that for all sufficiently large $t$,
$$
\|\rho(t) - \varrho_\mu\|_1 \leq \epsilon \qquad \mbox{and} \qquad
 \, \|\sqrt{\varrho_\mu} -\sqrt \rho \|_{2,R} :=\left( \int_{\{|x|\leq R\}}  |(\sqrt{\rho} - \sqrt{\varrho_\mu})(x)|^2
  \dd x\right)^{1/2} \leq \epsilon.
$$
However, for any $R>0$, by Minkowskii's inequality and
\eqref{endis},
\begin{align*}
\left(\int_{\{|x|\leq R\}}|\sqrt{\varrho_\mu} -\sqrt{\rhohls}|^2
\rhohls^{-1/2}\dd x \right)^{1/2} &\leq
\left(\int_{\{|x|\leq R\}}|\sqrt{\varrho_\mu} -\sqrt \rho |^2 \rhohls^{-1/2}\dd x \right)^{1/2}  + \sqrt{\Fcfd[\rho]}\\
&\leq \sqrt{8\lambda}(\lambda+R^2)\|\sqrt{\varrho_\mu} -\sqrt \rho
\|_{2,R} + \sqrt{\crho}\ .
\end{align*}
Since the left hand side diverges as $R$ increases, uniformly for
$|\mu - \lambda| > \delta>0$, we readily conclude that
$\mu=\lambda$ and $\lim_{t\to\infty}\|\rho(t) - \rhohls\|_1 = 0$.
\lanbox


\section*{Appendix: proof of Lemma~\ref{displconvexitydelta}}
By Lemma~\ref{discobnd} applied to $u_0=\rho$ and $u_1=\rho_0$
\begin{equation*}
   \Fcfdr[\rho] \le \Fcfdr[\rho_0] - \int_{\R^2}\left[  2\,x\,f_\delta'(|x|^2)
   +\frac{1}{2}\frac{\nabla \rho}{(\rho+\delta)^{3/2}}\right]\cdot  (\nabla\varphi(x) -x)\,\rho\dd x  - K_\delta(\rho,\rho_0)\
\end{equation*}
with $f_\delta'(s)=8\lambda\left[8\lambda  + \delta(\lambda +
s)^{2}\right]^{-3/2}$. Using~\eqref{disgrad}, {\it i.e.
$x-\nabla\varphi(x) = \tau \left(\nabla c_\epsilon - \nabla
\rho/\rho \right)$} and expanding we can rewrite it as
\begin{equation*}
\Fcfdr[\rho]  \le  \Fcfdr[\rho_0] + \tau \left( \frac12{\rm (I)} +
\frac12\,{\rm (II)} +2\,{\rm (III)} + 2\,{\rm (IV)} \right) -
K_\delta(\rho,\rho_0)\ ,
\end{equation*}
where
\begin{align*}
&{\rm (I)}:=-\int_{\R^2} \frac{\left|\nabla \rho\right|^2}{(\rho +
\delta)^{3/2}}\dd x,\qquad \qquad\,\,\,\,\,\,\,\,\,\,\,\,\,{\rm
(II)}:=\int_{\R^2}\frac{\rho\nabla
c_{\epsilon} \cdot\nabla \rho}{(\rho + \delta)^{3/2}}\dd x \\
&{\rm (III)}:=\int_{\R^2} f_\delta'(|x|^2)\,x \cdot\nabla
c_\epsilon\,\rho\dd x\,\quad\mbox{and}\quad {\rm
(IV)}:=-\int_{\R^2} f_\delta'(|x|^2)\,x \cdot  \nabla \rho\dd x\ ,
\end{align*}

We will keep the term ${\rm (I)}$ and we need to perform some
integration-by-parts in the other terms:

\noindent{\bf Control of ${\rm (II)}$:} We can rewrite this term
as
$${\displaystyle
  {\rm (II)}= 2\,\int_{\R^2} \nabla\left(\frac{\rho+2\delta}{\sqrt{\rho+\delta}}\right)\cdot \nabla c_{\epsilon}\dd
  x}\,.
$$
Integrating by parts on the ball of radius $R$ and noticing
that $-\Delta c_\epsilon=\gamma_\epsilon*\rho*\gamma_\epsilon$, we
obtain
\begin{align}
  2\,\int_{|x|\le R}\!\!\! \nabla\left(\frac{\rho+2\delta}{\sqrt{\rho+\delta}}\right)\cdot \nabla c_{\epsilon}\dd
  x= \,& 2\int_{|x|\le R} \left(\frac{\rho+2\delta}{\sqrt{\rho+\delta}}\right)\;(\gamma_\epsilon*\rho*\gamma_\epsilon)\dd x\nonumber\\
&+ 2\int_{|x|= R} \frac{\rho+2\delta}{\sqrt{\rho+\delta}}\nabla c_{\epsilon}\cdot n\dd \sigma\nonumber\\
\le \, & 2\int_{|x|\le R}\!\!
\sqrt{\rho}\;(\gamma_\epsilon*\rho*\gamma_\epsilon)\dd x +
4\,\sqrt{\delta}\int_{|x|\le R}\!\!\gamma_\epsilon*\rho*\gamma_\epsilon\dd x\nonumber\\
&+ 2\int_{|x|= R}
\frac{\rho+2\delta}{\sqrt{\rho+\delta}}\;\left|\nabla
c_{\epsilon}\right|\dd \sigma\nonumber\\ \le \, & 2\int_{|x|\le R}
\sqrt{\rho}\;(\gamma_\epsilon*\rho*\gamma_\epsilon)\dd x +
32\,\pi\,\sqrt{\delta} \nonumber\\
&+ 2\int_{|x|= R} (\sqrt{\rho}+2\sqrt{\delta})\;\left|\nabla
c_{\epsilon}\right|\dd \sigma \ ,\label{tech0}
\end{align}
where we used twice the estimate
\begin{equation*}
  \frac{\rho+2\delta}{\sqrt{\rho+\delta}} \le \sqrt{\rho} + 2\sqrt{\delta}\ .
\end{equation*}

Let us deal first with the second boundary term. By
\eqref{regconc2} in Lemma \ref{regconc}, we have
\begin{align}\label{tech1}
4\sqrt{\delta}\int_{|x|= R}| \nabla c_{\epsilon}| \dd \sigma & \le
8\pi\sqrt{\delta}\left(8C_{{\rm HLS}}\,\|\gamma\|_{4/3}
\||x|\gamma\|_{4/3} + 4 +\,\frac{C_{{\rm
HLS}}}{2\pi\epsilon}\,\|\gamma\|_{4/3}^2\,\|
|x|\rho\|_{1}\right)\\
& := 16\pi\sqrt{\delta} (J_\gamma + 1 + \hat C_\epsilon) \
,\nonumber
\end{align}
for all $R>0$. To cope with the first boundary term, we observe
that taking any natural $N>1$, and considering
\begin{align*}
\left(\int_{N-1}^N\int_{|x|= r}\sqrt{\rho}\,| \nabla c_{\epsilon}|
\dd \sigma\dd r\right)^2 &\le \pi\left(N^2 - (N-1)^2\right) \,\|\nabla c_{\epsilon}\|_\infty^2 \int_{N-1\le |x|\le N}\rho \dd x \\
&\le \pi\,\| \nabla c_{\epsilon}\|_\infty^2 \,\frac{2\,N-1}{N-1}\int_{N-1\le |x|\le N}|x|\,\rho\dd x
 \le \widetilde C_\epsilon^2\, \int_{N-1\le |x|\le N}|x|\,\rho\dd x
\end{align*}
where $\widetilde C_\epsilon:={\sqrt{2\,\pi}}\,\| \nabla
c_{\epsilon}\|_\infty$. Since
${\displaystyle
\sum_{N=3}^\infty \int_{N-1}^N\int_{|x|= r}\sqrt{\rho}\,| \nabla
c_{\epsilon}| \dd \sigma\dd r  \le
\widetilde C_\epsilon\,\left(\int_{\R^2}|x|\,\rho\dd x \right)^{1/2}}$,
it follows that
${\displaystyle
\lim_{N\to\infty}\int_{N-1}^N\int_{|x|= r}\sqrt{\rho}\,| \nabla
c_{\epsilon}| \dd \sigma\dd r  =0},$
and thus, there exists a sequence $\{R_j\}\nearrow \infty$ such
that
\begin{equation}\label{tech2}
\lim_{j\to\infty} \int_{|x|= R_j}\sqrt{\rho}\,| \nabla
c_{\epsilon}| \dd \sigma = 0.
\end{equation}
Plugging \eqref{tech1} and \eqref{tech2} into
\eqref{tech0}, we get
\begin{equation*}
  {\rm (II)} \le 2\int_{\R^2} \sqrt{\rho}\;(\gamma_\epsilon*\rho*\gamma_\epsilon)\dd x +
  32\,\pi\,\sqrt{\delta}+4\sqrt{\delta}( J_\gamma + 2\pi)\ .
\end{equation*}
Finally, a simple application of H\"older's inequality gives
$$
\int_{\R^2}\sqrt{\rho}\;(\gamma_\epsilon*\rho*\gamma_\epsilon)\dd
x \leq \|\sqrt{\rho}\|_3\,
\|\gamma_\epsilon*\rho*\gamma_\epsilon\|_{3/2}\leq \int_{\R^2}
\rho^{3/2}\dd x\ ,
$$
to conclude
\begin{equation}\label{tech6}
  {\rm (II)} \le 2\int_{\R^2} \rho^{3/2}\dd x +
  32\,\pi\,\sqrt{\delta}+4\sqrt{\delta}( J_\gamma + 2\pi)\ .
\end{equation}

\noindent{\bf Control of ${\rm (III)}$:} Remind that
$f_\delta'(s)\leq (8\lambda)^{-1/2}:=\kappa_\lambda/2$ and that
$2f_\delta'(s)\to \kappa_\lambda$ as $\delta\to 0$, see Proposition~\ref{entap}. By definition of $c_\epsilon$ and $G_\epsilon$ and by
symmetry of $\gamma$
\begin{align*}
 {\rm (III)} &= \int_{\R^2}\rho(x)\,f_\delta'(x) \,x\cdot \left(\nabla G_\epsilon *\rho\right)(x)\dd x =\iint_{\R^2\times\R^2}\rho(x)\,f_\delta'(x)\,x\,\gamma_\epsilon(x-z)\,\left(\nabla G*\gamma_\epsilon*\rho\right)(z)\dd z \dd x\\
& =\int_{\R^2}(\rho\,f_\delta'\,{\rm id}*\gamma_\epsilon)(z)\,\left(\nabla G*\gamma_\epsilon*\rho\right)(z)\dd z\ .
\end{align*}
By definition of the convolution, we have
\begin{align*}
 \rho\,f_\delta'\,{\rm id}*\gamma_\epsilon(z) &=\int_{\R^2}\rho(z-x)\,f_\delta'(z-x)\,(z-x)\,\gamma_\epsilon(x)\dd x\nonumber\\
&=z\,\int_{\R^2}\rho(z-x)\,f_\delta'(z-x)\,\gamma_\epsilon(x)\dd x- \int_{\R^2}\rho(z-x)\,f_\delta'(z-x)\,x\,\gamma_\epsilon(x)\dd x\nonumber\\
&=z\,(\rho f_\delta'*\gamma_\epsilon)(z)-(\rho f_\delta'*{\rm id}\,\gamma_\epsilon)(z)\;.
\end{align*}
As a consequence, we infer
\begin{equation}\label{key3b}
{\rm (III)} = \int_{\R^2} \left[z\,(\rho
f_\delta'*\gamma_\epsilon)(z)-(\rho f_\delta'*{\rm
id}\,\gamma_\epsilon)(z) \right]\cdot \left(\nabla
G*\gamma_\epsilon*\rho\right)(z)\dd z:={\rm (III)}_1-{\rm
(III)}_2\;.
\end{equation}
By the symmetrization argument just as in~\eqref{key2}, the first term of the right hand side of~\eqref{key3b} reads
\begin{align*}
{\rm (III)}_1 &=\frac12 \int_{\R^2} z\,(\rho
(2f_\delta'-\kappa_\lambda)*\gamma_\epsilon)(z)\cdot \left(\nabla
G*\gamma_\epsilon*\rho\right)(z)\dd z-8\pi\kappa_\lambda\\ &:={\rm
(III)}_{11}-8\pi\kappa_\lambda\, .
\end{align*}
We now control ${\rm (III)}_{11}$  using the
HLS and Young inequalities, see~\eqref{spclhls} to obtain
\begin{align*}
|{\rm (III)}_{11}| &\leq \frac12 \iint_{\R^2\times\R^2}
|z|\,|(\rho (2f_\delta'-\kappa_\lambda)*\gamma_\epsilon)(z)|
\frac{1}{|x-z|}|(\gamma_\epsilon*\rho)(x)|\dd z\dd x\\
&\leq \frac{C_{{\rm HLS}}}{4\pi} \,\||z|\,|(\rho
|2f_\delta'-\kappa_\lambda|*\gamma_\epsilon)\|_{4/3}\,
\|\rho*\gamma_\epsilon\|_{4/3} \\
&\leq \frac{C_{{\rm HLS}}}{4\pi} \,\left(\|(|x|\,\rho
|2f_\delta'-\kappa_\lambda|)*\gamma_\epsilon\|_{4/3} + \|(\rho
|2f_\delta'-\kappa_\lambda|)*(|x|\,\gamma_\epsilon)\|_{4/3}\right)\,
\|\rho*\gamma_\epsilon\|_{4/3}\\
&\leq 2C_{{\rm HLS}} \,\left(\||x|\,\rho
|2f_\delta'-\kappa_\lambda|\|_1 \, \|\gamma_\epsilon\|_{4/3} +
\|\rho |2f_\delta'-\kappa_\lambda|\|_1 \,
\||x|\,\gamma_\epsilon\|_{4/3}\right)\,
\|\gamma_\epsilon\|_{4/3}\, ,
\end{align*}
from which
\begin{equation}\label{eq:tocollectIII}
{\rm (III)}_1 \leq -8\pi\kappa_\lambda + C_\epsilon \,
\|\,|2f_\delta'-\kappa_\lambda|\,(1+|x|)\rho \|_1\, .
\end{equation}

To estimate the second term of the right hand side
of~\eqref{key3b}, we make again use of the
HLS inequality, see~\eqref{spclhls}:
\begin{align*}
   \left|{\rm (III)}_2\right| &\le
 \frac{1}{4\,\sqrt{2\lambda}\,\pi}\iint_{\R^2\times\R^2} (\rho*|z|\,\gamma_\epsilon)(z)\,\frac{1}{|z-y|}\, (\gamma_{\epsilon}*\rho)(y)\dd y \dd z\\
 &\le \frac{C_{{\rm HLS}}}{4\,\sqrt{2\lambda}\,\pi} \, \|\rho*(|z|\,\gamma_\epsilon)\|_{4/3}\,\|\gamma_{\epsilon}*\rho\|_{4/3}\
.
\end{align*}
By the Young inequality, and a direct calculation,
${\displaystyle
 \|\rho*(|z|\,\gamma_\epsilon)\|_{4/3} \le 8\pi \,\||x|\,\gamma_\epsilon\|_{4/3} = 8\pi \,\sqrt\epsilon\,\||x|\,\gamma\|_{4/3}}$.
and in the same way
${\displaystyle
\| \gamma_{\epsilon}*\rho\|_{4/3} \le 8\pi \,\|
\gamma_{\epsilon}\|_{4/3} = 8\pi \,\frac{1}{\sqrt\epsilon}\,\|
\gamma\|_{4/3}}.$
The positive and negative powers of $\epsilon$ cancel, and using
\eqref{eq:tocollectIII}, we conclude
\begin{align}\label{tech3}
{\rm (III)} &\leq -8\pi\kappa_\lambda + C_\epsilon \,
\|\,|2f_\delta'-\kappa_\lambda|\,(1+|x|)\rho \|_1 +
\frac{16\pi\,C_{{\rm HLS}}}{\sqrt{2\lambda}} \,
\||x|\,\gamma\|_{4/3}\,\|\gamma\|_{4/3}\, \nonumber\\
&= -8\pi\kappa_\lambda + C_\epsilon \,
\|\,|2f_\delta'-\kappa_\lambda|\,(1+|x|)\rho \|_1 +
\frac{8\pi\,J_\gamma}{\sqrt{2\lambda}}.
\end{align}

Let us estimate this third term in a different way that will be
useful later on. Using again the Young inequality, but this time
eliminating $\gamma$ instead of $\rho$, \emph{ i.e.}
$\|\gamma_{\epsilon}*\rho\|_{4/3} \le  \| \rho \|_{4/3}$, we get
$$
\|\rho*(|z|\,\gamma_\epsilon)\|_{4/3}\,\|\gamma_{\epsilon}*\rho\|_{4/3}\le
64\pi^2 \sqrt\epsilon\, \||x|\,\gamma\|_{4/3} \,\| \rho \|_{4/3} \
.
$$
As a consequence, we get this other control on ${\rm (III)}$ by
\begin{equation}\label{tech4}
{\rm (III)} \leq -8\pi\kappa_\lambda + C_\epsilon \,
\|\,|2f_\delta'-\kappa_\lambda|\,(1+|x|)\rho \|_1 +
\frac{16\pi\,C_{{\rm HLS}}}{\sqrt{2\lambda}} \, \sqrt\epsilon\,
\||x|\,\gamma\|_{4/3} \,\| \rho \|_{4/3}\ .
\end{equation}

\noindent{\bf Control of ${\rm (IV)}$:} By integrating by parts
for any $R>0$, we have
$$
\int_{|x|\le R} x\,f_\delta'(|x|^2) \cdot \nabla \rho(x)\dd x =
\int_{|x|= R} \rho(x)\,f_\delta'(|x|^2)\,x \cdot n \dd \sigma -
\int_{|x|\le R}\nabla\cdot[x\,f_\delta'(|x|^2)] \,\rho(x)\dd x\ ,
$$
where $n$ denotes the outward normal to the disk $D_R$. Taking
into account that
\begin{equation*}
  \sum_{N=3}^\infty \int_{N-1}^N\int_{|x|= r}|x|\,|f_\delta'(|x|^2)|\,{\rho(x)} \dd \sigma\dd r
  \le \frac{1}{2\sqrt{2\lambda}} \sum_{N=3}^\infty \int_{N-1\le |x|\le N}|x|\,\rho\dd x
  \le \frac{1}{2\sqrt{2\lambda}} \int_{\R^2}|x|\,\rho\dd x <\infty
\end{equation*}
we have
${\displaystyle
  \lim_{N\to\infty} \int_{N-1}^N\int_{|x|= r} \rho(x)\,|f_\delta'(|x|^2)|\,|x| \dd \sigma\dd r =0}$.
As a consequence, there exists a sequence $\{R_j\}\nearrow \infty$
such that
${\displaystyle
\lim_{j\to\infty} \int_{|x|= R_j} \rho(x)\,|f_\delta'(|x|^2)|\,|x|
\dd \sigma = 0}$,
and thus, we conclude
\begin{equation} \label{tech5}
{\rm (IV)}= \int_{\R^2}\nabla\cdot[x\,f_\delta'(|x|^2)] \rho(x)\dd
x\ .
\end{equation}

The desired estimates are obtained by putting together estimates
\eqref{tech6}, \eqref{tech3}, \eqref{tech4} and \eqref{tech5}.
\lanbox


\subsection*{Acknowledgments}
All authors thank CRM-Barcelona where this work was started during
the visit of the first and second authors in 2007. We also thank
I.P.A.M. where all three authors were guests at various times in
2008 and 2009, and finally the University of Toulouse Paul
Sabatier where the second and third authors were guests in 2008,
2009 and 2010, and where work was continued and concluded.
Finally, we thank the referee for helpful comments.


\end{document}